%% file: BoxBoseGas_ArXiv.tex
                  \def\version{13 January, 2022}		       %

\documentclass[reqno,12pt]{amsart} 
\usepackage[utf8]{inputenc}
\usepackage{amsmath} 
\usepackage[mathscr]{eucal}
\usepackage{amssymb}
\usepackage{srcltx} 
\usepackage{dsfont}
\usepackage[pagebackref, colorlinks=true,linkcolor=blue,citecolor=blue]{hyperref}
\usepackage{color} 
\usepackage{tikz}
\usepackage{caption,subcaption}
\usetikzlibrary{decorations.pathreplacing}
\usetikzlibrary{shapes.geometric,decorations,decorations.pathmorphing}
\usetikzlibrary{fit}
\pgfmathsetseed{1}

\numberwithin{equation}{section}
 
\def\emptyset{\varnothing}

\newfam\Bbbfam 
\font\tenBbb=msbm10 
\font\sevenBbb=msbm7 
\font\fiveBbb=msbm5 
\textfont\Bbbfam=\tenBbb 
\scriptfont\Bbbfam=\sevenBbb 
\scriptscriptfont\Bbbfam=\fiveBbb

\newcommand{\R}     {\mathbb{R}} 
\newcommand{\Z}     {\mathbb{Z}} 
\newcommand{\N}     {\mathbb{N}} 
\renewcommand{\P}   {\mathbb{P}} 
 
\newcommand{\E}     {\mathbb{E}}

\newcommand{\ct}  {{\mathfrak c}}
\def\L{\Lambda} 
\newcommand{\smfrac}[2]{\textstyle{\frac {#1}{#2}}}

\def\1{{\mathchoice {1\mskip-4mu\mathrm l}      
{1\mskip-4mu\mathrm l} 
{1\mskip-4.5mu\mathrm l} {1\mskip-5mu\mathrm l}}} 
\newcommand{\ssup}[1] {{{\scriptscriptstyle{({#1}})}}} 
\def\comment#1{} 
\newtheoremstyle{thm}{2ex}{2ex}{\itshape\rmfamily}{} 
{\bfseries\rmfamily}{}{1.7ex}{} 
 
\newtheoremstyle{rem}{1.3ex}{1.3ex}{\rmfamily}{} 
{\itshape\rmfamily}{}{1.5ex}{}

\newtheorem{theorem}{Theorem}[section] 
\newtheorem{lemma}[theorem]{Lemma} 
\newtheorem{prop}[theorem] {Proposition} 
\newtheorem{cor}[theorem]  {Corollary}

\newtheorem{step}{STEP} 
 
\theoremstyle{definition}

\renewcommand{\v}{{\bar v}} 
\renewcommand{\d}{{{\rm d}}} 
 
\newcommand{\eps}{\varepsilon} 
\newcommand{\Lambdaeb}{{\rm Leb}}

\newcommand{\supp}{{\operatorname {supp}}}

\newcommand{\diam}{{\operatorname {diam}}} 
\newcommand{\Poi}{{\operatorname {Poi}}}

\newcommand{\Acal}  {{\mathcal A}}
\newcommand{\Bcal}  {{\mathcal B}}
\newcommand{\Ccal}   {{\mathcal C }}

\newcommand{\Hcal}   {{\mathcal H }}

\newcommand{\Mcal}   {{\mathcal M }}

\newcommand{\Rcal}   {{\mathcal R }}

\newcommand{\Ucal}   {{\mathcal U }}

\newcommand{\e}   {{\operatorname e }}
\newcommand{\one}   {{\mathds{1}}}

\definecolor{Red}{rgb}{1,0,0}

\begin{document} 
 
\title[Box-version of the Bose gas]{\large The free energy of a box-version\\ \medskip of the interacting Bose gas}

\author[Orph\'ee Collin, Benedikt Jahnel and Wolfgang  K{\"o}nig]{} 

\maketitle

\thispagestyle{empty} 
\vspace{0.2cm}

\centerline {\sc By Orph\'ee Collin\renewcommand{\thefootnote}{1}\footnote{DMA, \'Ecole Normale Sup\'erieure, Universit\'e PSL, CNRS, 75005 Paris, France, {\tt orphee.collin@normalesup.org}}, Benedikt Jahnel\renewcommand{\thefootnote}{2}\footnote{Weierstrass Institute for Applied Analysis and Stochastics,
Mohrenstr.\ 39, 10117 Berlin, Germany, {\tt jahnel@wias-berlin.de}} and  Wolfgang
K{\"o}nig\renewcommand{\thefootnote}{3}\footnote{Technische Universit\"at Berlin, Str. des 17. Juni 136,
10623 Berlin, and Weierstrass Institute for Applied Analysis and Stochastics,
Mohrenstr.\ 39, 10117 Berlin, Germany, {\tt koenig@wias-berlin.de}}}
\medskip

\vspace{0.4cm}

\centerline{\small(\version)} 
\vspace{.5cm}

\begin{abstract}The interacting quantum Bose gas is a random ensemble of many Brownian bridges (cycles) of various lengths with interactions between any pair of legs of the cycles. It is one of the standard  mathematical models in which a proof for the famous Bose--Einstein condensation phase transition is sought for. We introduce a simplified version of the model with an organisation of the particles in deterministic boxes instead of Brownian cycles as the marks of a reference Poisson point process (for simplicity, in $\Z^d$ instead of $\R^d$). 

We derive an explicit and interpretable variational formula in the thermodynamic limit for the limiting free energy of the canonical ensemble for any value of the particle density. This formula features all relevant physical quantities of the model, like the microscopic and the macroscopic particle densities, together with  their mutual and self-energies and their entropies. The proof method comprises a two-step large-deviation approach for marked Poisson point processes and an explicit distinction into small and large marks.

In the characteristic formula, each of the microscopic particles and the statistics of the macroscopic part of the configuration are seen explicitly; the latter receives the interpretation of the condensate.  The formula enables us to prove a number of properties of the limiting free energy as a function of the particle density, like differentiability and explicit upper and lower bounds, and a qualitative picture below and above the critical threshold (if it is finite). This proves a modified saturation nature of the phase  transition. However, we have not yet succeeded in proving the existence of this phase transition.
\end{abstract}

\bigskip



\tableofcontents

\setcounter{section}{0}

\section{Bosonic systems via point processes: our main purpose}\label{Intro}

\noindent One of the greatest unsolved problems in mathematical physics is a mathematical understanding of the famous {\em Bose--Einstein condensation (BEC)} phase transition in  the {\em interacting quantum Bose gas} in the thermodynamic limit at sufficiently low, but positive temperature. This is a large symmetrised system of $N$ particles in a box of volume $\asymp N$; each particle is equipped with a kinetic energy, and the system is subject to a pair energy. 

Feynman \cite{F53} explained and interpreted the Bose gas in terms of a large interacting {\em ensemble of Brownian cycles} of various lengths, in each of which a random number of particles is spatially organised. The {\em condensate} is interpreted as the part of the particles that lie in very long cycles, i.e., in cycles of lengths that diverge with the particle number, $N$. The {\em condensation phase} is characterised by the appearance of a {\em macroscopic} part of the particles in these long  cycles, i.e., a number of particles that is $\asymp N$. We call each cycle with a fixed length {\em microscopic} and the ensemble of these cycles the microscopic part of the system. The prominent BEC conjecture (initiated by Bose's and Einstein's seminal papers in 1924 and 1925) is that, in dimensions $d\geq 3$, but not in $d\leq 2$, for sufficiently low temperature (or equivalently, for sufficiently high particle density), this macroscopic structure indeed would  emerge. Furthermore, this emergence is predicted to be triggered by a saturation effect.

A rigorous mathematical formulation of Feynman's picture, which goes back to Ginibre's work in the 1960s (see \cite{G70}) is in terms of a {\em Feynman--Kac formula}, a random interacting ensemble of many Brownian cycles of various lengths. This model is sometimes called the {\em interacting quantum Bose gas}. A further reformulation step was made in \cite{ACK10} and describes the system in terms of a {\em random Poisson point process with marks}, the marks being the cycles that start and terminate at the Poisson points, see the summary in Section \ref{sec-Bosegas}.

In \cite{ACK10}, this formulation was taken as the base of the following strategy to prove the occurrence of BEC in this model:
\begin{enumerate} 
\item Use the theory of marked random point processes to rewrite the partition function in terms of the {\em empirical stationary field},

\item adapt and apply large-deviation theory for the ergodic behaviour of random point processes to find an explicit formula for the limiting free energy,

\item reformulate within the frame of that formula what BEC means,

\item find a criterion under which the occurrence of BEC can be proved.
\end{enumerate}

However, this programme could not be completed in \cite{ACK10}. While a characteristic formula could be derived as an upper bound for the limiting free energy, a lower bound could be derived only for sufficiently small particle densities and was in terms of a slightly different formula, and a correct interpretation, not to mention a proof for, BEC could not be attained. The formula is not able to describe a macroscopic structure, and it was nevertheless unclear if it would be able to yield a formula for the limiting free energy in general. Furthermore, it seemed out of reach to prove this or to identify a clear criterion for the emergence of BEC from that formula.

In this paper, we  make decisive progress with regard to these open questions in a slightly simplified model that we introduce here and call a {\em box version of the interacting quantum Bose gas}. The main difference to the interacting Bose gas is that the marks are not taken as random Brownian cycles, but as deterministic centred boxes. This strongly simplifies the complexity of the underlying probability space, but keeps the most important characteristics of the model: an interacting Poisson point process with marks of unbounded sizes that have a spatial extent. For technical reasons, we decided to use a $\Z^d$-model rather than an $\R^d$-model. For this model we here derive a characteristic variational formula for its limiting free energy that is significantly extended and reveals much more and much more explicit information. We firmly believe that the progress that we make here for this model will enable us in future work to make an analogous progress for the interacting Bose gas.

One main novelty of our ansatz in this paper is a certain extension of the frame of the formula, such that both the non-condensate part (the {\it microscopic} structure) and the condensate part (the {\it macroscopic} structure) are explicitly seen in the formula, together with all interaction. This is even a great advantage over the current state of description of the free (i.e., non-interacting) Bose gas, where the phase transition is revealed just by detecting a loss of mass in the finite cycles. Instead, we create an enlargement of the description space in which also the macroscopic part of the gas (if it exists) is characterised. Indeed, it creates a global environment for the microscopic marks  labelled by $a\in\N$, which is the number of copies of the grid $\Z^d$ that is locally created by macroscopic boxes. This global environment is fully characterised by the percentage $\psi(a)$ of the space that is locally covered by precisely $a$ copies of $\Z^d$. 

This approach has a number of advantages. First, we are able to identify the free energy for {\em any} value of the particle density, not only in the non-condensate phase. Second, this formula possesses always a minimiser, since it has much better continuity and compactness properties. By exploring these advantages, we succeed in deriving a number of interesting properties of the limiting free energy, like differentiability and asymptotics for large and for small particle densities. Furthermore, all objects appearing in the formula admit clear interpretations and give in particular a clear criterion for the occurrence of the micro-macro phase transition that is analogous to BEC. However, we do not prove the present paper that this phase transition indeed occurs. This is devoted to future work; it seems to require the application of much finer methods to the variational formula than we are capable of yet. Nevertheless, if we  assume that the transition occurs, then we can prove a very detailed picture; in particular it surprisingly turns out that it is not a phase transition of {\em saturation type}, but quite close.

The organisation of the paper is as follows. In Section \ref{sec-simplified} we introduce our model, the box-version of the interacting Bose gas, and state and discuss the main results of this paper. In Section \ref{sec-phasetrans+Lit}, we explain why the proof of the existence of the phase transition is difficult, and we give a small literature survey. In Section \ref{sec-proofgrid} we derive the variational formula describing the limiting free energy, in Section \ref{sec-AnaForm} we prove our results on some analytical properties of the variational formula, and in Section \ref{sec-ProofBEC} we prove differentiability of the free energy with respect to the particle density and derive explicit formulas. In the Appendix, Section \ref{sec-Bosegas}, we recall the work of \cite{ACK10} on the interacting quantum Bose gas for comparison.

\section{The box version of the Bose gas, and main results}\label{sec-simplified}

\noindent In this section, we introduce a simplified model of the well-known interacting Bose gas and identify its free energy in terms of  a characteristic variational formula with explicit control on the microscopic and the macroscopic components. We introduce the model in Section \ref{sec-model}, formulate our identification of its free energy in Section \ref{sec-results} and some results on existence of minimisers in Section \ref{sec-Discussion}; then we discuss the nature of the phase transition in Section \ref{sec-BEC}, providing it exists.

\subsection{The box version of the Bose gas}\label{sec-model} 

The model that we are going to introduce has the following characteristics:
\begin{itemize}
\item It is defined as a marked Poisson point process in the $d$-dimensional Euclidean space.

\item Each mark is a particle configuration of $k$ particles for some $k\in\N$, centred at the Poisson point; the density $q_k$ of the size-$k$ marks is summable on $k\in\N$.

\item Any two particles in the system underly a pair interaction with an arbitrary nonnegative interaction functional having compact support.

\item We look at the {\em thermodynamic limit}, i.e., we have in total precisely $N$ particles in a box of volume $\asymp N$.

\end{itemize}

In these respects, the model is of the same type as the interacting Bose gas, which we review in Section \ref{sec-Bosegas}. However, the following feature makes the model different, and basically {\em only} this:

\begin{itemize}
\item The marks are deterministic boxes instead of random Brownian cycles.

\end{itemize}

We decided to work in the $\Z^d$-setting rather than in the $\R^d$-setting, which we consider a minor difference. We keep the model simple in order to concentrate on our main goal, the derivation of an interpretable variational formula for the free energy, and not to overburden the derivation with technicalities.

\begin{figure}[!htpb]
	\centering
\includegraphics[scale=0.485]{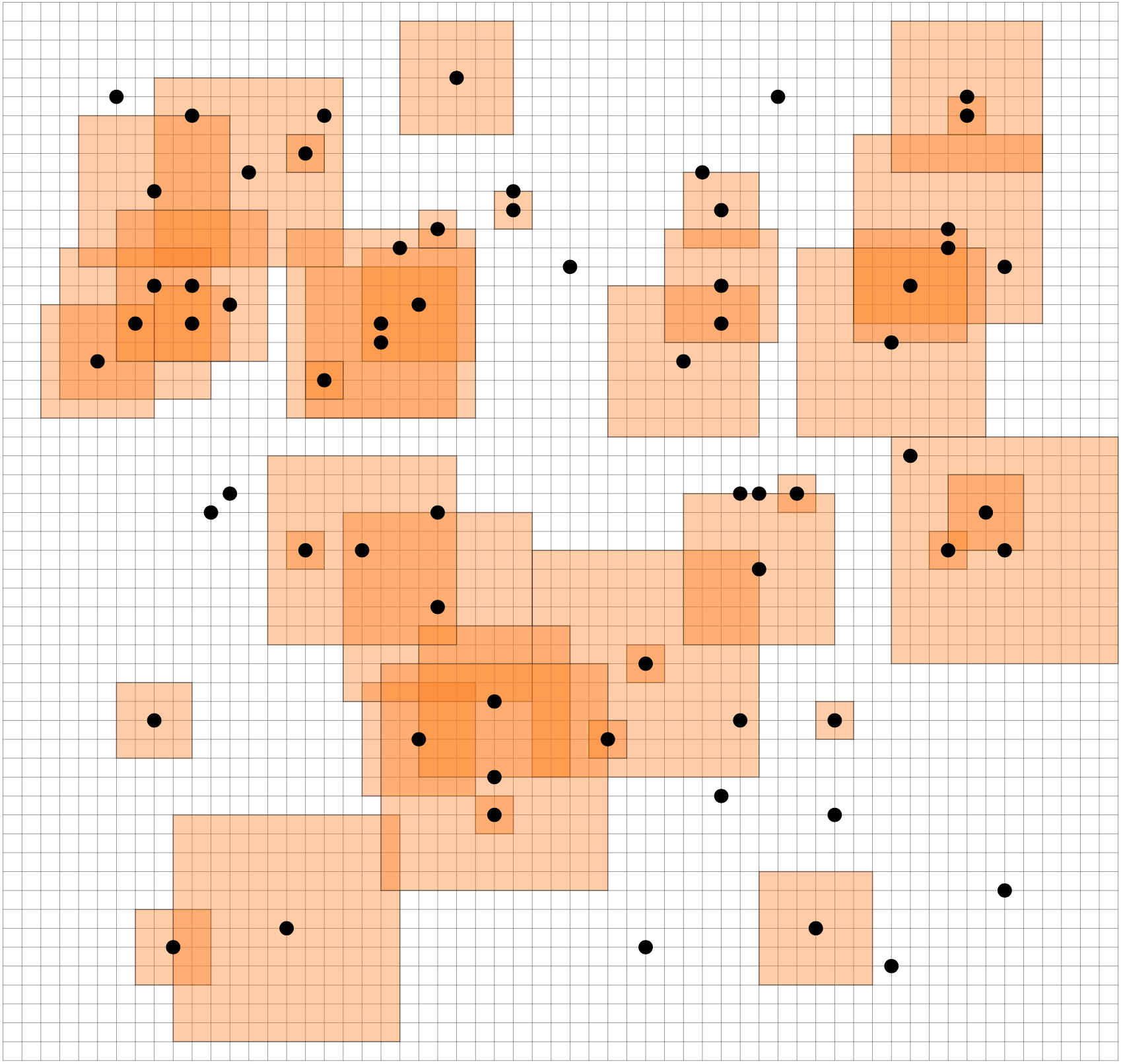}
	\caption{A realisation of the box version of the Bose gas in a finite container with Dirichlet boundary conditions.}
	\label{Pix_GridBoseGas}
\end{figure}

We consider configurations consisting of {\em points} in $\Z^d$ with {\em marks} that are subsets of $\Z^d$. For any $k\in\N$ write $\xi^{\ssup k}(x)\in\N_0$ for the number of points of a configuration with parameter $k$ at site $x\in\Z^d$. To each such point we attach a copy of a {\em mark} $G_k$, which is a deterministic subset of $\Z^d$ (approaching a large box for large $k$) with 
$$
\big([-L_k,L_k]^d\cap \Z^d\big)\subset G_k\subset \big([-L_k-1,L_k+1]^d\cap \Z^d\big)\qquad\text{ and }\qquad|G_k|=k,
$$
for some $L_k\in\N_0$. Hence, the configuration is uniquely determined by the collection $(\xi^{\ssup k})_{k\in\N}$ with $\xi^{\ssup k}=(\xi^{\ssup k}(x))_{x\in\Z^d}$. For any fixed $k\in\N$, this gives rise to a process of {\em marked points}
$$
\omega^{\ssup k}=\sum_{x\in\Z^d}\xi^{\ssup k}(x)\delta_{(x,G_k)}. 
$$
We call the elements of $x+G_k$ the {\em particles} of the site $x$ and note that several particles and several points may be at the same site. Then, we consider the {\em particle configuration}
$$\omega=\sum_{k\in\N}\omega^{\ssup k}=\sum_{k\in\N}\sum_{x\in\Z^d}\xi^{\ssup k}(x)\delta_{(x,G_k)}$$
as superpositions of the configurations of marked points with fixed $k$. See Figure~\ref{Pix_GridBoseGas} for an illustration. We write $\Omega$ for the set of all such configurations and  equip it with the usual evaluation sigma algebra.

It will be convenient to use the following notation. For any $\L,\L'\subset\Z^d$, we denote by 
$$M_{\L,\L'}^{\ssup{\delta_k}}(\omega)=\sum_{x\in \L}\sum_{y\in \L'}\xi^{\ssup k}(x)\one\{y\in x+G_k\}$$
the number of particles in $\L'$ that are attached to points in $\L$ via marks of size $k$. Using this, we denote by 
$$N_\L^{\ssup{\delta_k}}(\omega)=M_{\L,\Z^d}^{\ssup{\delta_k}}(\omega)$$
the number of points in $\L$ with the mark $G_k$, and we write $N_\Lambda^{\ssup f}(\omega)=\sum_{k\in\N}f(k)N_\L^{\ssup{\delta_k}}(\omega)$ for $f\colon\N\to [0,\infty$. Specifically, we denote by $N_\Lambda^{\ssup \ell}(\omega)=\sum_{k\in\N} k N_\L^{\ssup{\delta_k}}(\omega)$ the number of particles belonging to points in $\Lambda$, that is, we write $\ell(k)=k$ for the identity map. On the other hand, we write 
$$\widetilde N_\L^{\ssup{\delta_k}}(\omega)=M_{\Z^d,\L}^{\ssup{\delta_k}}(\omega)$$
for the number of particles in $\L$ that come from marks of size $k$, and abbreviate $\widetilde N_\L(\omega)=\sum_{k\in \N}\widetilde N_\L^{\ssup{\delta_k}}(\omega)$. In general, we write $M_{x,y}, N_x$ and $\widetilde N_x$ instead of $M_{\{x\},\{y\}}, N_{\{x\}}$ and $\widetilde N_{\{x\}}$ respectively.

We introduce a particle-to-particle interaction, including also all self-interactions, namely
\begin{equation}\label{energysimple}
\Phi_{\Lambda,\Lambda'}(\omega)= \sum_{x\in\Lambda, y\in\Lambda'}\sum_{k,l\in\N}\xi^{\ssup k}(x)\xi^{\ssup{l} }(y) T_{x,y}(G_k,G_{l}),\qquad \Lambda,\Lambda'\subset\Z^d,
\end{equation}
where
\begin{equation}\label{Tfunctiondef}
T_{x,y}(G,G')=\sum_{i\in G}\sum_{j\in G'}v(x+i-y-j),\qquad G,G'\subset\Z^d.
\end{equation}
Here $v\colon \Z^d \to[0,\infty)$ is some function with compact support. We assume that $v$ is symmetric in the sense that $v(x)=v(-x)$ for any $x$. If $\Lambda=\{z\}$ is a singleton, then we write $\Phi_{z,\Lambda'}$ instead of $\Phi_{\{z\},\Lambda'}$, analogously with $\Lambda'$. Note that we consider the particles of different marks of the process $\omega$ as different, even though they might be on the same spot. We also put
\begin{equation}\label{Cvdef}
{\v}=\sum_{i\in \Z^d} v(i),
\end{equation}
which is the interaction of a particle at the origin with the deterministic homogeneous grid $\Z^d$.

We now introduce a reference probability measure on the set $\Omega$ of point configurations.  More precisely, we denote by $\P^{\ssup k}$ the Poisson point process (PPP) in $\Z^d$ with intensity measures $q_k \ct$, where $\ct$ is the counting measure on $\Z^d$. Then, any point configuration $\xi^{\ssup k}$ can be seen as a realisation of an i.i.d.~field of Poisson random variables with parameter $q_k$. Moreover, we denote by $\P$ the independent superposition of the PPPs $\P^{\ssup k}$ and assume that $(q_k)_{k\in\N}$ is a summable sequence of positive numbers. 
Using this, any particle configuration $\omega$ can be seen as realisation of a PPP on $\Z^d\times\{G_k\colon k\in\N\}$ with intensity measure $\sum_{k\in\N}(q_k\ct\otimes \delta_{G_k})$.

For a finite set $\Lambda\subset\Z^d$, we write $\Omega_\L$ for the set of restrictions of configurations of $\omega\in\Omega$ to $\L$, i.e., the image of $\Omega$ under the projection $\omega\mapsto \omega_{\L}=\sum_{x\in\L}\sum_{k\in\N}\xi^{\ssup k}(x)\delta_{(x,G_k)}$. The image measure of $\P$ under this projection is denoted by $\P_\L$. We consider zero Dirichlet boundary conditions in $\L$. It is denoted \lq Dir\rq\ and means that all particles of the marks are contained in $\Lambda$. We denote the corresponding probability measure and partition function by 
$$
\P^{\ssup{\rm{ Dir}}}_{\L}(\cdot)=\P_\L(\cdot\,|\, M_{\L,\L^{\rm c}}=0)\qquad\mbox{and}\qquad Z_{N,\L}^{\ssup{\rm Dir}}=\E_\L^{\ssup{\rm Dir}}\Big[\e^{-\Phi_{\Lambda,\Lambda} }\1\{N_{\Lambda}^{\ssup\ell}=N\}\Big].
$$
We are sure that periodic boundary condition can be used well and will lead to the same results, but we abstain from including this in our analysis, to avoid a further blow up of the paper. However, we believe that open boundary conditions (where the points are restricted to $\L$, but the particles may project beyond $\L$) will leave to a different behaviour.

Note that we do not introduce any temperature parameter in this model. 
Also observe that relaxing the assumption of $v$ being symmetric would not make the model more general, since the model remains unchanged when $v$ is replaced by its symmetrised version $i\mapsto \frac 12(v(i)+v(-i))$.

As we announced, this model is analogous to the well-known interacting Bose gas at positive temperature with deterministic boxes instead of Brownian cycles, see Section \ref{sec-Bosegas} and in particular Proposition \ref{lem-rewrite}.

\subsection{Main result: a variational formula for the free energy}\label{sec-results}

We are going to formulate our main result about the limiting free energy of this model in the thermodynamic limit: a description in terms of a variational formula, valid for any particle density $\rho\in(0,\infty)$. We denote by $\Mcal_1^{\ssup{\rm s}}(\Omega)$ the set of all shift-invariant probability measures on $\Omega$, where we recall that $\Omega$ is the set of configurations of points on $\Z^d$ carrying marks in $\{G_k\colon k\in\N\}$. By 
\begin{equation}\label{Igriddef}
I(P)=\lim_{Q\uparrow\Z^d}\frac 1{|Q|}H(P_Q|\P_Q),\qquad P\in\Mcal_1^{\ssup{\rm s}}(\Omega),
\end{equation} 
we denote the entropy density function with respect to the reference  distribution $\P$, where $P_Q$ is the projection of $P$ from $\Z^d$ to  $Q$ (more precisely, from $\Omega$ to $\Omega_Q$), and the limit is w.r.t.\ diverging radius of  centred boxes $Q$. By $H(\mu|\nu)$ we denote the relative entropy of a finite measure $\mu$ with respect to another one, $\nu$, on a discrete space $\mathcal X$, defined by
\begin{equation}\label{Entropydef}
H(\mu|\nu)= \sum_{x\in\mathcal X}\Big[\nu(x)-\mu(x)+ \mu(x)\log \frac{\mu(x)}{\nu(x)}\Big],
\end{equation}
if $\mu\ll\nu$, and otherwise $H(\mu|\nu)=\infty$. According to~\cite{G88,GZ93}, the limit  in \eqref{Igriddef}  exists, and $I$ is an affine and lower-semi-continuous function with compact level sets $\{P\colon I(P)\leq \alpha\}$ for any $\alpha\in\R$ in the topology of local tame convergence, the topology on $\Mcal_1^{\ssup{\rm s}}(\Omega)$ that is induced by test integrals against local functions $f\colon \Omega\to\R$ that are bounded as $|f(\omega)|\leq C(1+N_\L^{\ssup{\1}}(\omega))$ for some finite $\L\subset\Z^d$ and some $C>0$, for any $\omega\in\Omega$. It is an easy exercise to show that the maps $P\mapsto P(N_0^{\ssup\ell})$ and $P\mapsto P(\Phi_{0,\Z^d})$ are lower semi-continuous in this topology. We write $\mu(f)=\int f\,\d \mu=\sum_x f(x)\,\mu(x)=\langle \mu, f\rangle$ for the integral of a integrable function $f$ with respect to a measure $\mu$ on a discrete space.  

Here is our main result. 

\begin{theorem}\label{thm-freeenergygrid}
Assume that the $k$-box densities of the reference PPP satisfy $q_k=\e^{o(k)}$ as $k\to\infty$. Fix $\rho\in(0,\infty)$ and a symmetric  interaction functional $v\colon\Z^d\to[0,\infty)$ having a compact support. Then, for the centred boxes $\Lambda_N$ with volumes $N/\rho$, 
\begin{equation}\label{freeenergyvarform}
\lim_{N\to\infty}\frac1{|\Lambda_N|}\log Z_{N,\L_N}^{\ssup{\rm Dir}}=-\inf_{ \rho_{\rm mi}, \rho_{\rm ma}\geq 0\colon  \rho_{\rm mi}+ \rho_{\rm ma}=\rho}\chi( \rho_{\rm mi}, \rho_{\rm ma}),
\end{equation}
where 
\begin{eqnarray}
\chi( \rho_{\rm mi}, \rho_{\rm ma})&=&\inf\Big\{\varphi(m,\psi)\colon m\in [0,\infty)^\N,\psi\in\Mcal_1(\N_0),\sum_{k\in\N}km_k= \rho_{\rm mi}, \sum_{a\in\N_0}a\psi(a)= \rho_{\rm ma}\Big\},\qquad\label{chi(m,rho2)def}\\
\varphi(m,\psi)&=&\inf\Big\{\sum_{a\in\N_0}\psi(a)\Big[I(P_a)+P_a(\Phi_{0,\Z^d})+2{\v}aP_a(N_0^{\ssup \ell})+{\v}a^2 \Big]\colon\notag \\
\qquad\qquad & &\qquad P_0,P_1,P_2,\dots\in  \Mcal_1^{\ssup{\rm s}}(\Omega), \sum_{a\in\N_0}\psi(a)P_a(N_0^{\ssup {\delta_k}})=m_k\,\, \forall k\in\N\Big\}.\label{varphidef}
\end{eqnarray}
\end{theorem}

The proof is presented in Section~\ref{sec-proofgrid}. Even though the most important object  here is the empirical stationary field of the reference Poisson point process and $I$ is the large-deviation rate function for this, we are not using this large-deviation principle (which is well-known from \cite{G93, GZ93, G94}), but we  go via another route. Instead, we decompose the box $\L_N$ regularly  into mesoscopic boxes, neglect all interaction between them and apply a large-deviation principle in the spirit of Sanov's theorem. Afterwards, we let the mesoscopic box approach $\Z^d$ and use the spatial ergodic theorem and compactness arguments. This method can be seen as an alternate route for deriving the LDP by Georgii/Zessin. Its application is even necessary here because of the disordered appearance of macroscopic marks.

\begin{figure}[!htpb]
	\centering
\includegraphics[scale=0.57]{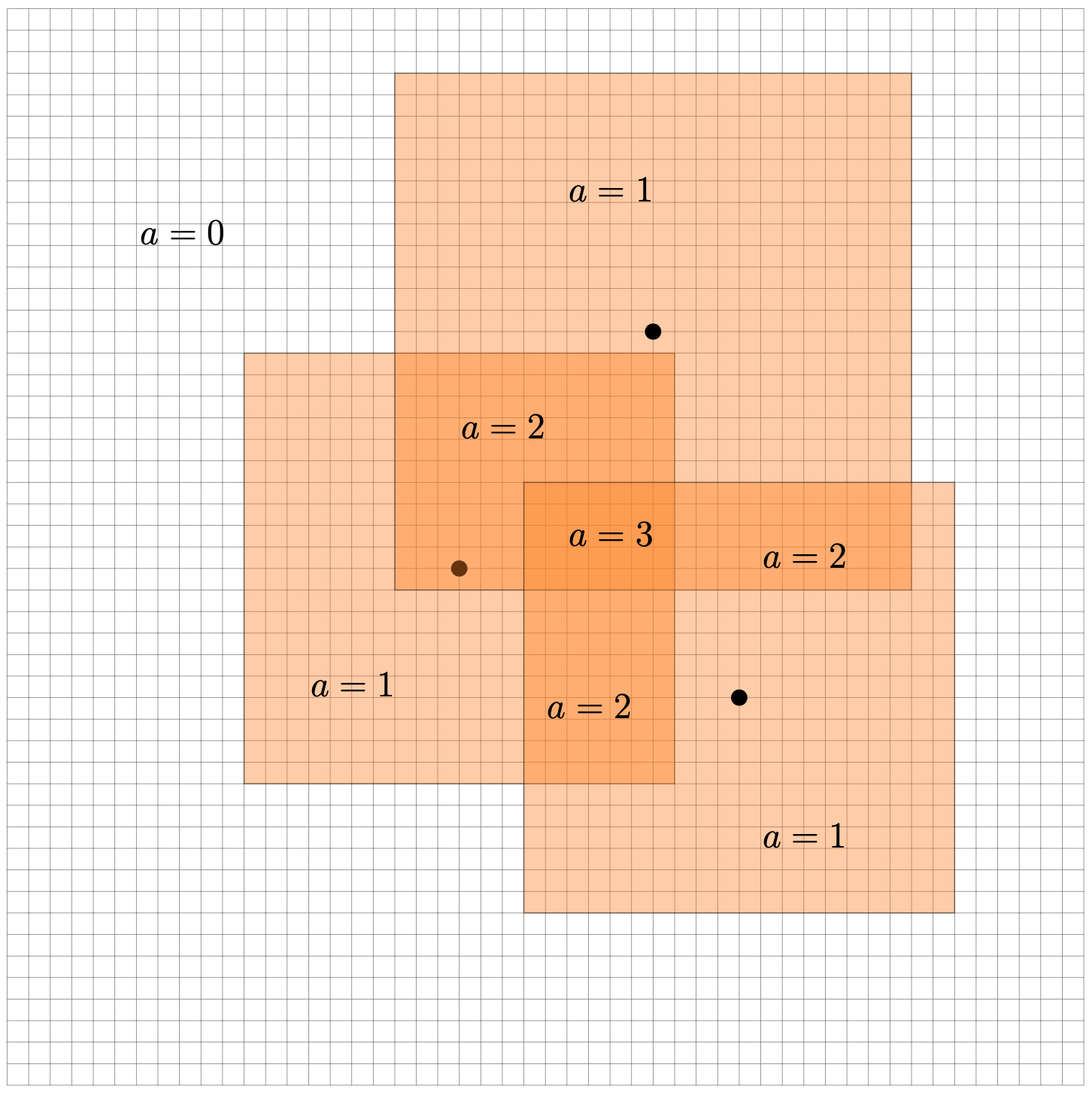}
	\caption{An environment of three macroscopic boxes, creating regions of various overlap numbers.}
	\label{Pix-macroconf}
\end{figure}

Let us give now a non-technical interpretation of Theorem \ref{thm-freeenergygrid}, see Figure \ref{Pix-macroconf} for an illustration. It is important to note that any of the objects appearing in the characteristic formula on the right-hand sides of \eqref{chi(m,rho2)def}--\eqref{varphidef} contains information about the particle ensemble, even though nothing of this is explicitly formulated nor proved. Making exact statements would require a two-step limiting procedure and involve auxiliary parameters.

The most prominent quantities are the {\em microscopic particle density}, $ \rho_{\rm mi}$, and the {\em macroscopic particle density}, $ \rho_{\rm ma}$, of the configuration, that is, the number of particles in microscopic, i.e., finite-size, marks $G_k$, $k\in\N$, per unit volume, and the number of particles in macroscopic marks per unit volume, i.e., in marks of sizes that depend on $N$ and diverge as $N\to\infty$. The $m_k$'s give a more precise information about the microscopic part; it is the spatial rate of size-$k$ marks. Hence, $\rho_{\rm mi}=\sum_k k m_k$ is the microscopic particle density. Likewise, $\psi$ gives a more precise information about the macroscopic marks: for any $a\in\N_0$, the number $\psi(a)$ is the percentage of the volume of the box $\L_N$ in which precisely $a$ of the macroscopic marks overlap each other. Hence $\rho_{\rm ma}=\sum_a a\psi(a)$ is the macroscopic particle density. 

In analogy of the famous {\em Bose--Einstein condensate} phase transition in the analogous interacting Bose gas, we sometimes refer to the macroscopic marks as to {\em condensate marks} and to their union as to the condensate. This is entirely for analogy reasons and does not imply any assertion about the physics of the model.

One should see (at least our proof suggests that) $\L_N$ as a regular decomposition into many large boxes in which the condensate density is constant, and $\psi(a)$ is the percentage of the number of those boxes in which precisely $a$ condensate marks overlap. In such a box, the spatial distribution of all the microscopic points is given by the stationary marked point process distribution $P_a$. Interestingly, the microscopic particles are randomly distributed, while the macroscopic ones are deterministic; it is an $a$-fold superposition of $\Z^d$. The latter builds a {\em condensate environment}, in which the microscopic part of the configuration floats; see Figure~\ref{Pix-macroconf}. The last condition in \eqref{varphidef} says that, averaged over all condensate environments, the density of $k$-sized marks is equal to $m_k$ everywhere.

The first term in the first line of \eqref{varphidef} is the entropy of the spatial distribution of all the microscopic points with respect to the reference measure, the second is the internal energy of all their microscopic marks, the third term is the interaction between the microscopic particles and the condensate environment (the $a$ marks), and the  fourth and last term is the internal energy of all the macroscopic marks. More precisely, the three energy terms are the interaction between the origin and the respective remainder.

Then  \eqref{freeenergyvarform} says that the main contribution to the partition function comes from those particle configurations that are represented by the minimising objects, provided they exist. Therefore, it will be of high importance to get clear information about the existence or non-existence of minimisers.

Our highest interest is in the question under what conditions a condensate occurs, i.e., the question about the existence of a phase transition of condensation type. Condensation occurs if a minimising configuration $(m,\psi)$ with $\psi\not=\delta_0$ exists or if even any minimising configuration has this property. More about that in Section \ref{sec-BEC}.

Before we enter questions about minimisers, let us give a number of properties of $I$, $\varphi$ and $\chi$ that can be easily deduced from their defining formulas. For $m\in[0,\infty)^{\N}$, we denote by $\P^m$ the process defined as the reference process $\P$ with $q$ replaced by $m$, by $I_m$ the entropy density function with respect to the process $\P^m$, defined as in \eqref{Igriddef} with $\P^m$ instead of $\P$, and we recall that $H(m|q)=\sum_{k\in \N} (q_k-m_k+m_k\log\frac{m_k}{q_k})$. By $\chi^{\ssup{v=0}}(\rho,0)$ we denote the value for the non-interacting model. Let us note that $\chi(0,0)$ is the free energy for the process restricted to having no marks in the box $\L_N$, with $\L_N\uparrow \Z^d$, extending Theorem~\ref{thm-freeenergygrid} to the case where $\rho=0$ in a natural manner. We write $\leq_{\rm st}$ for stochastic ordering on $\Mcal_1(\N_0)$, i.e., $\psi\leq_{\rm st} \psi'$ holds if and only if $\psi([a,\infty))\leq \psi'([a,\infty))$ for all $a\in\N_0$.

\begin{lemma}[Simple properties of $I$, $\varphi$ and $\chi$]\label{lem-properties}We have the following.
\begin{enumerate}
\item The function $\varphi$ is convex jointly in $m$ and $\psi$.
\item The function $\chi$ is convex and continuous jointly in $(\rho_{\rm mi}, \rho_{\rm ma})\in [0,\infty)^2$. In particular, $\rho\mapsto\chi(\rho,0)$ is convex.
\item For any $m$, the function $\varphi(m,\cdot)$ is non-decreasing in $\psi$ with respect to the stochastic ordering.

\item The function $\chi$ is non-decreasing in $\rho_{\rm ma}$.
\item For any $m\in [0,\infty)^{\N}$ and any $P\in\Mcal_1^{\ssup{\rm s}}(\Omega)$ satisfying $ P(N_0^{\ssup {\delta_k}})=m_k$ for all $k\in\N$, we have 
\begin{equation}
I(P)=H(m|q)+I_m(P).
\end{equation} 
In particular, $I(\P^m)=H(m|q)$ for any $m\in[0,\infty)^{\N}$.
\item We have $\chi(0,0)=\chi^{\ssup{v=0}}(0,0)=\sum_k q_k$ and  $\partial_{ \rho_{\rm mi}}\chi(0,0)=-\infty$. 
\item We have the bounds
\begin{equation}\label{chibounds}
\chi^{\ssup{v=0}}(\rho,0)+ \v \rho^2\leq \chi(\rho,0)\leq \chi^{\ssup{v=0}}(\rho,0)+\v (\rho^2+\rho),\qquad \rho\in[0,\infty),
\end{equation}
under the assumption $2v(0)\geq \v$ for the first inequality. 
\end{enumerate}
\end{lemma}

Lemma \ref{lem-properties} is used in Section \ref{sec-proofgrid}, but its proof in Section \ref{sec-simplepropvarphi} is independent and self-contained.  

It is known and a standard task to show that the map $\rho\mapsto \chi^{\ssup{v=0}}(\rho,0)=\inf_{m\in[0,\infty)^\N\colon \sum_k km_k=\rho}H(m|q)$ has a phase transition at $\rho^{\ssup{v=0}}_{\rm c}=\sum_k k q_k$, if this series has a finite value. Indeed, $\chi^{\ssup{v=0}}(\cdot,0)$ is positive and strictly decreasing left of $\rho^{\ssup{v=0}}_{\rm c}$ (with slope $-\infty$ at $\rho=0$) and there is a minimiser $m$, while for supercritical $\rho$, there is none, and it is constantly equal to zero.

\subsection{Further results: existence and regularity of minimiser(s)}\label{sec-Discussion}

Let us discuss the existence of minimisers in the variational formulas on the right-hand sides of \eqref{freeenergyvarform},  \eqref{chi(m,rho2)def} and \eqref{varphidef}. First we turn to \eqref{freeenergyvarform}. Indeed, we will be proving in Section \ref{sec-VarformProof} the following. Recall that we assume that the $k$-box densities of the reference PPP satisfy $q_k=\e^{o(k)}$ as $k\to\infty$. 

\begin{lemma}\label{lem-Varformident} For any $\rho\in(0,\infty)$, the minimum on the right-hand side  of \eqref{freeenergyvarform} is achieved; more precisely,
\begin{equation}\label{chiidentvarphi}
\begin{aligned}
\inf_{ \rho_{\rm mi}, \rho_{\rm ma}\geq 0\colon  \rho_{\rm mi}+ \rho_{\rm ma}=\rho}\chi( \rho_{\rm mi}, \rho_{\rm ma})&=\chi(\rho,0)\\&=\inf\Big\{I(P)+P(\Phi_{0,\Z^d})\colon P\in\Mcal_1^{\ssup{\rm s}}(\Omega), P(N_0^{\ssup{\ell}})=\rho\Big\}.
\end{aligned}
\end{equation}
\end{lemma}

This may be a bit surprising, as it says that the free energy can be described by exclusively looking at limiting configurations without condensate part, no matter if the condensation phase transition takes place or not. In particular, the existence of this phase transition cannot be discussed by exclusively looking at the value of the limiting free energy.

It is clear that \lq$\leq$\rq\ is trivial in \eqref{chiidentvarphi}. We did not find any direct, analytical proof of Lemma \ref{lem-Varformident}. Instead, the proof of \lq $\geq$\rq\ uses some inspiration from the proof of Theorem \ref{thm-freeenergygrid}; it constructs from a constrained partition function with an explicit macro-part in the configuration a constrained partition function without macro-part, but with increased micro-part, and shows that this manipulation does not increase the latter on the exponential scale. 

In the sequel we will abbreviate $\chi(\rho)=\chi(\rho,0)$ for the free energy with particle density $\rho$. We turn now to the question about whether or not a minimising particle configuration exists in terms of the sequence $m=(m_k)_{k\in\N}$ of the $k$-mark densities, the family of marked point processes $(P_a)_{a\in\N_0}$ and the distribution $\psi$ of the macro parts. The answer is positive:

\begin{lemma}[Minimisers of $\varphi$]\label{lem-ProofExistMin} For any $\rho\in[0,\infty)$, there is at least one minimiser $(m,\psi)$ for the variational problem 
\begin{equation}\label{crucialVP}
\chi(\rho)=\min\Big\{\varphi(m,\psi)\colon m\in [0,\infty)^\N,\psi\in\Mcal_1(\N_0), \rho=\sum_{k\in\N}k m_k+\sum_{a\in\N_0}a\psi(a)\Big\},
\end{equation}
i.e., at least one minimiser $(m,\psi)$ of $\varphi$ defined in~\eqref{varphidef} under the constraint $\sum_{k\in\N}k m_k+\sum_{a\in\N_0}a\psi(a)=\rho$. Furthermore, there is at least one minimiser $(m,\psi)$ such that $\psi$ has no more than two atoms.

Additionally, for any $(m,\psi)$, there is at least one minimiser $(P_a)_{a\in\N_0}$ for the variational formula in \eqref{varphidef}.
\end{lemma}

The proof of the existence of a minimiser is already at the end of Section \ref{sec-proofgrid}; it is a by-product of our proof of Theorem \ref{thm-freeenergygrid}. Along subsequences of approximately minimising sequences of $(m,\psi)$'s, the total microscopic particle density $\sum_k k m_k$ can in principle become smaller; there is the possibility of losing mass. However, our proof shows that such a loss can be dispended with by assigning it to the macro part with the help of some manipulations on the level of  particle configurations.
Much more information about properties of a particular minimiser is provided in Lemma~\ref{lem-saturation}.

A closer look at all the minimisers, in particular the proof that $\psi$ can be taken as a Dirac measure or as a mixture of two Dirac measures, is in Section \ref{sec-minimivarphi}.
It uses convexity arguments for a reformulation of the characteristic variational formula in \eqref{crucialVP}: carrying out only the infimum over $m$, making the substitution $\rho=\sum_{a\in\N_0}\psi(a)(a+P_a(N_0^{\ssup\ell}))$ and $\rho_a= P_a(N_0^{\ssup\ell})$, we see (also using Lemma \ref{lem-Varformident}) that
\begin{equation}\label{A(psi)def}
\chi(\rho)=\inf_{\psi\in\Mcal_1(\N_0)}\Big[ \v\sum_{a\in\N_0}\psi(a) a^2+\inf_{(\rho_a)_{a\in\N_0}\colon \rho=\sum_{a\in\N_0}\psi(a)[a+\rho_a]}\sum_{a\in\N_0}\psi(a)[2\v a\rho_a+\chi(\rho_a)]\Big].
\end{equation}

Now we turn to regularity properties of the map $m\mapsto \varphi(m,\psi)$:

\begin{lemma}[Differentiability of $\varphi$]\label{lem-DerivativeMicro}
Fix $\psi\in\Mcal_1(\N_0)$ satisfying $\sum_{a\in\N_0}a\psi(a)<\infty$. Then, for any $m\in[0,\infty)^\N$ and for any $k$ satisfying $m_k>0$, $\varphi(\cdot,\psi)$ is differentiable in $m_k$, and the partial derivative satisfies
\begin{equation}\label{DerivativeMicro_Eq}
\partial_{m_k}\varphi(m,\psi)=\log\frac{m_k}{q_k}+ t_k- \sup_{(P_a)_a}\log\sum_{a\in \N_0} \psi(a)\e^{-2\v ak}P_a(\e^{-2\Phi^{\ssup k}}), 
\end{equation}
where the supremum is over all minimisers $(P_a)_{a\in\N_0}$ in the formula \eqref{varphidef} of $\varphi(m,\psi)$ and
\begin{equation}
t_k=T_{0,0}(G_k,G_k)\qquad\text{and}\qquad
\Phi^{\ssup k}(\omega)=\sum_{x\in \Z^d}\sum_{l\in\N}\xi^{\ssup l}(x)T_{0,x}(G_k,G_l),
\end{equation} 
are the self-interaction of a $k$-mark and the interaction between $\omega$ and a $k$-mark at the origin, i.e., the configuration $\delta_{(0,G_k)}$, respectively.
\end{lemma}

The proof of Lemma \ref{lem-DerivativeMicro} is in Section \ref{sec-varphidiffble}. Let us note that we have no information about uniqueness or non-uniqueness of minimisers in \eqref{varphidef}, since the entropy density $I$ is affine hence not strictly convex. 
With the help of Lemma~\ref{lem-DerivativeMicro} we can derive the variational equations for minimisers of $\varphi(m,\psi)$: 
\begin{lemma}[Euler--Lagrange equations]\label{lem-ELeq}
Fix $\rho\in(0,\infty)$, then, for any minimiser $(m,\psi)$ of \eqref{crucialVP}, $m$ satisfies the Euler--Lagrange equations
\begin{equation}\label{ELequations}
m_k=q_k \e^{\alpha k-t_k}\sup_{(P_a)_a}\sum_{a\in\N_0}\psi(a)\e^{-2{\v}ak}P_a\big(\e^{-2\Phi^{\ssup k}}\big),\qquad k\in\N,
\end{equation}
where $\alpha\in\R$ is the Lagrange multiplier. 
\end{lemma}

The proof of Lemma \ref{lem-ELeq} is in Section \ref{sec-EL-eqs}.

\begin{cor}[Differentiability of $\chi(\cdot,0)$]\label{cor-chidiff}
The map $\rho\mapsto\chi(\rho)=\chi(\rho,0)$ is differentiable in $(0,\infty)$ with $\chi'(\rho)=\alpha$, the Euler--Lagrange parameter of Lemma \ref{lem-ELeq}.
\end{cor}

The proof of Corollary \ref{cor-chidiff} is in Section \ref{sec-EL-eqs}. As a consequence, a possible phase transition (if it exists) cannot be of first order. 

Alternatively to the formula in \eqref{DerivativeMicro_Eq}, in Section \ref{sec-thinning} we present another formula for the derivative of $\varphi$ with respect to $m$:
\begin{equation}\label{DerivativeMicro_EqAlter}
\partial_{m_k}\varphi(m,\psi)=\log\frac{m_k}{q_k}-t_k+\sup_{(P_a)_a}\log\sum_{a\in \N_0} \psi(a)\e^{2\v ak}\frac{P_a(N_0^{\ssup{\delta_k}}\,\e^{2\Phi^{\ssup k}})}{m_k}.
\end{equation}
Since we  do not use this formula for further study and since a proof would be quite technical, we decided to restrict to giving a heuristic argument for how to derive \eqref{DerivativeMicro_EqAlter}. Roughly speaking, this strategy goes via a uniform random thinning procedure, while the proof of \eqref{DerivativeMicro_Eq} is via adding a uniformly distributed Poisson point.

\subsection{On the phase transition}\label{sec-BEC}

Let us now discuss the nature of a possible condensation phase transition, assuming that it exists. This hinges on the minimisers $(m,\psi)$ that we established in  Lemma \ref{lem-ProofExistMin}, i.e., the minimisers of the problem \eqref{crucialVP} for a given $\rho$. It is convenient to introduce the corresponding  microscopic and macroscopic particle densities:
\begin{equation}\label{micromacromass}
\rho_{\rm mi}(m)=\sum_{k\in\N}k m_k,\qquad \rho_{\rm ma}(\psi)=\sum_{a\in\N_0}a\psi(a),\qquad \rho_{\rm mi}(m)+\rho_{\rm ma}(\psi)=\rho.
\end{equation}
We do not know if $\rho_{\rm mi}$ and $\rho_{\rm ma}$ are independent of the choice of the minimiser $(m,\psi)$. Actually, we somehow characterised all the minimisers $(m,\psi)$ for a given $\rho$ in the proof of Lemma \ref{lem-ProofExistMin}, but have no information about their uniqueness.

The occurrence of a non-trivial macroscopic particle density is characterised by the existence of a minimiser $(m,\psi)$ such that $\rho_{\rm ma}(\psi)>0$, i.e., $\psi\not=\delta_0$. 
We define the {\em critical particle density} for the emergence of a macroscopic part as follows:
\begin{equation}\label{critvalue}
\begin{aligned}
\rho_{\rm c}
&=\sup\Big\{\rho\in(0,\infty) \colon \mbox{\eqref{crucialVP} has a minimiser }(m,\delta_0), \mbox{ i.e., }\rho_{\rm mi}(m)=\rho\Big\}.
\end{aligned}
\end{equation}
Let us first note that there is indeed a non-condensate phase:

\begin{lemma}[Positivity of $\rho_{\rm c}$]\label{lem-rhocpos}
The critical particle density $\rho_{\rm c}$ is positive. More precisely, $\rho_{\rm c}\geq P(N_0^{\ssup\ell})$ for any $P\in \Mcal_1^{\ssup{\rm s}}(\Omega)$ that minimises $P\mapsto I(P)+P(\Phi_{0,\Z^d})$. In other words, $\rho_{\rm c}$ is not smaller than the smallest minimiser of $\rho\mapsto\chi(\rho)$.
\end{lemma}
Note that there is at least one such minimiser $P$ since the level sets of $I$ are compact and $P\mapsto P(\Phi_{0,\Z^d})$ is lower semi-continuous.
 
\begin{proof} Note that the existence of a minimiser $(m,\delta_0)$ in  \eqref{crucialVP} is equivalent to the existence of a minimiser $P$ in \eqref{chiidentvarphi} via taking $m_k=P(N_0^{\ssup{\delta_k}})$, see the last sentence in Lemma~\ref{lem-ProofExistMin}.

Consider the variational formula $\widetilde \chi(\rho)=\inf\{I(P)+P(\Phi_{0,\Z^d})\colon P\in\Mcal_1^{\ssup{\rm s}}(\Omega), P(N_0^{\ssup\ell})\leq  \rho\}$. For any $\rho\geq 0$, this formula possesses a minimising $P$, since $I$ has compact level sets  and the maps $P\mapsto P(\Phi_{0,\Z^d})$ and $P\mapsto P(N_0^{\ssup\ell})$ are lower semi-continuous. Indeed, if $(P_n)_{n\in \N}$ is a sequence of admissible approximate minimisers, then $(I(P_n))_{n\in \N}$ is bounded, since $\Phi_{0,\Z^d}\ge 0$. Hence, $(P_n)_{n\in \N}$ has a convergent subsequence with admissible limit. By lower semi-continuity of $P\mapsto I(P)$ and $P\mapsto P(\Phi_{0,\Z^d})$, the limit is a minimiser. 

Recall from Lemma \ref{lem-properties}(6) that $\rho\mapsto \chi(\rho,0)=\chi(\rho)$ strictly decays in a neighbourhood of $0$. By convexity, it is even strictly decreasing precisely in the interval $(0,\rho_{\rm min}]$, with $\rho_{\rm min}$ being the smallest minimiser of $\chi$. On this interval, $\chi$ and $\widetilde \chi$ coincide. Indeed, pick a minimiser $P$ for $\widetilde \chi(\rho)$ satisfying $\widetilde \rho=P(N_0^{\ssup\ell})<\rho$, then $P$ would be admissible also for $\chi(\widetilde \rho)$ and would imply that $\chi(\widetilde \rho)\leq \widetilde \chi(\rho)\leq \chi(\rho)$, which contradicts the strict monotonicity. Hence, every minimiser $P$ for $\widetilde \chi(\rho)$ satisfies $P(N_0^{\ssup\ell})=\rho$ and is  therefore also a minimiser for $\chi(\rho)$.  This implies that $\rho_{\rm c}\geq \rho_{\rm min}$.
\end{proof}

We say that a condensation phase transition occurs if $\rho_{\rm c}$ is finite, i.e., if for any sufficiently large $\rho$ any minimiser $(m,\psi)$ satisfies $\rho_{\rm ma}(\psi)>0$. This notion of a condensate phase transition is analogous to the famous Bose--Einstein condensation, which is conjectured to occur in dimensions $d\geq 3$ (but not in $d\in\{1,2\}$), and has been proved to occur in the non-interacting case where $v=0$ (see also Section \ref{sec-Bosegas}). There the dependence on the dimension is clearly seen to hinge on the summability of $k q_k$ over $k\in\N$; note that $kq_k=(4\pi\beta k)^{-d/2}$ in the interacting Bose gas. In our model, where we admit an arbitrary summable sequence $(q_k)_{k\in\N}$, the conjecture is suggested that the occurrence of the phase transition needs the summability of $k q_k$ as well. As a consequence, we would see that in our model the occurrence is not a dimensionality question, but only a summability question.

In the interacting Bose gas, the condensation phase transition is conjectured to be of {\em saturation} type, by which we mean that the particles organise in a microscopic part of density $\rho\wedge \rho_{\rm c}$ and a condensate of density $[\rho-\rho_{\rm c}]_+$, 
where $x_+$ denotes the positive part of $x$. This reflects the understanding that, if we consider increasing $\rho$, for all small values, the entire particle cloud is organised in microscopic marks, and as soon as $\rho$ exceeds $\rho_{\rm c}$, then additional particles are put into the macro part, but the total mass $\rho_{\rm c}$ of the micro part is not changed anymore. Hence, in the box version we might expect that we should have
\begin{equation}\label{saturation}
\sum_{k\in\N}k m_k=\rho\wedge \rho_{\rm c}\quad\mbox{and}\quad\sum_{a\in\N_0}a\psi(a)=\big[\rho-\rho_{\rm c}\big]_+,\quad\mbox{for any minimiser }(m,\psi).
\end{equation}

However, it turns out that this is not the case:

\begin{lemma}[Qualitative description in case of a phase transition]\label{lem-saturation}
Assume that $\rho_{\rm c}$ is finite. Then \eqref{saturation} is false. Instead, the following is true. 
\begin{enumerate} 
\item There exists $\rho_{\rm t}\in[\rho_{\rm c}\lor 1, \rho_{\rm c}+1)$ such that
\begin{equation}\label{rho_t}
\chi(\rho_{\rm t}-1)+(2\rho_{\rm t}-1)\v+(\rho_{\rm c}-\rho_{\rm t})(\chi'(\rho_{\rm t}-1)+2\v)=\chi(\rho_{\rm c}).
\end{equation}

\item For any $a\in\N_0$, for any $\rho\in[\rho_{\rm c},\rho_{\rm c}+1)$,
\begin{equation}\label{chiidentnew}
 \chi(\rho+a)=\begin{cases}
 \frac{\rho-\rho_{\rm c}}{\rho_{\rm t}-\rho_{\rm c}}(\chi(\rho_{\rm t}-1)+\v(2\rho_{\rm t}-1))+\frac{\rho_{\rm t}-\rho}{\rho_{\rm t}-\rho_{\rm c}}\chi(\rho_{\rm c})+2\v a \rho +\v a^2&\mbox{if }\rho\leq \rho_{\rm t},\\
 \chi(\rho-1)+2\v(\rho-1)(a+1)+\v (a+1)^2&\mbox{if }\rho\geq  \rho_{\rm t}.
\end{cases}
\end{equation}
As a consequence, we have $\chi'(\rho_{\rm c})=\chi'(\rho_{\rm t}-1)+2\v$.

\item We assume that for any $\rho\in[0,\rho_{\rm c}]$, there exists a minimiser $(m_{\rho},\delta_0)$ with density $\rho$. Then, for any $\rho\in[\rho_{\rm c},\infty)$ there is a minimiser $(m,\psi)$ of $\varphi$ with density $\rho$ such that, for  $a\in\N_0$:
\begin{itemize}
    \item if $\rho\in[\rho_{\rm c}+a, \rho_{\rm t}+a]$, the minimiser is a proper convex combination of  $(m_{\rho_{\rm c}},\delta_a)$ and  $(m_{\rho_{\rm t}-1},\delta_{a+1})$, with $(m_{\rho_{\rm c}},\delta_0)$ and $(m_{\rho_{\rm t}-1},\delta_0)$  minimisers at densities $\rho_{\rm c}$ respectively $\rho_{\rm t}-1$;
    \item if $\rho\in [\rho_{\rm t}+a, \rho_{\rm c}+a+1]$, the minimiser is $(m_{\rho-a-1},\delta_{a+1})$ with $(m_{\rho-a-1},\delta_0)$ a minimiser with density $\rho-a-1$. 
\end{itemize}

\item The microscopic total mass $\rho_{\rm mi}$ corresponding to the minimiser $(m,\psi)$ of (3) is one-periodic in $\rho$ in the interval $[\rho_{\rm t}-1,\infty)$. It increases linearly from $\rho_{\rm t}-1$ to $\rho_{\rm c}$ in $[\rho_{\rm t}-1,\rho_{\rm c}]$. It decreases linearly from $\rho_{\rm c}$ to $\rho_{\rm t}-1$ in $[\rho_{\rm c},\rho_{\rm t}]$. Correspondingly, the macroscopic total mass $\rho_{\rm ma}$ corresponding to the minimiser $(m,\psi)$ of (3) is constant equal to $a$ on $[\rho_{\rm t}-1+a, \rho_{\rm c}+a]$, and increases linearly from $a$ to $a+1$ on $[\rho_{\rm c}+a, \rho_{\rm t}+a]$, for any $a\in \N_0$.
\end{enumerate}
\end{lemma}

The proof of Lemma \ref{lem-saturation} is in Section \ref{sec-naturephasetrans}; it uses that, under \eqref{saturation}, $\chi(\cdot,0)$ can be shown to be not differentiable in any point of $\rho_{\rm c}+\N$, in contradiction to Lemma \ref{cor-chidiff}.
(1) follows from the facts that $\chi'(\rho_{\rm c})>(2\rho_{\rm c}-1)\v$ and $\chi(\rho_{\rm c})\leq\chi(\rho_{\rm c}-1)+(2\rho_{\rm c}-1)\v$. Furthermore, the analysis of the minimisers of the variational formula that was done in Section \ref{sec-minimivarphi} is crucial as well.

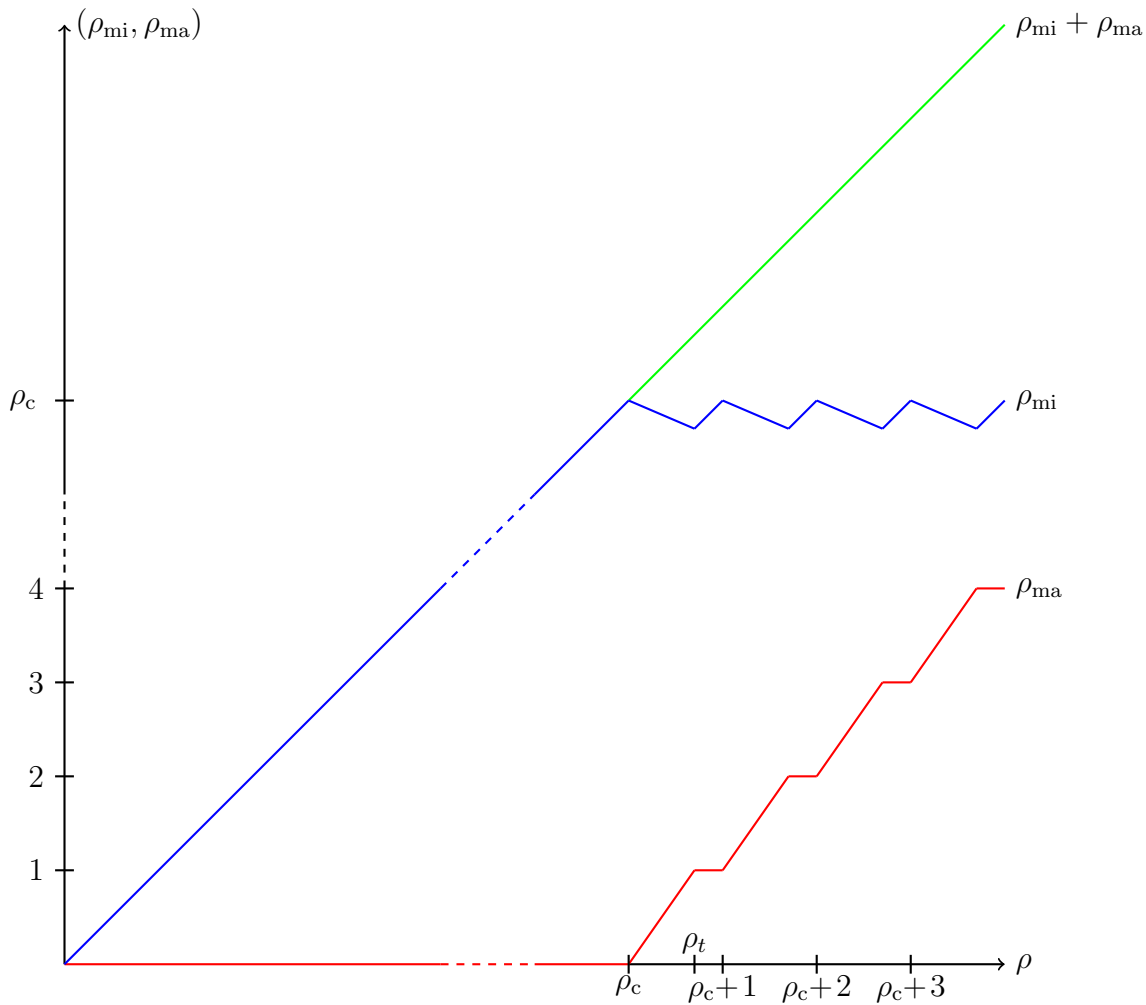
\begin{figure}[!htpb]
	\centering
\input{Pix_rho.tex}
	\caption{Plot of the microscopic total mass $\rho_{\rm mi}=\sum_{k\in\N} k m_k$ (blue) and the macroscopic mass $\rho_{\rm ma}=\sum_{a\in\N}a\psi(a)$ (red) as functions of $\rho$  if the critical threshold $\rho_{\rm c}$ is finite.}
	\label{chiplot}
\end{figure}

An  illustration of the micro- and macroscopic total masses as functions of $\rho$ is in Figure~\ref{chiplot}. In words, as $\rho$ increases from zero to infinity, then, first each optimal strategy organises all particles in microscopic boxes. If this changes at some finite $\rho_{\rm c}$, then, if $\rho$ further increases, it is an optimal strategy to cover a certain percentage of the space with one macroscopic box and to reduce the microscopic particle density linearly in $\rho$, until a second critical threshold $\rho_{\rm t}$ is reached, at which the whole space is covered by one macroscopic box. Further increasing $\rho$, additional microscopic mass is added without changing the macroscopic part until $\rho_{\rm c}+1$ is reached. This procedure is then iterated by further adding macroscopic boxes. 

\section{Discussion}\label{sec-phasetrans+Lit}

In Section \ref{sec-phasetransexist} we explain the difficulty in finding a proof for a phase transition, and in Section \ref{sec-literature} we mention some related works.

\subsection{Does the phase transition occur?}\label{sec-phasetransexist}
We did not yet touch the most interesting question: under what conditions on $(q_k)_{k\in\N}$ and $v$ does the micro-macro phase transition exist that is analogous to BEC, i.e., under what conditions is  $\rho_{\rm c}$ finite? We do not give any answer to this question in the present paper and leave this open problem to future work. But we would like to comment on that now.

Note that the question about minimisers decomposes into many independent questions about the existence or non-existence of a minimising $(m,\psi)$ with $\psi=\delta_0$ or $\psi\not=\delta_0$. Even though we somehow characterised all the minimisers $(m,\psi)$ in the proof of Lemma~\ref{lem-ProofExistMin}, all the answers that we can give are summarised in Lemma~\ref{lem-rhocpos}. Beyond this, we are not able to say anything descriptive about the set of $\rho$'s for which there is or there is not a minimiser with or without nontrivial macro part, not even whether or not it is convex, i.e., an interval. One could define another critical density $\widetilde \rho_{\rm c}$ as the supremum over all $\rho$ such that {\em every} minimiser $(m,\psi)$ satisfies $\rho_{\rm mi}(m)=\rho$, then $0\leq \widetilde\rho_{\rm c}\leq\rho_{\rm c}$. The finiteness of $\widetilde \rho_{\rm c}$ also would imply the existence of a minimising configuration with non-trivial condensate but possibly an additional minimising configuration without condensate. 

One way to attack the question is by looking at the Euler--Lagrange equations for a possible minimiser and giving arguments in favour or against its existence.  The analogous equation for the free Bose gas (see the end of Section \ref{sec-Bosegas}) reads $m_k=q_k\e^{\alpha k}$ with some Lagrange multiplier $\alpha$, subject to the constraint $\rho=\sum_{k\in\N}k q_k\e^{\alpha k}$. If $\sum_k k q_k$ is finite, then this is the largest value that can be reached by proper choice of $\alpha$. The conclusion is that, for $\rho>\rho_{\rm c}=\sum_k k q_k$, there is no Lagrange multiplier $\alpha$ and therefore no minimiser $m$, but for $\rho\leq \rho_{\rm c}$, there is one.

In the box version of the interacting Bose gas, there is always a minimiser $(m,\psi)$ of the variational formula for $\chi(\rho)$, but the question is now about the existence of a  minimiser that has no macroscopic part, i.e., is of the form $(m,\delta_0)$.  The characteristic equation for that reads
\begin{equation}\label{ELdeltazero}
m_k=q_k \e^{\alpha k-t_k}\sup_P P\big(\e^{-2\Phi^{\ssup k}}\big),\qquad k\in\N,
\end{equation}
where the supremum is taken over all minimisers $P$ of $I(P)+P(\Phi_{0,\Z^d})$ subject to $\rho=\sum_{k\in\N}kP(N_0^{\ssup{\delta_k}})$, and the Lagrange multiplier $\alpha$ needs to satisfy $\rho=\sum_{k\in\N}k q_k \e^{\alpha k-t_k}\sup_P P(\e^{-2\Phi^{\ssup k}})$. Recall that $t_k=T_{0,0}(G_k,G_k)$ is the self-interaction of a mark $G_k$. Note that $P$ may {\em a priori} depend on $k$, since we have no information about uniqueness of the minimiser $P$. 

The most important difference to the Euler--Lagrange equation for the free Bose gas is the appearance of the last term, $\sup_P P(\e^{-2\Phi^{\ssup k}})$. We found no way to utilise this formula for deriving interesting information about the existence of the phase transition or other details. It is likely that the critical $\alpha$ for the largest $\rho$ satisfies that the large-$k$ exponential rate of the summands $q_k \e^{\alpha k-t_k}\sup_P P(\e^{-2\Phi^{\ssup k}})$ is equal to zero, and then some explicit information about the second-order term of $P(\e^{-2\Phi^{\ssup k}})$ is necessary, and here we do not see any ansatz to identify or derive that.

We remark that the internal interaction of $G_k$ behaves like
\begin{equation}\label{tkasy}
t_k= \v k- (C+o(1)) k^{1-\frac 1d},\qquad k\to\infty, \mbox{ for some }C>0,
\end{equation}
where the second-order term comes from boundary effects. This is a clear difference to the interacting Bose gas, as cycles have no beginning nor an end and therefore no boundary effect, but the box $G_k$ has, and its internal energy explicitly appears here.

In \cite{BKM21}, a comparable situation, where only interactions within the marks are considered could partially be solved in this respect. In comparison to \eqref{ELdeltazero}, the term $t_k$ is missing there, and instead of $P(\e^{-2\Phi^{\ssup k}})$, there is the expectation of a single Brownian bridge with time interval $[0,k \beta]$ with exponential interaction between any two legs, like in the interacting Bose gas. With the help of an extension of the lace expansion technique it is proved in \cite{BKM21} that this expectation behaves, for large $k$, as $\e^{C k} k^{-d/2}(1+o(1))$ with $C$ a characteristic quantity, however, only for sufficiently small interaction potential and only in dimensions $d\geq 5$. This (very fine!) asymptotics made it possible to derive the existence of that phase transition in $d\geq 5$.

Another possible route to prove the occurrence of the phase transition might be to prove that, for some large $\rho$, for any $m\in[0,\infty)^\N$ satisfying $\sum_k k m_k=\rho$, there is $\widetilde m\in[0,\infty)^\N$ satisfying $\sum_k k \widetilde m_k=\rho-1$ such that $\varphi(m,\delta_0)>\varphi(\widetilde m,\delta_0)+(2\rho-1)\v$, where we note that the right-hand side is equal to $\varphi(\widetilde m,\delta_1)$. This would show explicitly that it is not optimal to organise the entire total particle mass in microscopic boxes, but one part of it in the regular grid $\Z^d$ and the remaining part (i.e., total mass $\rho-1$) in microscopic boxes. We tried to prove this assertion with the help of several of the techniques that proved successful in Sections \ref{sec-EL-eqs} and \ref{sec-thinning} for handling derivatives with respect to $m$ (de-Poissonisation and thinning), but our ansatzes were not fine enough, partially since we perturbed only with independent processes, which seem to be not well enough adapted.

There might be a phenomenological connection between the finiteness of $\rho_{\rm c}$ and the uniqueness of minimisers $P$ in the variational formula. Indeed, the conjecture is tempting that the minimiser $P$ is a Gibbs measure for a related potential, and that the Gibbs measures are unique precisely in the sub-critical case. However, we have no clue about existence nor uniqueness of Gibbs measures nor about uniqueness of a minimising $P$ nor about how this knowledge could help in the analysis of the free energy.

\subsection{Literature remarks}\label{sec-literature}

Let us give a small survey on the literature on treatments of the Bose gas with the help of the theories of random point processes and of large deviations.

The starting point of this line of research is the Feynman-Kac formula for $N$ interacting Brownian cycles (bridges) with uniformly distributed starting sites in a box and a symmetrisation. Using the Markov property, this formula  can be turned into a random ensemble of closed cycles with various lengths and independent starting/terminating sites. See e.g.~\cite{U06a} for a pedagogical explanation. In \cite{ACK10}, an additional step is made by rewriting this ensemble explicitly in  terms of an interacting marked Poisson point process.

It has been noticed in \cite{F91} that the probabilistic concept of a {\em random point process} is highly appropriate. There are some efforts undertaken to construct interacting marked Gibbs point measures that show the suggested characteristics of the interacting Bose gas in the thermodynamic limit, however without any clear relation to the limiting free energy of the partition function of the gas. It is very likely that some of these target Gibbs measures will sooner or later turn out to describe the microscopic part of the interacting Bose gas and to be minimisers in a characteristic variational formula as in \cite{ACK10}, but this is widely open. Here we would like to mention H.~Zessin and co-workers and students, see \cite{NPZ13,RZ20}. Also in several papers by H.~Tamura and co-workers,  point-process descriptions are employed to gain an understanding of the Bose gas, see for example \cite{TI06}. In a series of papers by a team around J.~Fröhlich (see the summary \cite{FKSS20}), various rescalings and limiting regimes of the interacting Bose gas are examined with mathematically highly involved methods, however these techniques do not have much to do  with point process theory. In \cite{V20}, local limits of the trace of the non-interacting Brownian cycle loop soup towards the Brownian interlacement process (a Poisson point process on the set of infinitely long Brownian paths) is proved, which is a non-trivial step towards an understanding of the condensate, but still far away from handling the free energy. Another work in this vein was recently done in \cite{AFY19}.

However, for handling the thermodynamic limit of such an interacting marked Poisson point process, there are only very few  investigations in the literature. Based on the theory of so-called {\em level-three large deviation principles}, introduced by Donsker and Varadhan in the 1970s, adaptations to marked Poisson point processes both in $\Z^d$ and $\R^d$ were developed in \cite{G93, GZ93, G94}. This concept is rather suitable for handling limiting free energies of partition functions like the one for the interacting Bose gas, which was also mentioned in these papers, but not carried out. This was done for the first time in \cite{ACK10}, the starting point of the present paper (see Section \ref{sec-Bosegas} for a summary). There a {\em characteristic variational formula} was derived for the limiting free energy, which reflects  and encodes all thermodynamic quantities that are relevant for understanding the gas, the most prominent of which are the {\em energy} and the {\em entropy} and the {\em effective density} of the particle configuration. The method used in the present paper does not rely on that large-deviation principle, but carries out the arguments via another route; it is very much in spirit of \cite{G93, GZ93, G94}. However, for handling also the macroscopic part, some new techniques had to be found, and this is one of the new contributions of the present paper.

The model introduced and studied in the present paper is in the spirit of models that are called {\it polydisperse mixture models} in the physics literature. Recently, in \cite{J20} a mutually repellent version was studied under the additional assumption of hierarchy of the droplet configuration, the droplets being deterministic, randomly placed discrete cubes. The assumption of hierarchy made it possible to derive a formula for the limiting free energy via a decomposition according to the hierarchies. It is certainly desirable to overcome this assumption. The method developed in the present paper seems to give a suitable ansatz to do so in future work.

\section{Derivation of the variational formula: proof of Theorem \ref{thm-freeenergygrid}}\label{sec-proofgrid}

\noindent In this section, we prove Theorem \ref{thm-freeenergygrid}. It will be sufficient to prove the assertion for the restriction to Dirichlet boundary condition, i.e., 
\begin{eqnarray}
\liminf_{N\to\infty}\frac1{|\Lambda_N|}\log Z_{N,\L_N,{\rm Dir}}&\geq& -\inf_{ \rho_{\rm mi}, \rho_{\rm ma}\ge 0\colon \rho_{\rm mi}+ \rho_{\rm ma}=\rho}\chi( \rho_{\rm mi}, \rho_{\rm ma}),\label{ProofThm1lower}
\\
\limsup_{N\to\infty}\frac1{|\Lambda_N|}\log Z_{N,\L_N,{\rm Dir}}&\leq& -\inf_{ \rho_{\rm mi}, \rho_{\rm ma}\ge 0\colon \rho_{\rm mi}+ \rho_{\rm ma}=\rho}\chi( \rho_{\rm mi}, \rho_{\rm ma}),\label{ProofThm1}
\end{eqnarray}
for centred boxes $\Lambda_N$ with volume $N/\rho$, where we put
\begin{equation}
Z_{N,\L_N,{\rm Dir}}=\E_{\L_N}\big[{\rm e}^{-\Phi_{\Lambda_N,\Lambda_N}}\1\{N^{\ssup\ell}_{\Lambda_N}=N\}\1\{\widetilde N_{\L_N^{\rm{c}}}=0\}\big].
\end{equation}
We prove \eqref{ProofThm1lower} in Section \ref{sec-lowbound} and \eqref{ProofThm1} in Section~\ref{sec-uppbound}. Indeed, the assertion of Theorem \ref{thm-freeenergygrid}  for zero Dirichlet boundary condition follows since the probability of the conditioning event is $\e^{o(|\L_N|)}$: 
$$
\begin{aligned}
\P_{\L_n}(\widetilde N_{\L_N^{\rm{c}}}=0)=\P\Big(\bigcap_{k\in\N}\{M^{\ssup{\delta_k}}_{\L_N,\L_N^{\rm{c}}}=0\}\Big)&=\prod_{x\in\L_N}\prod_{k\in\N\colon x+ G_k\not\subset \L_N}\P(N^{\ssup{\delta_k}}_x=0)\\
&=\prod_{k\in\N}\exp\Big\{-q_k\#\{x\in\L_N\colon  x+ G_k\not\subset \L_N\}\Big\}.
\end{aligned}
$$
Now it is easy to see that the right-hand side is $\e^{-o(|\L_N|)}$, since the cardinality is $o(|\L_N|)$ for any $k\in\N$, and $\sum_k q_k<\infty$.

\subsection{Proof of the lower bound in Theorem \ref{thm-freeenergygrid}}\label{sec-lowbound}

Let $ \rho_{\rm mi}, \rho_{\rm ma}\geq 0$ be given such that $ \rho_{\rm mi}+ \rho_{\rm ma}=\rho$. Our goal is to show that, for any $m\in[0,\infty)^\N$ and $\psi\in\Mcal_1(\N_0)$ satisfying $\sum_k k m_k= \rho_{\rm mi}$ and $\sum_{a} a\psi(a)= \rho_{\rm ma}$, that
\begin{equation}\label{lowbound1}
\liminf_{N\to\infty}\frac 1{|\Lambda_N|}\log Z_{N,\L_N,{\rm Dir}}\geq -\varphi(m,\psi),
\end{equation}
with $\varphi$ as in \eqref{varphidef}. This implies \eqref{ProofThm1lower}.

\begin{step}\label{step:compactsupp}
We reduce the problem to $m$ and $\psi$ having finite support.
\end{step}

We define $\widetilde \chi( \rho_{\rm mi}, \rho_{\rm ma})$ as $\chi( \rho_{\rm mi},  \rho_{\rm ma})$ with the additional constraint that $m$ and $\psi$ must have finite support. We will show that $\widetilde \chi( \rho_{\rm mi}, \rho_{\rm ma})\leq\chi( \rho_{\rm mi},  \rho_{\rm ma})$. Pick $m=(m_k)_{k\in\N}$ satisfying $\sum_{k\in\N}k m_k= \rho_{\rm mi}$, and a probability measure $\psi$ on $\N_0$ satisfying $\sum_{a\in\N_0} a\psi(a)= \rho_{\rm ma}$. We introduce cut-off versions of $m$ and $\psi$ using large auxiliary parameters $K,A\in\N$ by putting $m^{\ssup{\leq K}}=(m_1,m_2,\dots,m_K,0,0,\dots)$ and $\psi^{\ssup A}=\psi+\sum_{a>A} \psi(a)(\delta_0-\delta_a)$. Let us show that $\varphi(m^{\ssup K},\psi^{\ssup A})\leq\varphi(m^{\ssup K},\psi)\leq \varphi(m,\psi)+\sum_{k>K} q_k$. 

The first inequality comes from Lemma \ref{lem-properties}. For the second one, assume that $(P_a)_{a\in\N_0}$ is admissible for the formula for $\varphi(m,\psi)$. Then it is clear that $(P_a^{\ssup{\leq K}})_{a\in\N_0}$, with $P_a^{\ssup{\leq K}}$ defined from $P_a$ by suppressing all marks with size $>K$, is admissible for the formula for $\varphi(m^{\ssup {\leq K}},\psi)$, furthermore $I(P_a^{\ssup{\leq K}})\leq I(P_a)+\sum_{k>K} q_k$, and also the energy terms of $(P_a^{\ssup{\leq K}})_{a\in\N_0}$ are not bigger than those of $(P_a)_{a\in\N_0}$. Hence the inequality.

We deduce that
$\widetilde \chi(\sum_{k\in[K]}k m_k,\sum_{a\in[A]} a\psi(a))\leq \varphi(m^{\ssup {\leq K}},\psi^{\ssup A})\leq \varphi(m,\psi)+\sum_{k>K} q_k$. Letting $K,A\to \infty$ and using continuity of $\widetilde \chi$ (see Lemma \ref{lem-properties}) and summability of $(q_k)_{k\in\N}$, we get $\widetilde \chi( \rho_{\rm mi}, \rho_{\rm ma})\leq \varphi(m,\psi)$. Taking the infimum over $m$ and $\psi$, we obtain the desired inequality $\widetilde \chi( \rho_{\rm mi}, \rho_{\rm ma})\leq\chi( \rho_{\rm mi},  \rho_{\rm ma})$.

\begin{step}\label{Step_2}
We construct a configuration adapted to $m$ and $\psi$.
\end{step}

Pick $m=(m_k)_{k\in\N}$ having support included in $[K]=\{1,\dots,K\}$  for some $K\in \N$ and satisfying $\sum_{k\in\N}k m_k= \rho_{\rm mi}$, and a probability measure $\psi$ on $\N_0$ having support included in $\{0, \dots, A\}$ for some $A\in\N_0$ with $\psi(A)\neq 0$ and satisfying $\sum_{a\in\N_0} a\psi(a)= \rho_{\rm ma}$. According to Step \ref{step:compactsupp}, it is sufficient to prove \eqref{lowbound1} for these $m$ and $\psi$. In the following, we consider the (more interesting) case that $\psi\not=\delta_0$, i.e., $\rho_{\rm ma}>0$; the remaining case $\rho_{\rm ma}=0$ needs some minor modifications, whose details we leave to the reader.

Introducing an auxiliary parameter $\delta\in(0,1)$, we will show for $\delta$ close to zero that
\begin{equation}\label{lowboundAKdelta}
\liminf_{N\to\infty}\frac 1{|\Lambda_N|}\log Z_{N,\L_N,{\rm Dir}}\geq -\varphi^{\ssup \delta}(m,\psi^{\ssup \delta}),
\end{equation}
where $\psi^\ssup\delta =\psi + 2\delta  \rho_{\rm mi} (\delta_{A+1}-\delta_A)$, and
\begin{equation}\label{varphi3}
\begin{aligned}
\varphi^{\ssup \delta}(m,\psi)&:=\inf\Big\{\sum_{a\in\N_0}\psi(a)\Big[I(P_a)+P_a(\Phi_{0,\Z^d})+2{\v}aP_a(N_0^{\ssup \ell})+{\v}a^2 \Big]\colon \\
&\qquad P_0,P_1,\dots\in\Mcal_1^{\ssup{\rm s}}(\Omega), \forall k\in\N\colon \sum_{a\in\N_0}\psi(a) P_a(N_{0}^{\ssup {\delta_k}})\in m_k(1-\delta,1+\delta)\Big\},
\end{aligned}
\end{equation}
where we interpret $0(1-\delta,1+\delta)$ as $\{0\}$.

To do this, we restrict to a configuration that has a predescribed microscopic part (depending on $m$), located well away from the boundary of $\L_N$, and additionally some macroscopic points at the origin (whose marks are contained in the box) and no other points. More precisely, we insert in $Z_{N,\L_N,{\rm Dir}}$ the indicator on the event that the configuration in $\Lambda_N$ has a microscopic random part with $N^{\ssup{\delta_k}}_{\L_N}\in m_k(1-\delta,1+\delta)|\L_N|$ for all $k\in [K]$, additionally the indicator $\1\{N_{\L_N\setminus \widetilde \L_N}^{\ssup\ell}=0\}$ on the event that a certain inner boundary  of $\L_N$ is empty of points, and we require that the remaining part of the configuration is equal to $\omega^{\ssup{\omega}}_\psi$ (to be defined later; it depends on $\omega^{\ssup{\leq K}}_{\L_N}$), that has only points at the origin and completes the entire configuration in such a way that it has precisely $N$ particles, all of which are contained in $\widetilde\L_N$. In other words, we estimate
$$
\begin{aligned}
\1\{N^{\ssup\ell}_{\Lambda_N}=N\}\1\{\widetilde N_{\L_N^{\rm{c}}}=0\}
&\geq \prod_{k\in[K]}\1\big\{N_{\widetilde \Lambda_N}^{\ssup{\delta_k}}\in m_k|\widetilde \L_N|(1-\delta,1+\delta)\big\}\\
&\qquad\times\1\{N^{\ssup\ell}_{\Lambda_N\setminus \widetilde\L_N}=0\}\1\{\omega^{\ssup{>K}}_{\widetilde \L_N}=\omega^{\ssup{\omega}}_\psi\},
\end{aligned}
$$
where $\widetilde\L_N$ is the largest centred box such that $x+G_k\subset\L_N$ for all $x\in \widetilde\L_N$ and $k\in [K]$. Clearly $|\widetilde\L_N|\sim|\L_N|$. Furthermore, we switch from $\E_{\L_N}$ to $\E_{\widetilde \L_N}^{\ssup{\leq K}}$ and note that $\E_{\L_N}=\E_{\widetilde \L_N}^{\ssup{\leq K}}\otimes \E_{\L_N\setminus\widetilde \L_N}\otimes \E_{\widetilde \L_N}^{\ssup{>K}}$ and that the expectation of the two last indicators is lower bounded by $\e^{-|\L_N\setminus \widetilde\L_N|\sum_k q_k}\e^{-|\L_N|\sum_{k>K}q_k+o(|\L_n|)}= \e^{-|\L_N|\sum_{k>K}q_k}\e^{o(|\L_N|)}$, since $\omega^{\ssup{\omega}}_\psi$ is macroscopic and $q_k=\e^{o(k)}$. This gives that
\begin{equation*}
\begin{split}
Z_{N,\L_N,{\rm Dir}}&\ge \e^{-|\L_N|\sum_{k>K}q_k}\e^{o(|\L_N|)}\\ &\qquad\times\E^{\ssup{\le K}}_{\widetilde\L_N}\Big[{\rm e}^{-\Phi_{\widetilde\Lambda_N,\widetilde\Lambda_N}(\cdot+\omega^{\ssup{\omega}}_\psi)}\prod_{k\in[K]}\1\big\{N_{\widetilde \Lambda_N}^{\ssup{\delta_k}}\in m_k|\widetilde \L_N|(1-\delta,1+\delta)\big\}\Big],
\end{split}
\end{equation*}
Here is the definition of $\omega^{\ssup{\omega}}_\psi$. We pick deterministic integers $N_0,N_1,\dots,N_{A-1}$ such that $N_a\le |\L_N|$ and $N_a\sim |\L_N|\psi([a,\infty))$, as $N\to\infty$, for all $a\in \{0,1,\dots,A-1\}$.
Then, we define the random integers $N_A=(N-\sum_{a\in[A-1]} N_a-\sum_{k\in[K]} k N_{\widetilde\Lambda_N}^{\ssup{\delta_k}})\wedge \lfloor |\L_N|\rfloor$ and $N_{A+1}=N-\sum_{a\in[A]} N_a-\sum_{k\in[K]} k N_{\widetilde\Lambda_N}^{\ssup{\delta_k}}$, and note that 
$$
-\delta\rho_{\rm ma}+\psi([A,\infty))+o(1)\le \frac{N_A}{|\L_N|}\le\big( \psi([A,\infty))+\delta\rho_{\rm ma}+o(1)\big)\wedge 1,
$$
and $0\le N_{A+1}/|\L_N|\le \delta\rho_{\rm ma}+o(1)$. From now on we assume that $0< \delta\le \psi([A,\infty))/4\rho_{\rm ma}$, which guarantees that $N_A\ge 0$. We put $\omega^{\ssup{\omega}}_\psi=\sum_{a\in[A+1]\colon N_a\neq 0}\delta_{(0,G_{N_a})}$, i.e., we put at the origin the marks $G_{N_1},\dots,G_{N_{A+1}}$ on top of each other. This is a macroscopic configuration (with a possible exception of the mark $G_{N_{A+1}}$) that satisfies the zero-Dirichlet boundary conditions since $N_a\le |\L_N|$ for all $a\in [A+1]$. Note that, for any $a\in\N_0$, the number of sites in $\L_N$ that carry precisely $a$ particles of $\omega^{\ssup{\omega}}_\psi$ is  $\psi(a)|\L_N |(1+o(1))$. Furthermore,  the total number of particles in the configuration $\omega+\omega^{\ssup{\omega}}_\psi$ is equal to $N$. Note that $\omega^{\ssup{\omega}}_\psi$ depends on $\omega$ and is therefore random.

Now we construct a deterministic macroscopic marked point configuration $\omega_{\psi^{\ssup\delta}}$ in $\L_N$ such that the number of sites in $\L_N$ that are covered by precisely $a$ particles from this configuration is $\sim |\L_N|\psi^{\ssup\delta}(a)$ as $N\to\infty$, for any $a\in\N_0$. Indeed, we replace $N_A$ and $N_{A+1}$ by the deterministic values $\widetilde N_A=\lfloor |\widetilde\L_N|[(\psi([A,\infty))+2\delta\rho_{\rm ma})\wedge 1]\rfloor$ and $\widetilde N_{A+1}=\lfloor 2\delta \rho_{\rm ma}|\L_N|\rfloor$ and denote by $\omega_{\psi^{\ssup\delta}}$ the corresponding macroscopic configuration, with $\psi^{\ssup{\delta}}$ defined below \eqref{lowboundAKdelta}. Since $N_A\leq\widetilde N_A$ and $N_{A+1}\leq\widetilde N_{A+1}$, we have
$$
\Phi_{\widetilde\Lambda_N,\widetilde\Lambda_N}(\omega+\omega^{\ssup{\omega}}_\psi)\leq \Phi_{\widetilde\Lambda_N,\widetilde\Lambda_N}(\omega+\omega_{\psi^{\ssup \delta}})
$$
and we can lower bound: 
\begin{equation}\label{lowbound2}
Z_{N,\L_N,{\rm Dir}}\geq \e^{-|\L_N|\sum_{k>K}q_k}\e^{o(|\L_N|)} Z_{N,\widetilde\L_N}^{\ssup\delta}(m,\psi^\ssup\delta),
\end{equation}
where 
\begin{equation}\label{ZNdeltaKAdef}
Z_{N,\widetilde\L_N}^{\ssup \delta}(m,\psi)=\E^{\ssup{\le K}}_{\widetilde\L_N}\Big[\e^{-\Phi_{\widetilde\L_N,\widetilde\L_N}(\cdot+\omega^{\ssup{\omega}}_{\psi})}\prod_{k\in[K]}\1\big\{N_{\widetilde\Lambda_N}^{\ssup{\delta_k}}\in m_k(1-\delta,1+\delta)|\widetilde\L_N|\big\}\Big].
\end{equation}
Since $\psi^{\ssup\delta}\geq_{\rm st}\psi$, there are more that $N$ particles in the expectation on the right-hand side of \eqref{ZNdeltaKAdef} for $\psi^{\ssup\delta}$ in place of $\psi$, but this is no problem at all. 
In a small abuse of notation, we write  from now $\L_N$ instead of $\widetilde\L_N$.

Recall that we picked $m=(m_k)_{k\in\N}$ with support included in $[K]=\{1,\dots,K\}$ for some $K\in \N$, and a probability measure $\psi$ on $\N_0$ having support included in $\{0, \dots, A\}$ for some $A\in\N_0$. Steps \ref{step-decompLN} to \ref{step-QtoZd} are devoted to showing the following result: for any $\delta\in(0,1)$, for any $m$ and $\psi$ having finite support,
\begin{equation}\label{lowboundZdeltavarphidelta}
\lim_{N\to \infty} \frac{1}{|\L_N|}\log Z_{N,\L_N}^{\ssup\delta}(m,\psi_N)\geq -\varphi^{\ssup{\delta}}(m,\psi),\qquad\delta\in (0,1),
\end{equation}
with $\varphi^{\ssup{\delta}}(m,\psi)$ defined in \eqref{varphi3}, and any sequence $(\psi_N)_{N\in\N}$ of measures on $\N_0$ such that $\psi_N([a,\infty))\to\psi([a,\infty))$ as $N\to \infty$, and $\psi_N([A+1,\infty))=0$ for all $N$. We will conclude the proof by applying in \eqref{lowbound2} this result with $(m, \psi^{\ssup \delta})$ instead of $(m, \psi)$ and, in Step \ref{step-lowbounddeltato0}, by taking $\delta$ to 0.

\begin{step}\label{step-decompLN}
We decompose $\L_N$ in $R$-boxes.
\end{step}

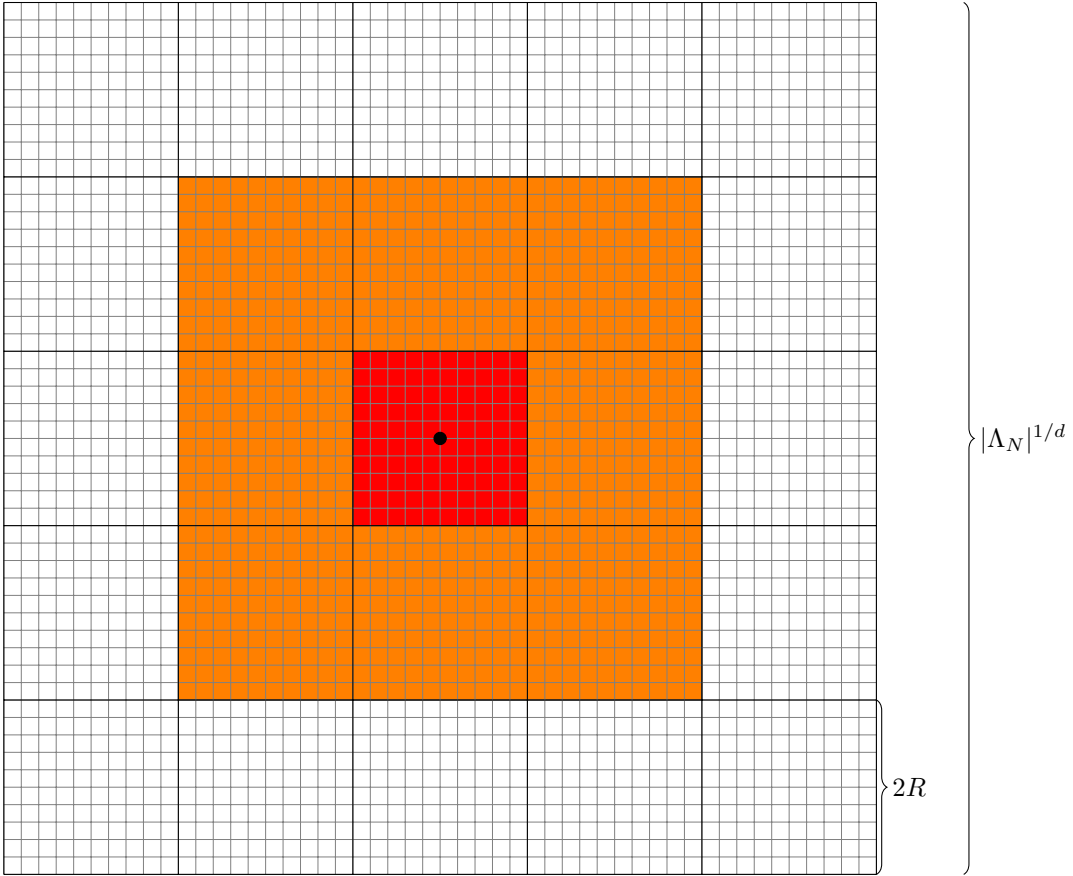
\begin{figure}[!htpb]
\centering
\input{Pix_LowerBound}
	\caption{Illustration for the construction of macroscopic marks. We attach two Poisson points at the origin with the desired marks indicated by the red and orange boxes. In this case, we arrive at $\psi(0)=16/25$, $\psi(1)=8/25$ and $\psi(2)=1/25$.}
	\label{Pix_LowerBound}
\end{figure}

Introduce a new large auxiliary parameter $R\in\N$ and consider the box $Q=[-R,R)^d\cap\Z^d$. We denote by $Y_N=Y_N(R)$ the set of all $z\in 2R\Z^d$ such that $z+Q\subset \L_N$. We decompose the box $\Lambda_N$ into the boxes $Q_z=z+Q$ with $z\in Y_N$. We call the boxes $Q_z$ sometimes {\em mesoboxes}. We may  assume that $\L_N$ is equal to the union of these boxes, since otherwise we insert the indicator on the event that the configuration $\omega_{\rm P}$ has no particles in the difference between $\L_N$ and the union (call it $\widetilde\L_N$), then we can replace the interaction in $\L_N$ by the interaction in $\widetilde \L_N$ and can separate the entire expectation into the probability that the difference is empty and the same expectation with $\L_N$ replaced by $\widetilde \L_N$. It is easy to see that the exponential rate of the former vanishes as $N\to\infty$ for any $R$. Hence, we assume from now that $\L_N$ is equal to the union of the boxes $Q_z$.
Observe that for each $a\in\{0,\dots, A\}$, the number of boxes $Q_z$ that are hit by precisely $a$ of the macroscopic marks is deterministic and approximately equal to $\psi_N(a)|Y_N|$, where we already note that $|Y_N|=|\Lambda_N|(2R)^{-d}$.
It is guaranteed that
\begin{equation}\label{Achoice}
\lim_{N\uparrow\infty}\frac1{|Y_N|}\big|\{z\in Y_{N}\colon Q_z\mbox{ is hit by precisely }a\mbox{ marks}\}\big|=\psi(a),\qquad a\in \{0,\dots, A\},
\end{equation}
see Figure~\ref{Pix_LowerBound} for an illustration.

\begin{step}
We isolate the $Q_z$'s.
\end{step}

In this step, we further lower bound $Z_{N,\L_N}^{\ssup\delta}(m,\psi_N)$ in such a way that there is no mutual interaction between any two different mesoboxes $Q_z$. For this, we insert an indicator on the event that all the Poisson points are located sufficiently far away from the boundaries of the $Q_z$. For this sake, we recall that the interaction potential $v$ is assumed to have  bounded support, and pick $L$ such  that $v(r)=0$ for any $r>L$. Pick the box $\widetilde Q=[-R+S,R-S)^d\cap\Z^d$ with some $S>L+\max_{k\in[K]}\diam(G_k)$ and put $\widetilde Q_z=z+\widetilde Q$, assuming that $R>4S$. We require that the PPP has no particles in the region $z+(Q\setminus \widetilde Q)=Q_z\setminus \widetilde Q_z$. Then any two particles in microscopic marks in different $Q_z$'s have no interaction with each other, and a microscopic particle of a point in $Q_z$ has no interaction  with any macroscopic particle that lies in a different $Q_z$. 
Then the total energy is now equal to the sum over $z$ of the self-energy of the microscopic marks in $Q_z$ plus their energy with those macroscopic particles that lie in the same $Q_z$ plus the self-energy of the latter particles, plus the energy between the macroscopic particles in different $Q_z$'s. The total energy coming from macroscopic particles in different boxes is bounded by some constant (depending on $v$, $d$ and $A$) times the number of points in $\Lambda_N$ whose distance to the boundary their box $Q_z$ is smaller than the diameter of the support of $v$. This region has a volume that is not larger than a constant times $|\Lambda_N|/R$.

We write $a^{\ssup z}\in\{0,1,\dots,A\}$ for the number of macroscopic marks that hit $Q_z$. Note that the cross-energy between a particle in a microscopic mark with all the macroscopic particles of one such grid is not greater than the number ${\v}$ defined in \eqref{Cvdef}; hence this part of the energy is equal to ${\v}$ times the number of microscopic particles in $Q_z$, i.e., times $N_{Q_z}^{\ssup{\ell}}$. Hence, we estimate, on the above indicator,
$$
\begin{aligned}
\Phi_{\Lambda_N,\Lambda_N}&\big(\omega,\psi_N\big)\leq \sum_{z\in Y_{N}} \widetilde \Phi_{\widetilde Q}(\omega^{\ssup z}, a^{\ssup z})+ C |\Lambda_N|/R,
\end{aligned}
$$
where $C$ is a constant that depends on $d$, $v$ and $A$ only, and we introduced the shifted restriction $\omega^{\ssup z}=\theta_{-z}(\omega_{\widetilde Q_z})\in\Omega_{\widetilde Q}$, where we recall that $\theta_z$ is the shift operator by $z$, and $\Omega_{\widetilde Q}$ denotes the set of all marked point processes with points in $\widetilde Q$ and marks in $\{G_k\colon k\in\N\}$. Also, we introduced the energy functional
\begin{equation}
\widetilde \Phi_{\widetilde Q}(\omega,a)=\Phi_{\widetilde Q,\widetilde Q}(\omega)+2{\v}a N^{\ssup\ell}_{\widetilde Q}(\omega)+|Q|\v_Q a^2, \qquad \omega\in \Omega_{\widetilde Q},a\in\N_0, 
\end{equation}
where  and $\v_{Q}=\frac 1{|Q|}\sum_{i,i'\in Q}v(|i-i'|)$ is a constant that depends only on $v$ and $R$ and we note that it converges towards ${\v}$ as $R\to\infty$.

This gives the estimate, for any sufficiently large $N$, 
\begin{equation}\label{lowbound3}
\begin{aligned}
Z_{N,\L_N}^{\ssup\delta}(m,\psi_N)&\geq \e^{-\frac CR |\Lambda_N|}\E^{\ssup{\leq K}}_{\L_N}\Big[\e^{-\sum_{z\in Y_{N}}\widetilde\Phi_{\widetilde Q}(\cdot^{\ssup z},a^{\ssup z})}\\
&\qquad \times\Big[\prod_{k\in[K]}\1\big\{N_{\Lambda_N}^{\ssup{\delta_k}}\in |\Lambda_N|m_k(1-\delta,1+\delta)\big\}\Big]\Big[\prod_{z\in Y_{N}} \1\big\{N_{Q_z\setminus\widetilde Q_z)}=0\big\}\Big]\Big].
\end{aligned}
\end{equation}
We can write the event $\{N_{\Lambda_N}^{\ssup{\delta_k}}\in m_k|\Lambda_N|(1-\delta,1+\delta)\}$ as the event $\{\sum_{z} N^{\ssup{\delta_k}}_{\widetilde Q}(\omega_z)\in |\Lambda_N|m_k(1-\delta,1+\delta)\}$, which is independent of the void probability in the end. Also the interaction term is independent of the void  probability. The latter is not smaller than $\prod_z\e^{-|Q\setminus \widetilde Q|\sum_{k\in\N} q_k}\geq \e^{-C|\Lambda_N|/R}$ for any $N$ and $R$, if $C$ is sufficiently large, depending only on $v$, $d$, $K$ and $\sum_k q_k$.

\begin{step}
We rewrite the expectation in terms of an expectation with respect to  a crucial empirical measure.
\end{step}

Now we introduce an empirical measure that we will need for our large-deviations arguments:
$$
\widetilde \eta_{N,R}=\frac 1 {|Y_N|}\sum_{z\in Y_{N}}\delta_{(\omega^{\ssup z},a^{\ssup z})}\in\Mcal_1(\Omega_{\widetilde Q}\times \{0,\dots,A\}).
$$
Then we have
\begin{equation}\label{lowbound4}
\begin{aligned}
Z_{N,\L_N}^{\ssup\delta}(m,\psi_N)&\geq \e^{-2\frac CR |\Lambda_N|}\\
&\times \E^{\ssup{\leq K}}_{\L_N}\Big[\e^{-|Y_N|\langle \widetilde\eta_{N,R},\widetilde\Phi_{\widetilde Q}\rangle}\prod_{k\in[K]}\1\Big\{\Big\langle\widetilde \eta_{N,R},\frac1{|Q|} N_{\widetilde Q}^{\ssup{\delta_k}}\Big\rangle\in m_k(1-\delta,1+\delta)\Big\}
\Big],
\end{aligned}
\end{equation}
where we wrote short $N_{\widetilde Q}^{\ssup{\delta_k}}$ for the map $(\omega,a)\mapsto N_{\widetilde Q}^{\ssup{\delta_k}}(\omega)$. Note that, by \eqref{Achoice}, $\widetilde \eta_{N,R}(\Omega_{\widetilde Q}\times\{a\})$ is deterministic and converges to $\psi(a)$.

\begin{step}
We carry out the large-$N$ asymptotics with the help of large-deviation arguments. 
\end{step}

Hence, the application of a multi-type variant of Sanov's theorem  yields that $(\widetilde \eta_{N,R})_N$ satisfies a large-deviation principle as $N\to\infty$ on $\Mcal_1(\Omega^{\ssup{\leq K}}_{\widetilde Q}\times \{0,\dots,A\})$ with scale $|\Lambda_N|$ and rate function
$$
\eta\mapsto  I_{Q,\widetilde Q}(\eta)=
\frac 1{|Q|}\sum_{a=0}^{A}\psi(a)H_{\widetilde Q}(\eta_a|\P^{\ssup{\leq K}}_{\widetilde Q}),
$$
if $\eta(\Omega^{\ssup{\leq K}}_{\widetilde Q}\times\{a\})=\psi(a)$ for all $a\in [A]$, and infinity otherwise. Here, $\eta_a(\cdot)=\eta(\cdot\times\{a\})/\eta(\Omega^{\ssup{\leq K}}_{\widetilde Q}\times\{a\})\in \Mcal_1(\Omega^{\ssup{\leq K}}_{\widetilde Q})$ is the conditional distribution on $\Omega^{\ssup{\leq K}}_{\widetilde Q}$ given the mark $a$.

Note that the state space $\Omega^{\ssup{\leq K}}_{\widetilde Q}\times \{0, \dots, A\}$ is countable, and on $\Mcal_1(\Omega^{\ssup{\leq K}}_{\widetilde Q}\times \{0,\dots,A\})$ we use the product topology. Note that the maps $\eta\mapsto \langle \eta,\widetilde\Phi_{\widetilde Q}\rangle$ and $\eta\mapsto \langle \eta,\frac 1{|Q|}N_{\widetilde Q}^{\ssup {\delta_k}}\rangle$ are continuous and also bounded on the event that we consider in \eqref{lowbound3}. 

Hence, recalling \eqref{lowbound2} and  using Varadhan's lemma, we obtain that
\begin{equation}\label{Varadhan1}
\liminf_{N\to \infty}\frac 1 {|\L_N|}\log Z^{\ssup{\delta}}_{N,\Lambda_N}(m,\psi_N)\geq -\frac {2C}R- \varphi^{\ssup{\delta, K}}_{Q,\widetilde Q}(m,\psi),
\end{equation}
where we recall that the support of $m$ is contained in $[K]$ and we define
\begin{equation}\label{varphi2}
\begin{aligned}
\varphi_{Q,\widetilde Q}^{\ssup{\delta, L}}(m,\psi)=\frac1{|Q|}\inf\Big\{\sum_{a\in \N_0}\psi(a)&\big[H_{\widetilde Q}(\eta_a|\P^{\ssup {\leq L}}_{\widetilde Q})+\langle \eta_a,\widetilde \Phi_{\widetilde Q}(\cdot,a)\rangle\big]\colon \eta_0,\eta_1,\dots\in\Mcal_1(\Omega^{\ssup {\leq L}}_{\widetilde Q}),\\
&\forall k\in[L]\colon \sum_{a\in \N_0}\psi(a)\big\langle \eta_a,\smfrac 1{|Q|}N_{\widetilde Q}^{\ssup {\delta_k}}\big\rangle\in m_k(1-\delta,1+\delta)\Big\}.
\end{aligned}
\end{equation}

\begin{step}\label{step-QtoZd}
We let $Q\uparrow \Z^d$.
\end{step}
Note that $\varphi_{Q,\widetilde Q}^{\ssup{\delta, K}}(m,\psi)=\varphi_{Q,\widetilde Q}^{\ssup{\delta, \infty}}(m,\psi)+\frac{|\widetilde Q|}{|Q|}\sum_{k> K}q_k$. 
Now we consider the limit as $R\to\infty$ and want to show that $\limsup_{R\to\infty}\varphi_{Q,\widetilde Q}^{\ssup{\delta,\infty}}(m,\psi)\leq \varphi^{\ssup{\delta}}(m,\psi)$. 
Recall that $Q=[-R,R)^d\cap \Z^d$ and $\widetilde Q=[-R+S,R-S)^d\cap\Z^d$ for some $S>0$, and hence $|\widetilde Q|/|Q|\to 1$ as $R\to\infty$. The main idea is to restrict the infimum in \eqref{varphi2} to those $\eta_a$ that are the restriction to $\widetilde Q$ (more precisely, the projection on $\Mcal_1(\Omega^{\ssup {\leq K}}_{\widetilde Q})$) of some $P_a\in\Mcal_1^{\ssup{\rm s}}(\Omega)$ satisfying $P_a(\Phi_{0,\Z^d})<\infty$, which leads to an upper estimate for $\varphi_{Q,\widetilde Q}^{\ssup{\delta,\infty}}(m,\psi)$. Hence, pick some $(P_a)_{a}$ that are admissible in the definition \eqref{varphi3} and put $\eta_a=(P_a)_{\widetilde Q}=P_a\circ \pi_{\widetilde Q,K}^{-1}$ for $a\in\N$, where $\pi_{\widetilde Q,K}\colon\Omega\to\Omega^{\ssup {\leq K}}_{\widetilde Q}$ is the canonical projection. Without loss of generality we conceive the $\eta_a$ as elements of $\Mcal_1(\Omega_{\widetilde Q})$ without mass on configurations with marks of size $>K$. In particular, $(\eta_a)_a$ then satisfies the last condition in \eqref{varphi2} for all $k> K$ since $m_k=0$ for these $k$. Moreover, $(\eta_a)_a$ is admissible on the right-hand side of \eqref{varphi2} for all sufficiently large $R$, since
$$
\big\langle \eta_a,\smfrac 1{|Q|}N_{\widetilde Q}^{\ssup {\delta_k}}\big\rangle=  \frac{|\widetilde Q|}{|Q|}\langle P_a,N_0^{\ssup{\delta_k}}\rangle \to \langle P_a,N_0^{\ssup{\delta_k}}\rangle,\quad\mbox{as }R\to\infty,\qquad a\in\N_0, k\in[K],
$$
by stationarity of $P_a$, where we wrote $N_{0}^{\ssup {\delta_k}}$ also for the map $\Omega\to\R$.

We will show that plugging this family $(\eta_a)_a$ into the functional on the right-hand side of \eqref{varphi2} gives, in the limit as $R\to\infty$, a value that is not larger than the value of the functional on the right-hand side of \eqref{varphi3} for $(P_a)_a$. Minimising over $(P_a)_a$ then gives the desired result.

We first argue that
\begin{equation}\label{ergothm}
\frac 1{|Q|}\langle\eta_a,\widetilde \Phi_{\widetilde Q}(\cdot,a)\rangle 
\leq\langle P_a,\Phi_{0,\Z^d}\rangle+2{\v}a\langle P_a,N_0^{\ssup \ell}\rangle+{\v} a^2 ,\quad a\in\N_0, R>0.
\end{equation} 
Indeed, the two last terms are easily understood, using stationarity and the above remark that $|\widetilde Q|/|Q|\to 1$; recall also that $\v_{Q}\to {\v}$ as $R\to\infty$. For understanding the first term in \eqref{ergothm}, we note that, since $v\ge0$ and hence $\Phi_{\widetilde Q,\widetilde Q}\le \Phi_{Q,\Z^d}$,
$$
\frac 1{|Q|}\langle \eta_a,\Phi_{\widetilde Q,\widetilde Q}\rangle\leq 
\frac 1{|Q|}\sum_{x\in Q} \langle P_a, \theta_x\Phi_{0,\Z^d}\rangle=\langle P_a,\Phi_{0,\Z^d}\rangle,\quad R>0,
$$
by stationarity, where  $\theta_x$ is the shift operator by $x\in\Z^d$. This explains also the first term in \eqref{ergothm}. 

Now we turn to the entropic term and recall \eqref{Igriddef}, which gives us that
\begin{equation}\label{entropRtoinfty}
\lim_{R\to\infty}\frac1{|Q|}H_{\widetilde Q}(\eta_a|\P_{\widetilde Q})=I(P_a),\qquad a\in\N_0.
\end{equation}
This, together with  \eqref{ergothm} gives that
$$
\limsup_{R\to\infty}\varphi^{\ssup{\delta,\infty}}_{Q,\widetilde Q}(m,\psi)\leq \sum_{a\in\N_0}\psi(a)\Big[I(P_a)+P_a(\Phi_{0,\Z^d})+2{\v}a \langle P_a,N_0^{\ssup \ell}\rangle+ {\v}a^2\Big].
$$
By minimising over $(P_a)_a$, the right-hand side approaches $\varphi^{\ssup\delta}(m,\psi)$. Hence, using \eqref{Varadhan1}, we have proved \eqref{lowboundZdeltavarphidelta}. Together with \eqref{lowbound2}, this yields \eqref{lowboundAKdelta}.

\begin{step}\label{step-lowbounddeltato0}
We take $\delta$ to zero.
\end{step}

Recall that \eqref{lowboundAKdelta} is in terms of $\psi^{\ssup\delta}$  instead of $\psi$, where we recall that $\psi^{\ssup\delta}=\psi+2\delta\rho_{\rm ma}(\delta_{A+1}-\delta_{A})$. To finish the proof of \eqref{lowbound1}, we have to prove that $\limsup_{\delta\downarrow0}\varphi^{\ssup\delta}(m,\psi^{\ssup\delta})\leq \varphi(m,\psi)$ for $m$ with support in $[K]$ and $\psi$ with support in $\{0,\dots,A\}$. Note first that, trivially, $\varphi^{\ssup\delta}\leq \varphi$. Now we show that $\limsup_{\delta\downarrow0}\varphi(m,\psi^{\ssup\delta})\leq \varphi(m,\psi)$. Pick $P_0,P_1,\dots,P_A\in\Mcal_1^{\ssup{\rm s}}(\Omega)$ admissible in the formula \eqref{varphidef} for $\varphi(m,\psi)$. Now  put $P_{A+1}=P_A$, then $(P_a)_{a\in\{0,\dots,A+1\}}$ is admissible in the formula \eqref{varphidef} for $\varphi(m,\psi^{\ssup{\delta}})$, since $\sum_a\psi^{\ssup\delta}(a)P_a(N_0^{\ssup{\delta_k}})=m_k+2\delta\rho_{\rm ma}[P_{A+1}(N_0^{\ssup{\delta_k}})-P_A(N_0^{\ssup{\delta_k}})]=m_k$ for any $k\in\N$. Inserting this in \eqref{varphidef} for $\varphi(m,\psi^{\ssup{\delta}})$ gives (abbreviating $\Xi_a(P)=I(P)+P(\Phi_{0,\Z^d})+2{\v}aP(N_0^{\ssup \ell})+{\v}a^2$),
$$
\varphi(m,\psi^{\ssup{\delta}})\leq\sum_{a}\psi^{\ssup\delta}(a)\Xi_a(P_a)=\sum_{a}\psi(a)\Xi_a(P_a)+2\delta\rho_{\rm ma}[\Xi_{A+1}(P_A)-\Xi_A(P_A)].
$$
Taking $\delta$ to 0 and then minimising over $(P_a)_a$ finishes the proof of \eqref{lowbound1} and hence the proof of the lower bound in Theorem \ref{thm-freeenergygrid}.

\setcounter{step}{0}

\subsection{Proof of the upper bound in Theorem~\ref{thm-freeenergygrid}}\label{sec-uppbound}

Now we turn to the proof of \eqref{ProofThm1}.

\begin{step}
We estimate against one maximal cardinality configuration.
\end{step}

Recall that  $\rho|\Lambda_N|=N$ and decompose according to the numbers of Poisson points for any mark size:
\begin{equation}
Z_{N,\L_N,{\rm Dir}}
=\sum_{l=(l_k)_{k\in\N}\in\N_0^{\N}\colon \sum_k k l_k=\rho|\Lambda_N|}Z_{N,\L_N,{\rm Dir}}(l)
\end{equation}
with 
$$
Z_{N,\L_N,{\rm Dir}}(l)=\E_{\L_N}\Big[\e^{-\Phi_{\Lambda_N,\Lambda_N}}\1\{\widetilde N_{\L_N^{\rm c}}=0\}\prod_{k\in\N}\1\{N_{\Lambda_N}^{\ssup{\delta_k}}=l_k\}\Big].
$$
Now we take, for any $N$, one sequence $l^{\ssup{*,N}}=(l_k^{\ssup{*,N}})_{k\in\N}$ that is maximal for $l\mapsto Z_{N,\L_N,{\rm Dir}}(l)$ under the constraint $\sum_k k l_k=\rho|\Lambda_N|$, and estimate
\begin{equation}\label{ZNgegenZN(l)}
\begin{aligned}
Z_{N,\Lambda_N,{\rm Dir}}&\leq \#\Big\{(l_k)_{k\in[N]}\in\N_0^{[N]}\colon \sum_k kl_k=N\Big\} Z_{N,\L_N,{\rm Dir}}(l^{\ssup{*,N}}).
\end{aligned}
\end{equation}
It is known that the counting term is not larger than $\e^{o(N)}$. It is clear that we can find a subsequence along which, for some sequence $(m^*_k)_{k\in\N}$, we have that $\frac 1{|\Lambda_N|}l_k^{\ssup{*,N}}\to m^*_k$ for any $k\in\N$. According to Fatou's lemma, we have that
\begin{equation}\label{rho1rho2}
 \rho_{\rm mi}=\sum_{k\in\N} k m^*_k\in[0,\rho]\mbox{, and we put }  \rho_{\rm ma}=\rho- \rho_{\rm mi}\in[0,\rho].
\end{equation}
The remainder of the proof shows that
\begin{equation}\label{upperbound}
\limsup_{N\to\infty}\frac1{|\L_N|}\log Z_{N,\L_N,{\rm Dir}}(l^{\ssup N})\leq -\inf_{\psi\in\Mcal_1(\N_0)\colon \sum_{a\in\N_0}a \psi(a)= \rho_{\rm ma}}\varphi(m,\psi),
\end{equation}
for any $m\in[0,\infty)^\N$ satisfying $\sum_{k\in\N} k m_k= \rho_{\rm mi}$ and any sequence $l^{\ssup N}$ in $\N_0^\N$ satisfying $\sum_k kl_k^{\ssup{N}}=N$ for all $N\in\N$, and $\frac 1{|\Lambda_N|}l_k^{\ssup{N}}\to m_k$ for any $k\in\N$.

\begin{step}
We integrate out macroscopic marks and decompose $\L_N$ in $R$-boxes.
\end{step}
We introduce a large parameter $K\in \N$ and introduce the cut-off version of $m$ by putting $m^{\ssup{\leq K}}=(m_1,m_2,\dots,m_K,0,0,\dots)$.  We call the marks $G_k$ with $k\leq K$ {\it microscopic} and the others {\it macroscopic}. We write now the expectation over the configuration as an integration over the part of the configuration that has only microscopic marks, and an explicit integral for the location of any macroscopic mark. Hereby, we use that, conditional on $\{N_{\Lambda_N}^{\ssup{\delta_k}}=l_k\}$, the $l_k$ Poisson points are independent and uniformly distributed over $\L_N$, for any $k>K$ and any $l_k\in\N$. Estimating these Poisson probabilities against one, this gives
\begin{equation}\label{Upperboundcitelater}
\begin{aligned}
Z_{N,\L_N,{\rm Dir}}(l^{\ssup{N}})&\leq\E_{\L_N}^{\ssup{\leq K}}\otimes\bigotimes_{k>K}\bigotimes_{j=1}^{l_k^{\ssup{N}}}
\Ucal_{\Lambda_N}\Big[\e^{-\Phi_{\Lambda_N,\Lambda_N}(\cdot+\omega_X)}\1\{\widetilde N_{\L_N^{\rm c}}=0\}\prod_{k=1}^K\1\{N_{\Lambda_N}^{\ssup{\delta_k}}=l_k^{\ssup{N}}\}\Big],
\end{aligned}
\end{equation}
where $X=(X_{k,j})_{k>K; j=1,\dots,l_k^{\ssup{N}}}$ is a collection of independent and uniformly over $\Lambda_N$ distributed random sites with distribution (and expectation) $\Ucal_{\L_N}$ and $\omega_X=\sum_{k>K}\sum_{j=1}^{l_k^{\ssup{N}}}\delta_{(X_{k,j},G_k)}$ is the superposition of the marked points. 

Let us fix a large auxiliary parameter $R\in\N$ and decompose $\Lambda_N$ regularly into auxiliary boxes $Q_z=z+Q=z+[-R,R)^d\cap\Z^d$, of radius $R$ centred at $z\in Y_N=Y_{N,R}=\{z\in 2R\Z^d \colon Q_z\subset \L_N\}$. If $|\Lambda_N|$ is not a multiple of $|Q|=(2R)^d$, then cut the overshoot away, obtaining a box $\widetilde \Lambda_N$, which is precisely equal to the disjoint union of the $Q_z$ with $z\in Y_{N}$. There are $|Y_N|\sim |\Lambda_N|/(2R)^d$ such boxes. The volume of the amended box satisfies $1-C R N^{-1/d}\leq|\widetilde\Lambda_N|/|\Lambda_N|\leq 1$ for some $C$ that depends only on $d$. 
Additionally, we pick a small $\delta\in(0,1)$ and estimate, for any $k\in[K]$ and any sufficiently large $N$,
$$
\1\{N_{\Lambda_N}^{\ssup{\delta_k}}=l_k^{\ssup{N}}\}
\leq \1\{N_{\widetilde\Lambda_N}^{\ssup{\delta_k}}\in m_k |\Lambda_N|[1-\delta,1+\delta]\}+\1\{N_{\Lambda_N\setminus \widetilde \Lambda_N}^{\ssup{\delta_k}}> \delta |\Lambda_N|\}.
$$
The expectation of the latter term is easily shown to have an exponential rate equal to $-\infty$ on the scale $|\L_N|$ for any $\delta>0$, using the exponential Chebyshev inequality and the fact that $N_{\Lambda_N\setminus \widetilde \Lambda_N}^{\ssup{\delta_k}}$ is Poisson-distributed with parameter of surface order of $|\L_N|$. So far, we have that, with some arbitrarily large $C>0$, for any $N$ large enough,
\begin{equation}\label{uppbounddeltaappears}
\begin{aligned}
&Z_{N,\L_N,{\rm Dir}}(l^{\ssup{N}})\leq \e^{-C|\Lambda_N|}\\
&+\E_{\L_N}^{\ssup{\leq K}}\otimes\bigotimes_{k>K}\bigotimes_{j=1}^{l_k^{\ssup{N}}}
\Ucal_{\L_n}\Big[\e^{-\Phi_{\widetilde\Lambda_N,\widetilde\Lambda_N}(\cdot+\omega_X)}\prod_{k=1}^K\1\{N_{\widetilde\Lambda_N}^{\ssup{\delta_k}}\in m_k|\Lambda_N|[1-\delta,1+\delta]\}\Big].
\end{aligned}
\end{equation}
We also used that $\Phi_{\Lambda_N,\Lambda_N}\geq \Phi_{\widetilde\Lambda_N,\widetilde\Lambda_N}$.

\begin{step}
We cut off overshoot.
\end{step}

We also would like to have that each macroscopic mark $G_{k}$ centred at $X_{k,j}$ for $k>K$ either covers any $Q_z$ entirely or does not intersect it. To achieve this, from each of the macroscopic marks $G_k$ centred at $X_{k,j}$ cut away the overshoot of the largest union of the $Q_z$'s that lie in the mark. Also this amendment does not increase the interaction. We are cutting away only a negligible amount of particles, i.e., a number of particles that is $\leq |\L_N|\delta_{K,R}$ with some $\delta_{K,R}$ that vanishes as $K\to\infty$, for fixed $R$. In order to see this, recall that $G_k$ approaches a regular discrete box of cardinality $\approx k$ and note that therefore the number of particles in a mark $G_k$ with $k>K$ that we cut away is at most
$$
\lceil k^{1/d}\rceil^d-\big(\lceil k^{1/d}\rceil-4R\big)^d
=\lceil k^{1/d}\rceil^d\Big[1-\Big(1-\frac{4R}{\lceil k^{1/d}\rceil}\Big)^d\Big]\leq \lceil k^{1/d}\rceil^d 2d \frac{4R}{\lceil k^{1/d}\rceil}\leq  k \delta_{K,R},
$$
with some $\delta_{K,R}$ as announced. We denote the amended mark centred at $X_{k,j }$ by $\widetilde G_{k,j}$.

Let $ a^{\ssup z}\in\N_0$ denote the number of macroscopic Poisson points $X_{k,j}$ with $k>K$ and $j\in[l_k^{\ssup{N}}]$ such that its mark $\widetilde G_{k,j}$ centred at $X_{k,j}$ contains $ Q_z$. We then say that $Q_z$ is of type $ a^{\ssup z}$. For any $(k,j)$ with $k>K$ and $j\in[l_k^{\ssup{N}}]$, there are $|Q|^{-1}|\widetilde G_{k,j}|$  boxes $Q_z$ that are covered by the amended mark $X_{k,j}+\widetilde G_{k,j}$. Hence, 
\begin{equation}\label{countingmacro}
\sum_{z\in Y_{N}} a^{\ssup z}=\sum_{k>K}\sum_{j=1}^{l_k^{\ssup{N}}}|Q|^{-1}|\widetilde G_{k,j}|\in |Q|^{-1}\sum_{k>K}\sum_{j=1}^{l_k^{\ssup{N}}}k\big[1-\delta_{K,R},1\big]=|Q|^{-1} N^{\ssup{\rm Ma}}[1-\delta_{K,R},1],
\end{equation}
where we write 
\begin{equation}\label{nzweidef}
N^{\ssup{\rm Ma}}=N-\sum_{k\in[K]} k l_k^{\ssup{N}}
\end{equation}
for the total number of particles in all the macroscopic marks. Recall  \eqref{rho1rho2} to see that, for any sufficiently large $N$, 
\begin{equation}\label{NMaasy}
N^{\ssup{\rm Ma}}\in \rho |\L_N|-|\Lambda_N|\Big(\sum_{k\in[K]} k m_k\Big)[1-\delta, 1+\delta]=|\Lambda_N|\Big( \rho_{\rm ma}+\sum_{k>K} k m_k+\Big(\sum_{k\in [K]}km_k\Big)[-\delta,\delta]\Big).
\end{equation}

\begin{step}
We drop all interaction between any two distinct $Q_z$'s. 
\end{step}

More precisely, for any $z$, we attach to $Q_z$ all interaction that (1) the microscopic marks at points $\in Q_z$ have with each other or (2) the microscopic marks $\subset Q_z$ with any macroscopic particle $\in Q_z$ (these appear twice) and (3) the interaction that all the macroscopic particles in $\in Q_z$ have with each other. The sum of these three interactions is not smaller than $\widetilde \Phi_{Q}(\omega^{\ssup z}, a^{\ssup z})$, where 
$$
\omega^{\ssup z}=\theta_{-z}(\omega_{Q_z})\in\Omega_Q^{\ssup {\leq K}},\qquad \omega\in \Omega_Q^{\ssup {\leq K}},
$$
is the restriction of $\omega\in\Omega$ to $Q_z$, shifted to $Q$, and 
\begin{equation}\label{energydecompose}
\widetilde \Phi_{Q}(\omega, a)=
\Phi_{Q,Q}(\omega)
+2 \v a N_{\widetilde Q}^{\ssup\ell}(\omega)+ |Q| \v_{Q} a^2 ,\qquad  a\in\N_0,\omega\in\Omega_Q^{\ssup{\leq K}},
\end{equation}
where we introduced $\v_Q= \frac 1{|Q|}\sum_{i,i'\in Q}v(i-i')$ and used $\widetilde Q=[-R+S,R-S]^d\subset Q$ as defined in Step~\ref{step-QtoZd} in Section~\ref{sec-lowbound} (where $S$ is chosen in such a way that there is no interaction between particles associated with points in $\widetilde Q$ and $Q^{\rm c}$). Recall that $\v=\sum_{i\in\Z^d} v(i)$, and $N_Q^{\ssup\ell}(\omega)=\sum_{k\in\N}\sum_{x\in Q}k\xi^{\ssup k}(x)$ is the number of particles in marks at points in $Q$ in the configuration $\omega$.

We have the following lower bound on the energy:
\begin{equation}\label{energylowbound}
\Phi_{\widetilde \Lambda_N,\widetilde \Lambda_N}\big(\omega+\omega_X\big)\geq \sum_{z\in Y_{N}}\widetilde \Phi_{Q}\big(\omega^{\ssup z}, a^{\ssup z}\big).
\end{equation}

Since we have estimated the influence of the macroscopic Poisson points $X_{k,j}$ with $k>K$  and $j\in[l_k^{\ssup{N}}]$ and their marks in terms of the $  a^{\ssup z}$, we can proceed by taking the maximum over all these $ a^{\ssup z}$:
\begin{equation}\label{UppBound3}
\begin{aligned}
&Z_{N,\L_N,{\rm Dir}}(l^{\ssup{N}})\leq\e^{-C|\Lambda_N|}+\\
&\max_{( a^{\ssup z})_z\in\N_0^{Y_{N}}\colon \sum_z  a^{\ssup z}\in |Q|^{-1}|\Lambda_N|J_{K,R,\delta}}\E_{\widetilde\L_N}^{\ssup{\leq K}}\Big[\e^{-\sum_{z\in Y_{N}}\widetilde \Phi_{Q}(\cdot^{\ssup z}, a^{\ssup z})}\prod_{k=1}^K\1\{N_{\widetilde\Lambda_N}^{\ssup{\delta_k}}\in m_k|\Lambda_N|[1-\delta,1+\delta]\}\Big],
\end{aligned}
\end{equation}
where 
\begin{equation}\label{eq_J}
J_{K,R,\delta}:=\Big[\Big( \rho_{\rm ma}+\sum_{k>K} k m_k-\delta\Big(\sum_{k\in [K]}km_k\Big)\Big)(1-\delta_{K,R}),\rho_{\rm ma}+\sum_{k>K} k m_k+\delta\Big(\sum_{k\in [K]}km_k\Big)\Big]\cap[0,\rho].
\end{equation}

\begin{step}
We rewrite the expectation in terms of an expectation with respect to  a crucial empirical measure.
\end{step}

Now we introduce the main tool in our large-deviation analysis, the empirical measure of the subconfigurations in the $Q_z$,
\begin{equation}\label{etadef}
\eta_{N,R}(\omega,\bar a)=\frac 1{|Y_N|}\sum_{z\in Y_{N}}\delta_{(\omega^{\ssup z}, a^{\ssup z})}\in\Mcal_1(\Omega_Q^{\ssup {\leq K}}\times\N_0),\qquad \omega \in\Omega_{\widetilde\L_N}, \bar a\in\N_0^{Y_{N}}.
\end{equation}
In terms of $\eta_{N,R}$, we may write
$$
\sum_{z\in Y_{N}}\widetilde \Phi_{Q}(\omega^{\ssup z}, a^{\ssup z})=|Y_N|\langle \eta_{N,R}(\omega,\bar a),\widetilde\Phi_Q\rangle
$$
and
$$
\{N_{\Lambda_N}^{\ssup{\delta_k}}(\omega)=l_k\}=\Big\{\langle \eta_{N,R}(\omega,\bar a),N_Q^{\ssup{\delta_k}}\rangle=\frac {l_k}{|Y_N|}\Big\},\qquad  k\in[K],l_k\in\N,\bar a\in\N_0^{Y_{N}},
$$
where we conceive $N_Q^{\ssup{\delta_k}}$ as the map $(\omega, a)\mapsto N_Q^{\ssup{\delta_k}}(\omega)=\sum_{x\in Q}\xi^{\ssup k}(x)$ if $\omega=\sum_{x\in Q}\sum_{k\in[K]}\xi^{\ssup k}(x)\delta_{(x,G_k)}$. (Again, we make no notational difference between application to $\omega$ and its restriction to $Q$.)

Furthermore, observe that the condition on the $ a^{\ssup z}$ in \eqref{countingmacro} under the maximum can be written in terms of the projection of $\eta_{N,R}$ on the second component. Indeed, denoting this projection by $\pi_{\N_0}(\omega, a)= a$ and recalling that $|Y_N|\sim |\L_N|/|Q|$ as $N\to\infty$ and the asymptotics in \eqref{NMaasy}, we see that \eqref{countingmacro} implies
\begin{equation}\label{JKRdef}
\begin{aligned}
\langle \eta_{N,R}(\omega,\bar a),\pi_{\N_0}\rangle
& \in\frac{N^{\ssup{\rm Ma}}}{|Y_N|}|Q|^{-1}\big[1-\delta_{K,R},1\big]\subset J_{K,R,\delta},
\end{aligned}
\end{equation}
for any sufficiently large $N$. 
Collecting everything, this means that
\begin{equation}\label{upperbound2}
\begin{aligned}
Z_{N,\L_N,{\rm Dir}}(l^{\ssup{N}})\leq \e^{-C|\Lambda_N|}+ 
\max_{\bar a\in\N_0^{Y_{N}}}&\E_{\widetilde\L_N}^{\ssup{\leq  K}}\Big[\e^{-|Y_N|\langle \eta_{N,R}(\cdot, \bar a),\widetilde\Phi_Q\rangle}\1\big\{\langle \eta_{N,R}(\cdot, \bar a),\pi_{\N_0}\rangle\in J_{K,R,\delta}\big\}\\
&\prod_{k=1}^K \1\Big\{\Big\langle \eta_{N,R}(\cdot, \bar a),\frac 1{|Q|}N_Q^{\ssup{\delta_k}}\Big\rangle\in m_k[1-\delta,1+\delta]\Big\}\,\Big].
\end{aligned}
\end{equation}

\begin{step}
We carry out the large-$N$ asymptotics with the help of large-deviation arguments. 
\end{step}

We introduce the rate function
\begin{equation}\label{IQleqKdef}
\Mcal_1(\Omega_Q^{\ssup {\leq K}}\times\N_0)\ni\eta\mapsto I_Q^{\ssup{\leq K}}(\eta)= \frac 1{|Q|}\sum_{ a\in\N_0}\eta\big(\Omega^{\ssup{\leq K}}_Q\times\{ a\}\big)\,H_Q\big(\eta(\cdot| a)| \P_Q^{\ssup{\leq K}}\big),
\end{equation}
where we wrote $\eta(\cdot| a)=\eta(\cdot\times\{ a\})/\eta(\Omega_Q^{\ssup{\leq K}}\times\{ a\})$ for the conditional distribution given the type $ a$.

\begin{lemma}[$N\to\infty$]\label{lem-Ntoinfty}
If $C$ is large enough, then, for any $K,R\in\N$ and $\delta\in(0,1)$,
$$
\begin{aligned}
\limsup_{N\to\infty}&\frac1{|\Lambda_N|}\log Z_{N,\L_N,{\rm Dir}}(l^{\ssup{N}})
&\leq -\varphi_{Q}^{\ssup{K,\delta}}(m^{\ssup {\leq K}},J_{K,R,\delta}),
\end{aligned}
$$
where the interval $J_{K,R,\delta}$ is defined in \eqref{eq_J}, and for any $J\subset[0,\infty)$ and for any $m=(m_1,\dots,m_K)\in[0,\infty)^K$, we define
\begin{equation}\label{chiQKdef}
\begin{aligned}
\varphi_{Q}^{\ssup{K,\delta}}(m,J)&=\inf\Big\{I^{\ssup{\leq K}}_Q(\eta)+\Big\langle \eta,\frac1{|Q|}\widetilde \Phi_{Q}\Big\rangle\colon \eta\in \Mcal_1(\Omega_Q^{\ssup{\leq K}}\times\N_0), \\
&\qquad\forall k\in[K]\colon \Big\langle \eta,\frac1{|Q|}N_Q^{\ssup{\delta_k}}\Big\rangle\in m_k [1-\delta,1+\delta],\langle \eta,\pi_{\N_0}\rangle\in J\Big\}.
\end{aligned}
\end{equation}
\end{lemma}

\begin{proof} 
We need a large-deviation principle (LDP) for $\eta_{N,R}$ as $N\to\infty$. By the properties of a PPP, if $\omega$ is an $\Omega_{\widetilde\L_N}^{\ssup{\leq K}}$-valued random variable under $\P_{\widetilde \L_N}^{\ssup{\leq K}}$, the family $(\omega^{\ssup z})_{z\in Y_{N}}$ is an i.i.d.~sequence with values in $\Omega_Q^{\ssup{\leq K}}$ with distribution $\P_Q^{\ssup{\leq K}}$, the distribution of the projection of the marked reference PPP $\P$ to $\Omega_Q^{\ssup{\leq K}}$. However, these $\omega^{\ssup z}$'s come with disordered values $a^{\ssup z}$, such that an application of the usual Sanov theorem is not possible. We use a variant of Sanov's theorem with a countable number of types of distributions labeled by $a\in\N_0$. This says that $\eta_{N,R}(\omega,\bar a)$ satisfies on $\Mcal_1(\Omega_Q^{\ssup{\leq K}}\times\N_0)$ an LDP with speed $|\Lambda_N|$ and rate function $I_Q^{\ssup{\leq K}}$ defined in \eqref{IQleqKdef}. 

Let us give some explanations, as we actually do not know an explicit reference for this LDP. Consider the $\N_0$-projection $\pi_{\N_0}\eta_{N,R}$ of $\eta_{N,R}$, then \eqref{upperbound2} may be written as
$$
\begin{aligned}
&Z_{N,\L_N,{\rm Dir}}(l^{\ssup{N}})\leq \e^{-C|\Lambda_N|}+ \sup_{\psi\in \Mcal_1(\N_0)\colon \sum_{a\in\N_0}a\psi(a)\in J_{K,R,\delta}}\max_{\bar a\in\N_0^{Y_{N}}}
\E_{\widetilde\L_N}^{\ssup{\leq  K}}\Big[\e^{-|Y_N|\langle \eta_{N,R}(\cdot, \bar a),\widetilde\Phi_Q\rangle}\\
&\qquad\Big[\prod_{k=1}^K \1\Big\{\Big\langle \eta_{N,R}(\cdot, \bar a),\frac 1{|Q|}N_Q^{\ssup{\delta_k}}\Big\rangle\in m_k[1-\delta,1+\delta]\Big\}\,\1\{\pi_{\N_0}\eta_{N,R}(\cdot, \bar a)=\psi\}\Big].
\end{aligned}
$$
For any fixed $\psi$, on the event $\{\pi_{\N_0}\eta_{N,R}(\cdot, \bar a)=\psi\}$, the measure $\eta_{N,R}$ is the empirical measure of $|Y_N|$ independent random variables $(\omega^{\ssup z},a^{\ssup z})$. Furthermore, for any $a\in\N_0$, the second argument $a^{\ssup z}$ is equal to $a$ for $\sim|Y_N| \psi(a)$ of them. Hence, $\eta_{N,R}$ is a convex combination of empirical measures indexed by $a\in\N_0$ mixed according to $\psi$, each of which satisfies Sanov's theorem, i.e., an LDP on the scale $|Y_N|$ with rate function equal to the entropy with respect to the reference distribution, $\P_Q^{\ssup{\leq K}}$. This implies the LDP on the scale $|\L_N|\sim |Y_N||Q|$ for $\eta_{N,R}$ with rate function equal to the convex combination in \eqref{IQleqKdef}.

Let us discuss the used topology. We are using the smallest topology on the countable set $\Omega_Q^{\ssup {\leq K}}$ that makes continuous any functional of the form $\omega=\sum_{x\in Q}\sum_{k\in[K]}\xi^{\ssup k}(x)\delta_{(x,G_{k})}\mapsto N_{Q'}^{\ssup f}(\omega)=\sum_{x\in Q'}\sum_{k\in[K]} f(k)\xi^{\ssup k}(x)$ with measurable $Q'\subset Q$ and  bounded $f\colon \N_0\to\R$. In particular, $N_Q^{\ssup{\delta_k}}$ belongs to this class. On $\Omega_Q^{\ssup {\leq K}}\times\N_0$ we use the product topology.

On the set $\Mcal_1(\Omega_Q^{\ssup {\leq K}}\times\N_0)$ we use the weak topology induced by test integrals against all continuous bounded functionals $\Omega_Q^{\ssup {\leq K}}\times\N_0\to\R$. On the scale $|\L_N|$, the functional  $\eta\mapsto -\langle \eta,\frac1{|Q|}\widetilde\Phi_Q\rangle$ in the exponent is clearly  bounded from above and upper semi-continuous ($\widetilde\Phi_Q$ can be monotonously increasing approximated by linear combinations of nonnegative functions as above), and the functional $\eta\mapsto \langle \eta,\frac1{|Q|}N_Q^{\ssup{\delta_k}}\rangle$ is continuous for any $k\in[K]$, hence the $k$-dependent indicators in the expectation on the right-hand side of \eqref{upperbound2} are on closed subsets. 

The state space of the LDP that we use is countable, and on the indicators on the right-hand side of \eqref{upperbound2},  the only unboundedness comes from the second factor of the state space, $\N_0$. However, there is a simple compactification argument with the help of the energy functional $\widetilde \Phi_Q$ defined in \eqref{energydecompose} as follows. Consider in the expectation on the right-hand side of \eqref{upperbound2}, for a large auxiliary parameter $L$, the event $\{\langle \eta_{N,R},\pi_{\N_0}^2\rangle>L\}=\{\sum_{a\in\N_0}\eta_{N,R}(\Omega_Q^{\ssup{\leq K}}\times\{a\})a^2>L\}$ and its complement; on the first event, the interaction  can be estimated by
$$
\e^{-|Y_N|\langle\eta_{N,R},\widetilde \Phi_Q\rangle}
\leq  \e^{-|Y_N| \v_{Q} \langle \eta_{N,R},\pi_{\N_0}^2\rangle}
\leq \e^{-\v_Q L|\Lambda_N|/|Q|}
$$
and the entire term gives a negligible exponential contribution on the scale $|\L_N|$ in the limit $N\to\infty$ if $L$ is sufficiently large. The second event makes the considered set of $\eta$'s compact, i.e., the set of all $\eta$ such that 
$$
\Big\langle \eta,\frac1{|Q|}N_Q^{\ssup{\delta_k}}\Big\rangle\in m_k [1-\delta,1+\delta]\, \forall k\in[K]\qquad\mbox{and}\qquad \langle\eta,\pi_{\N_0}\rangle\in J_{K,R,\delta}\qquad\mbox{and}\qquad \langle \eta,\pi_{\N_0}^2\rangle\leq L.
$$ 
This is seen from the estimate $\sum_{a>A}a \eta(\Omega_Q^{\ssup {\leq K}}\times\{a\})\leq \frac 1A \langle \eta,\pi_{\N_0}^2\rangle\leq L/A\to 0$ as $A\to\infty$. (Alternatively, one can also show that the used LDP for $\eta_{N,R}$ with fixed $\psi=\pi_{\N_0}\eta_{N,R}$ is locally uniform in $\psi$ and argue that the corresponding set of $\psi$'s is compact.)

Hence, we may apply the upper-bound part of Varadhan's lemma and obtain  the assertion of the lemma, also noting that the infimum over $\psi$ such that $\sum_a \psi(a)a\in J_{K,R,\delta}$ of the infimum over $\eta$ satisfying $\pi_{\N_0}\eta=\psi$ can be summarised as the infimum over $\eta$ such that $\langle \eta,\pi_{\N_0}\rangle \in J_{K,R,\delta}$.
\end{proof}

\begin{step}
We relax the restrictions of the constraints and on the boundedness of the mark size. 
\end{step}

Recall that $ \rho_{\rm ma}$ was defined in \eqref{rho1rho2} and the interval $J_{K,R,\delta}$ in \eqref{eq_J}.

\begin{lemma}[$K\to \infty$ and $\delta\to 0$]\label{lem-Ktoinfty}
For any $R\in(0,\infty)$ and  $m=(m_k)_{k\in\N}\in[0,\infty)^{\N}$,
\begin{equation}
\begin{aligned}
\liminf_{K\to\infty,\delta\to 0}&\varphi_Q^{\ssup{K,\delta}}(m^{\ssup{\leq K}},J_{K,R,\delta})\geq \inf_{\psi\in\Mcal_1(\N_0)\colon\sum_{a\in\N_0}a\psi(a)=\rho_{\rm ma}}\varphi_Q(m,\psi),
\end{aligned}
\end{equation}
where 
\begin{equation}\label{varphiKdeltaQpsidef}
\begin{aligned}
\varphi_Q(m,\psi)=&\inf\Big\{\sum_{a\in\N_0}\psi(a)\Big[\frac1{|Q|}H_Q(\eta_a|\P_Q)+ \frac1{|Q|}\langle\eta_a,\Phi_{Q,Q}\rangle +2{\v}a\frac1{|Q|}\langle\eta_a,N_{\widetilde Q}^{\ssup\ell}\rangle +\v_{Q} a^2 \Big]\colon\\
&\qquad \eta_0,\eta_1,\dots \in\Mcal_1(\Omega_Q), \sum_{a\in\N_0}\psi(a)\frac1{|Q|}\langle \eta_a ,N_Q^{\ssup{\delta_k}}\rangle = m_k\ \forall k\in \N\Big\}.
\end{aligned}
\end{equation}
\end{lemma}

\begin{proof}
First we isolate the $a$-dependence by substituting $\psi(a)=\eta(\Omega^{\ssup {\leq K}}_Q\times\{a\})$ and $\eta_a(\cdot)=\eta(\cdot\times\{a\})/\eta(\Omega^{\ssup {\leq K}}_Q\times\{a\})$, then we see that $\varphi_Q^{\ssup{K,\delta}}(m^{\ssup{\leq K}}, J)$ can be reformulated as 
$$
\begin{aligned}
\varphi^{\ssup{K,\delta}}_Q(m^{\ssup{\leq K}},J)=&\inf\Big\{\sum_{a\in\N_0}\psi(a)\Big[\frac1{|Q|}H_Q(\eta_a|\P^{\ssup {\leq K}}_Q)+ \frac1{|Q|}\langle\eta_a,\Phi_{Q,Q}\rangle +2{\v}a\frac1{|Q|}\langle\eta_a,N_{\widetilde Q}^{\ssup\ell}\rangle +\v_{Q}a^2 \Big]\colon\\
&\qquad \eta_0,\eta_1,\dots \in\Mcal_1(\Omega^{\ssup {\leq K}}_Q), \psi\in\Mcal_1(\N_0), \\
&\qquad\sum_{a\in\N_0}a\psi(a)\in J,\sum_{a\in\N_0}\psi(a) \frac1{|Q|}\langle \eta_a ,N_Q^{\ssup{\delta_k}}\rangle\in m_k[1-\delta,1+\delta]\ \forall k\in [K]\Big\}.
\end{aligned}
$$
Now use this for $J=J_{K,R,\delta}$ and pick $L$ large enough such that an insertion of the constraint $\sum_{a\in\N_0}a^2\psi(a)\leq L$  in the right-hand side does not change its value for any $K\in\N$ and $\delta\in(0,1)$.
Pick, for any $K\in\N$ and $\delta\in(0,1)$, an (approximative) minimiser $(\psi^{\ssup{K,\delta}}, (\eta^{\ssup{K,\delta}}_a)_{a\in\N_0})$ for this variational problem. Now pick some sequence $(K_n)_{n\in\N}$ and $(\delta_n)_{n\in\N}$ such that $K_n\to\infty$ and $\delta_n\downarrow 0$ as $n\to\infty$. Pick some accumulation point $\psi$ for the family $(\psi^{\ssup{K_n,\delta_n}})_{n\in\N}$, then $\sum_{a\in\N_0}a\psi(a)=\rho_{\rm ma}$, since $\bigcap_{\delta\in(0,1)}J_{K,R,\delta}=\{ \rho_{\rm ma}+\sum_{k>K}km_k\}$.

View $\eta_a^{\ssup{K_n,\delta_n}}$ as an element of $\Mcal_1(\Omega_Q)$ by adding an empty configuration of particles with marks of size $>K$. Then, the reference measure $\P_Q$ weights also these empty configurations with density $q_k$ for configurations with mark $G_k$. Since $\sum_k q_k$ is finite, it is easy to calculate that
$$
H_Q(\eta_a^{\ssup{K_n,\delta_n}}|\P^{\ssup{\leq K_n}}_Q)=H_Q(\eta_a^{\ssup{K_n,\delta_n}}|\P_Q)+\eps_n
$$
with some $\eps_n\to0$ as $n\to\infty$.  Now use that $H_Q(\cdot|\P_Q)$ has compact level sets, a fact that is proved in the proof of Sanov's theorem, since this is a good rate function. Hence, we find, jointly for all $a\in\N_0$ satisfying $\psi(a)>0$, a subsequence of $\eta^{\ssup{K_n,\delta_n}}_a$ as $n\to\infty$ with limit $\eta_a\in \Mcal_1(\Omega_Q)$. For notational convenience, we assume that $(\psi^{\ssup{K_n,\delta_n}},(\eta_a^{\ssup{K_n,\delta_n}})_{a\in\N})$ converges towards $(\psi,(\eta_a)_{a\in\N})$. (We tacitly dropped all $a$ with $\psi(a)=0$ from this sequence, which gives a lower bound for the functional.)

Then we see that $(\psi,(\eta_a)_{a\in\N_0})$ is admissible in the set on the right-hand side of \eqref{varphidef}, since the map $\Mcal_1(\Omega_Q)\ni \eta\mapsto \langle \eta,N_Q^{\ssup{\delta_k}}\rangle$ is continuous for any $k\in\N$ and the map $(\psi,(\eta_a)_{a\in\N_0})\mapsto \sum_{a\in\N_0}\psi(a) \frac1{|Q|}\langle \eta_a ,N_Q^{\ssup{\delta_k}}\rangle$ is continuous on the set where $\sum_{a\in\N_0}a^2\psi(a)\leq L$.

Furthermore, use that the functional that is minimised on the right-hand side of \eqref{varphiKdeltaQpsidef} is lower semi-continuous, to end the proof.
\end{proof}

\begin{step}
We let $Q\uparrow \Z^d$.
\end{step}

Here we will be using the spatial ergodic theorem and the definition of the limiting entropy density $I$ defined in \eqref{Igriddef}. Recall the notation of Theorem \ref{thm-freeenergygrid}.

\begin{lemma}[$Q\uparrow \Z^d$]\label{lem-Qtoinfty} For any $m=(m_k)_{k\in\N}\in[0,\infty)^{\N}$ and any $\rho_{\rm ma}\in[0,\infty)$, 
\begin{equation}\label{varphilowboundQ}
\liminf_{Q\uparrow \Z^d}\inf_{\psi\in\Mcal_1(\N_0)\colon \sum_a a\psi(a)=\rho_{\rm ma}}\varphi_{Q}(m,\psi)\geq\inf_{\psi\in\Mcal_1(\N_0)\colon \sum_a a\psi(a)=\rho_{\rm ma}}\varphi(m,\psi).
\end{equation}
\end{lemma}

\begin{proof}

Fix a large number $L>0$ such that the insertion of the condition $\sum_a \psi(a)a^2\leq L$ in the two infima in \eqref{varphilowboundQ} does not change anything in the values of the infima, and note that these two infima range over a compact set. Hence, it will be sufficient to prove \eqref{varphilowboundQ} for a fixed $\psi$, i.e., we will show only that $\liminf_{Q\uparrow \Z^d}\varphi_{Q}(m,\psi)\geq\varphi(m,\psi)$. We may assume that $\varphi_Q(m,\rho_{\rm ma})$ is bounded as $Q\uparrow \Z^d$ along the sequence of $Q_n$'s that we consider, otherwise there is nothing to be shown.

Let a small $\eps>0$ be given, and assume that $(\eta_{Q,a})_{a\in\N_0}$ is, for any $Q$, an $\eps$-approximate minimiser in the formula in the right-hand side of \eqref{varphiKdeltaQpsidef}. We construct now a measure $P_a^{\ssup Q}\in\Mcal_1^{\ssup{\rm s}}(\Omega)$  as follows.  Recall that $\Z^d$ is decomposed into the sets $Q_z=z+Q$ with $z\in 2R\Z^d$, and put in each of the $Q_z$'s an independent copy of a configuration with distribution equal to $\eta_{Q,a}$. We write the arising distribution as $\eta_{Q,a}^{\otimes 2R\Z^d}$. Now put
$$
P_a^{\ssup Q}=\frac 1{|Q|}\sum_{z\in Q}\eta_{Q,a}^{\otimes 2R\Z^d}
\circ \theta_z^{-1}\in \Mcal_1^{\ssup{\rm s}}(\Omega).
$$
We want to show that $P_a^{\ssup Q}$ has a converging subsequence as $Q\uparrow \Z^d$. For this we will be using that the level sets of $I$ are compact. For using this, we need to show that $I(P_a^{\ssup Q})$ is bounded in $Q$. This goes as follows. We have that
$$
\begin{aligned}
I(P_a^{\ssup Q})&=\inf_{\widetilde Q\subset\Z^d}\frac1{|\widetilde Q|}H_{\widetilde Q}\big((P_a^{\ssup Q})_{\widetilde Q}\big|\P_{\widetilde Q}\big)\\
&\leq \frac1{| Q|}H_{Q}\big((P_a^{\ssup Q})_Q\big|\P_{Q}\big)
\leq \frac 1{|Q|}\sum_{z\in Q}\frac1{| Q|}H_{Q}\big((\eta_{Q,a}^{\otimes 2R\Z^d}
\circ \theta_z^{-1})_Q\big|\P_{Q}\big)\\
&=\frac1{| Q|}H_{Q}(\eta_{Q,a}|\P_Q),
\end{aligned}
$$
using the definition of $I$ (the infimum ranges over all centred boxes $\widetilde Q$), the convexity of $H_Q$ and the shift-invariance of the reference measure, $\P$. Put $P^{\ssup Q}=\sum_{a\in\N_0}\psi(a) P_a^{\ssup Q}$. For any $a\in\N$ such that $\psi(a)>0$, we therefore have that 
$$
I(P_a^{\ssup Q})\leq \frac1{\psi(a)}\sum_{\widetilde a\in\N_0}\psi(\widetilde a)I(P_{\widetilde a}^{\ssup Q})\leq  \frac1{\psi(a)}\sum_{\widetilde a\in\N_0}\psi(\widetilde a)\frac1{| Q|}H_{Q}(\eta_{Q,\widetilde a}|\P_Q)\leq \frac1{\psi(a)}\big(\varphi^{\ssup{\leq K,\delta}}_Q(m,\psi)+\eps\big),
$$
since the energy terms are nonnegative and $(\eta_{Q,a})_{a\in\N_0}$ is an $\eps$-approximate minimiser, and $\varphi^{\ssup{\leq K,\delta}}_Q(m,\psi)$ is bounded in $Q$. Since this upper bound is bounded in $Q\subset\Z^d$, we now know that $P_a^{\ssup Q}$ has a convergent subsequence as $Q\uparrow \Z^d$. The topology used is the one that is induced by the test integrals against any local and tame function, i.e., against any function $\Omega\to\R$ that depends only on some bounded box $\Lambda\subset\Z^d$ and can in absolute value be upper estimated against a constant plus a constant times $N_\L^{\ssup\one}$, the number of points in $\L$.

Denote the limit by $P_a$ and put $P=\sum_{a\in\N_0}\psi(a)P_a$. By lower semicontinuity and affinity  of $I$, we see that 
\begin{equation}\label{Isemicont}
I(P)\leq \liminf_{Q\uparrow \Z^d}\sum_{a\in\N_0}\psi(a)\frac1{|Q|}H_Q(\eta_a|\P_Q^{\ssup{\leq K}}).
\end{equation}

Note that, by shift-invariance of $P^{\ssup Q}$, and since the maps $P\mapsto \langle P,N_0^{\ssup{\delta_k}}\rangle$ are continuous for any $k\in\N$, we have  
\begin{equation*}\label{countingfunctionals}
\langle P,N_0^{\ssup{\delta_k}}\rangle = m_k\ \forall k\in \N,
\end{equation*}
and in particular, $P$ is admissible in  the formula in \eqref{varphidef}. Now we turn to analogous estimates for the remaining terms in the first line of~\eqref{varphiKdeltaQpsidef}, where the last term is harmless. For the one but last term, note that, for all $a\in \N_0$, we have that
\begin{equation}\label{2ndLastTerm}
\langle P_a,N_0^{\ssup\ell}\rangle\le \liminf_{Q\uparrow\Z^d}\frac1{|Q|}\langle \eta_{Q,a},N_{\widetilde Q}^{\ssup\ell}\rangle,
\end{equation}
where we note that $\liminf_{Q\uparrow\Z^d}|\widetilde Q|/|Q|=1$. 
In order to see~\eqref{2ndLastTerm}, note that, for any $K\in\N$, since $P_a^{\ssup Q}\to P_a$ and $N_0^{\ssup{\delta_k}}$ is local and tame,
$$
\Big\langle P_a,\sum_{k\in[K]}kN_0^{\ssup{\delta_k}}\Big\rangle 
=\lim_{Q\uparrow\Z^d}
\Big\langle P_a^{\ssup Q},\sum_{k\in[K]}kN_0^{\ssup{\delta_k}}\Big\rangle;
$$
furthermore, for any box $Q$,
$$
\Big\langle P_a^{\ssup Q},\sum_{k\in[K]}kN_0^{\ssup{\delta_k}}\Big\rangle
\leq \langle P_a^{\ssup Q},N_0^{\ssup\ell}\rangle
=\frac1{|Q|}\Big\langle P_a^{\ssup Q},N_Q^{\ssup{\ell}}\rangle 
=\frac1{|Q|}\Big\langle \eta_{Q,a},N_Q^{\ssup{\ell}}\Big\rangle.
$$
where we used the shift invariance of $P_a^{\ssup Q}$. Now make $K\to\infty$ on the left-hand side of the one-but-last display to get the assertion.

In a similar fashion we show for the second term in the first line of~\eqref{varphiKdeltaQpsidef} that
\begin{equation}\label{ergodictheorem}
\langle P_a,\Phi_{0,\Z^d}\rangle\leq \liminf_{Q\uparrow\Z^d}\frac1{|Q|}\langle \eta_{Q,a},\Phi_{Q,Q}\rangle,\qquad a\in\N_0,
\end{equation}
where $\frac1{|Q|}\langle \eta_{Q,a},\Phi_{Q,Q}\rangle=  \frac1{|Q|}\sum_{z\in Q}\langle \eta_{Q,a},\Phi_{z,Q}\rangle$.
For this, we approximate $\Phi_{z,Q}$ from below with local tame functions. First note that for all $a\in \N_0$, $z\in Q$, bounded $\L\subset \Z^d$ and $S,K\in\N$, we have that $\Phi_{z,Q}\ge\Phi^{\ssup{\L,S,K}}_{z,Q}\rangle$, where $\Phi^{\ssup{\L,S,K}}_{z,Q}=\Phi_{z,Q}\one\{N_{\L+z}^{\ssup \one}\le S\}\prod_{k>K}\one\{N_{\L+z}^{\ssup{\delta_k}}=0\}$ is a local and tame function since marks have a maximal cardinality $K$ and the number of points involved is bounded by $S$. Further, 
\begin{align*}
\frac1{|Q|}\sum_{z\in Q}\langle \eta_{Q,a},\Phi^{\ssup{\L, S,K}}_{z,Q}\rangle-\langle P_a^{\ssup Q},\Phi^{\ssup{\L, S,K}}_{0,\Z^d}\rangle&=\frac1{|Q|}\sum_{z\in Q}\Big[\langle \eta_{Q,a},\Phi^{\ssup{\L, S,K}}_{z,Q}\rangle-\langle \eta_{Q,a}^{\otimes 2R\Z^d},\Phi_{z,\Z^d}\rangle\Big]\\
&=-\frac1{|Q|}\sum_{z\in Q}\langle \eta_{Q,a}^{\otimes 2R\Z^d},\Phi^{\ssup{\L, S,K}}_{z,\Z^d\setminus Q}\rangle.
\end{align*}
Now, since the marks have a maximal cardinality $K$ and the support of $v$ is assumed to be finite, there exists a bounded set $\L\subset\Z^d$ such that 
\begin{equation*}
\frac1{|Q|}\sum_{z\in Q}\langle \eta_{Q,a}^{\otimes 2R\Z^d},\Phi^{\ssup{\L, S,K}}_{z,\Z^d\setminus Q}\rangle
=\frac1{|Q|}\sum_{z\in Q\colon \L+z\not\subset Q}\Big[\langle \eta_{Q,a}^{\otimes 2R\Z^d},\Phi^{\ssup{\L, S,K}}_{z,\Z^d\setminus Q}\rangle\Big]\le \frac{|z\in Q\colon \L+z\not\subset Q|}{|Q|}(KS)^2\v,
\end{equation*}
which tends to zero as $Q$ tends to $\Z^d$. 

Now, since $\Phi^{\ssup{\L, S,K}}_{0,\Z^d}$ is local and tame, we have
\begin{equation*}
\liminf_{Q\uparrow\Z^d}P_a^{\ssup Q}(\Phi^{\ssup{\L, S,K}}_{0,\Z^d})= P_a(\Phi^{\ssup{\L, S,K}}_{0,\Z^d}).
\end{equation*}
Letting $\L$ tend to $\Z^d$ und $S,K$ tend to infinity we arrive at the desired result. 

Collecting \eqref{Isemicont}, \eqref{2ndLastTerm} and \eqref{ergodictheorem} shows that
$$
\sum_{a\in\N_0}\psi(a)\Big[I(P_a)+P_a(\Phi_{0,\Z^d})+2{\v}aP_a(N_0^{\ssup \ell})+{\v}a^2 \Big]\leq \eps+\liminf_{Q\uparrow\Z^d}\varphi_Q(m,\psi).
$$
Since $\varphi(m,\psi)$ is not larger than the left-hand side (since $P$ is admissible in \eqref{varphidef}), we arrived at the claim.
\end{proof}

\begin{step}\label{step-uppbound9}
We finish the proof of the upper bound in Theorem \ref{thm-freeenergygrid}.
\end{step} 

We apply first \eqref{ZNgegenZN(l)} (recalling that the counting term is not larger than $\e^{o(N)}$) and then Lemma~\ref{lem-Ntoinfty}, to see that, for any box $Q$ and any $K\in\N$ and $\delta\in(0,1)$,
\begin{equation}\label{Proofupbound1}
\limsup_{N\to\infty}\frac1{|\Lambda_N|}\log Z_{N,\Lambda_N,{\rm Dir}}
\leq - \varphi_Q^{\ssup{K,\delta}}(m^{*,\ssup{\leq K}},J_{K,R,\delta}).
\end{equation}
Here we recall that $m^{*}$ was defined via a convergent subsequence, and  $\rho_{\rm ma}$ via Fatou's lemma, see \eqref{rho1rho2}.

Applying Lemma \ref{lem-Ktoinfty} for making $K\to\infty$ and $\delta\downarrow 0$  and then  Lemma \ref{lem-Qtoinfty} for making $Q\uparrow \Z^d$, the left-hand side  of \eqref{Proofupbound1} can also be estimated against $-\inf_{\psi\in\Mcal_1(\N_0)\colon \sum_{a\in\N_0}a\psi(a)=\rho_{\rm ma}}\varphi(m^*,\psi)$.

Since the latter is obviously not larger than the right-hand side of \eqref{freeenergyvarform} (we have $\sum_{k\in\N}k m^*_k+\rho_{\rm ma}=\rho$ from \eqref{rho1rho2}), this finishes the proof of the upper bound in Theorem \ref{thm-freeenergygrid}, \eqref{ProofThm1}.

As a by-product of the preceding proof, we now also have a proof of the first part of Lemma \ref{lem-ProofExistMin}.

\begin{proof}[Proof of the first part of Lemma \ref{lem-ProofExistMin}]

Using again the device that the infimum $\inf_{\psi\in\Mcal_1(\N_0)\colon \sum_{a\in\N_0}a\psi(a)=\rho_{\rm ma}}\varphi(m^*,\psi)$ can be restricted to some compact set by adding the constraint $\sum_{a\in\N_0}\psi(a)a^2\leq L$ for some large $L$, and using lower semicontinuity, we see that this infimum has a minimiser.

Then, from Step \ref{step-uppbound9} and the lower bound in Theorem \ref{thm-freeenergygrid}, we have, for $m^*$ and $ \rho_{\rm ma}$ constructed in \eqref{rho1rho2} and $\psi$ taken as a minimiser of $\varphi(m^*,\cdot)$ under $\sum_{a\in\N_0}\psi(a)a= \rho_{\rm ma}$,
$$
\begin{aligned}
\varphi(m^*,\psi)&\leq -\limsup_{N\to\infty}\frac1{|\Lambda_N|}\log Z_{N,\Lambda_N,{\rm Dir}}\leq -\liminf_{N\to\infty}\frac1{|\Lambda_N|}\log Z_{N,\Lambda_N,{\rm Dir}}\\
&\leq \inf_{ \rho_{\rm mi}, \rho_{\rm ma}\colon  \rho_{\rm mi}+ \rho_{\rm ma}=\rho}\chi( \rho_{\rm mi}, \rho_{\rm ma}),
\end{aligned}
$$
that is, $(m^*,\psi)$ is a minimiser. 
\end{proof}

Another by-product of the preceding proof of Theorem~\ref{thm-freeenergygrid} is the following characterisation of $\varphi$.
Fix $\psi\in\Mcal_1(\N_0)$ with $\sum_a a\psi(a)<\infty$. With a cutting parameter $A\in\N$, put $\psi^{\ssup A}=\psi+\sum_{a>A}\psi(a)(\delta_0-\delta_a)$, as in Step~\ref{step:compactsupp} of Section~\ref{sec-lowbound}. Pick $N_a=|\L_N|\psi^{\ssup A}([a,\infty))$ for $a\in[A]$. Then $\omega_{\psi^{\ssup A}}=\sum_{a=1}^A \delta_{(0,G_{N_a})}$ is a (non-random!) distribution of macroscopic boxes whose rescaled empirical measure  approaches $\psi$ in the limit $N\to\infty$, followed by $A\to\infty$. We recall that $\E^{\ssup{\leq K}}$ denotes expectation with respect to the restriction of the reference process to $\Omega^{\ssup{\leq K}}$, the set of point processes with marks $G_1,\dots,G_K$ only.

\begin{lemma}\label{lem-varphirepr} For any $m\in[0,\infty)^{\N}$ and $\psi\in\Mcal_1(\N_0)$,
\begin{equation}\label{widetildevarphiformula}
\varphi(m,\psi)=-\lim_{A,K\to\infty,\, \delta\downarrow 0}\lim_{N\to\infty}\frac 1{|\L_N|}\log Z_{N,\L_N}^{\ssup{A,K,\delta}}(m,\psi),
\end{equation}
where
\begin{equation}\label{ZNmpsidef}
Z_{N,\L_N}^{\ssup{A,K,\delta}}(m,\psi)=\E_{\L_N}^{\ssup{\leq K}}\Big[\e^{-\Phi_{\Lambda_N,\Lambda_N}(\cdot+\omega_{\psi^{\ssup A}})}\prod_{k\in[K]}\1\{N_{\L_N}^{\ssup{\delta_k}}\in m_k |\L_N|(1-\delta,1+\delta)\}\Big],
\end{equation}
and $\L_N$ is a centred box with  volume $N/\rho$, and $\rho=\sum_{k\in\N}k m_k+\sum_{a\in\N}\psi(a)a$. 
\end{lemma}

\begin{proof} The proofs of the upper and the lower bound of \eqref{widetildevarphiformula} are versions of the proofs of the upper and lower bounds in Theorem \ref{thm-freeenergygrid} that we detailed in  Sections \ref{sec-lowbound} and \ref{sec-uppbound}, respectively.  We decided not to give details, but let us hint for the lower bound at Steps~\ref{step:compactsupp} and \ref{step-lowbounddeltato0} of Section~\ref{sec-lowbound} and at \eqref{lowboundZdeltavarphidelta}, and we give now some few exemplary hints for the upper bound. 

In the proof of the upper bound, we now need to make the following two main differences: (1) instead of the (random) macroscopic potential $ \omega_X$ in \eqref{Upperboundcitelater}, here we have the (deterministic) $\omega_{\psi^{\ssup A}}$, and (2) the target upper bound is now in terms of one fixed $\psi$ instead of an infimum over many such functions.

Indeed, that proof shows that the right-hand side of \eqref{uppbounddeltaappears} is not larger than the right-hand side of \eqref{UppBound3} when the maximum over all the vectors  $(a^{\ssup z})_{z}$ is replaced by just the vector that is induced by $\omega_\psi$, after applying the cutting procedure described in Steps 2 and 3. In particular, for a given small $\eps>0$ and $R$ large enough, the empirical measure $\eta_{N,R}(\omega)$ defined in \eqref{etadef} lies in the set  $\Acal_\eps(\psi^{\ssup A})=\{\eta\colon |\eta^{\ssup {2}}(a)-\psi^{\ssup A}(a)|\leq\eps\,\forall a\in\{0,\dots,A\}\}$, where $\eta^{\ssup{2}}$ denotes the marginal measure of $\eta$ on $\N_0$. Hence, we arrive at  \eqref{upperbound2} for $Z_{N,\L_N}^{\ssup{A,K,\delta}}(m,\psi)$ instead of $Z_{N,\Lambda_N,{\rm Dir}}(l^{\ssup N})$ on the left-hand side and with the indicator on $\{\langle\eta_{N,R}(\omega),\pi_{\N_0}\rangle \in J_{K,R,\delta}\}$ on the right-hand side replaced by the indicator on $\{\eta_{N,R}(\omega)\in \Acal_\eps(\psi^{\ssup A})\}$.

Now we apply large-deviations arguments as in the proof of Lemma \ref{lem-Ntoinfty} and obtain
\begin{equation}\label{Zcutuppbound}
\limsup_{N\to\infty}\frac 1{|\L_N|}\log  Z_{N,\L_N}^{\ssup{A,K,\delta}}(m,\psi)\leq -\varphi_Q^{\ssup{K,\delta}}(m^{\ssup{\leq K}},\Acal_\eps(\psi^{\ssup A})),
\end{equation}
where $\varphi_Q^{\ssup{K,\delta}}(m,\Acal)$ is defined analogously to \eqref{chiQKdef} with $\langle \eta,\pi_{\N_0}\rangle\in J$ replaced by $\eta^{\ssup 2}\in\Acal$ for $\Acal\subset \Mcal_1(\N_0)$. 

The remainder of the proof is, as in Section \ref{sec-uppbound}, to let $Q\uparrow \Z^d$, $K\to \infty$ and $\delta,\eps\downarrow 0$ and $A\to\infty$ to see that in these limits, the right-hand side of \eqref{Zcutuppbound} is not larger than $-\varphi(m,\psi)$. The details are left to the reader.
\end{proof}

\section{Analysis of the variational formulas}\label{sec-AnaForm}

We prove  Lemma \ref{lem-properties} in Section \ref{sec-simplepropvarphi}, Lemma \ref{lem-Varformident} in Section \ref{sec-VarformProof} and Lemma \ref{lem-ProofExistMin} in Section \ref{sec-minimivarphi}.

\subsection{Properties of $I$, $\varphi$ and $\chi$: proof of Lemma \ref{lem-properties}}\label{sec-simplepropvarphi}

\setcounter{step}{0}

\begin{step}\label{step-convex}
$\varphi$ is convex.
\end{step}
\begin{proof}
The claim is that, for any $m, m' \in [0,\infty)^\N,\psi, \psi'\in\Mcal_1(\N_0)$ and any $\lambda\in [0,1]$,
$$
\varphi(\lambda m + (1-\lambda) m',\lambda \psi+(1-\lambda)\psi')\leq \lambda \varphi(m,\psi)+(1-\lambda)\varphi(m',\psi').
$$
Indeed, pick families $(P_a)_{a\in\N_0}$ and $(P'_a)_{a\in\N_0}$ admissible respectively in the formulas \eqref{varphidef} of $\varphi(m, \psi)$ and of $\varphi(m',\psi')$. Take $\widetilde P_a=\frac{\lambda \psi(a)}{\lambda \psi(a)+(1-\lambda)\psi'(a)}P_a+\frac{(1-\lambda) \psi'(a)}{\lambda \psi(a)+(1-\lambda)\psi'(a)}P'_a$, then for any $k\in \N$:
$$
\sum_{a\in\N_0}(\lambda \psi(a)+(1-\lambda)\psi'(a))\widetilde P_a(N_0^{\ssup {\delta_k}})=\lambda m_k + (1-\lambda) m'_k,
$$
hence $(\widetilde P_a)_{a\in\N_0}$ is admissible in the formula \eqref{varphidef} of $\varphi(\lambda m + (1-\lambda) m',\lambda \psi+(1-\lambda)\psi')$. Using affinity of $I$, we get:
$$
\begin{aligned}
\varphi(\lambda m + &(1-\lambda) m',\lambda \psi+(1-\lambda)\psi')
\\
&\leq \sum_{a\in\N_0}(\lambda \psi(a)+(1-\lambda)\psi'(a))\Big[I(\widetilde P_a)+\widetilde P_a(\Phi_{0,\Z^d})+2{\v}a \widetilde P_a(N_0^{\ssup \ell})+{\v}a^2 \Big]\\
&= \lambda \sum_{a\in\N_0}\psi(a)\Big[I(P_a)+P_a(\Phi_{0,\Z^d})+2{\v}aP_a(N_0^{\ssup \ell})+{\v}a^2\Big]\\
&\qquad +(1-\lambda)\sum_{a\in \N_0}\psi'(a)\Big[I(P'_a)+P'_a(\Phi_{0,\Z^d})+2{\v}aP'_a(N_0^{\ssup \ell})+{\v}a^2 \Big].
\end{aligned}
$$
We conclude by taking the infimum over the families $(P_a)_{a\in\N_0}$ and $(P'_a)_{a\in\N_0}$.
\end{proof}

We now show the monotonicity of $\varphi$ in $\psi$ with respect to stochastic ordering.

\begin{step}\label{step-stochasnondec}
For any $m\in [0,\infty)^\N$, we have $\varphi(m,\psi)\leq \varphi(m,\psi')$  for any $\psi, \psi'\in\Mcal_1(\N_0)$ such that $\psi\leq_{\rm st} \psi'$.
\end{step}

\begin{proof}
Pick a family $(P_a)_{a\in \N_0}$ admissible in the definition of $\varphi(m,\psi')$. 
Since $\psi\leq_{\rm st} \psi'$, there exists $\pi\in \Mcal_1(\N_0\times\N_0)$ with marginals $\psi'$ and $\psi$ and such that $\pi(a,b)>0$ implies $a\geq b$.
Observe that 
\begin{align*}
     \sum_{a\in\N_0}&\psi'(a)\Big[I(P_a)+P_a(\Phi_{0,\Z^d})+2{\v}aP_a(N_0^{\ssup \ell})+{\v}a^2 \Big]\\
    & =\sum_{(a,b)\in\N_0^2}\pi(a,b)\Big[I(P_a)+P_a(\Phi_{0,\Z^d})+2{\v}aP_a(N_0^{\ssup \ell})+{\v}a^2 \Big]\\
    & \geq \sum_{(a,b)\in\N_0^2}\pi(a,b)\Big[I(P_a)+P_a(\Phi_{0,\Z^d})+2{\v}b P_a(N_0^{\ssup \ell})+{\v}b^2 \Big]\\
    & = \sum_{b\in\N_0}\psi(b)\Big[I(\widetilde P_b)+\widetilde P_b(\Phi_{0,\Z^d})+2{\v}b \widetilde P_b(N_0^{\ssup \ell})+{\v} b^2 \Big],
    \end{align*}
    where $\widetilde P_b=\sum_{a\in \N_0} \frac{\pi(a,b)}{\psi(b)}P_a$  and we used that $I$ is an affine function.
    Also note that for any $k\in \N$,
    $$\sum_{b\in\N_0}\psi(b)\widetilde P_b(N_0^{\ssup {\delta_k}})=\sum_{(a,b)\in \N_0^2} \pi(a,b)P_a(N_0^{\ssup {\delta_k}})=\sum_{a\in\N_0} \psi'(a)P_a(N_0^{\ssup {\delta_k}})=m_k,$$
    so $(\widetilde P_b)_{b\in \N_0}$ is admissible in the definition of $\varphi(m,\psi)$. Hence the inequality $\varphi(m,\psi')\geq \varphi(m,\psi)$ follows via a minimisation over $(P_a)_{a\in \N_0}$.
\end{proof}

Now we turn to the decomposition of the relative entropy. We recall that
\begin{equation}
H(m|q)=\sum_{k\in\N}\Big(q_k-m_k+m_k\log\frac{m_k}{q_k}\Big),\qquad m\in [0,\infty)^{\N_0}.
\end{equation}

\begin{step}\label{step_3a}
We have that $I(P)=H(m|q)+I_m(P)$, where $m_k=P(N_0^\ssup{\delta_k})$ and $I_m(P)$ is defined as in \eqref{Igriddef}, with $\P$ replaced by $\P^{m}$, the marked Poisson point process with $q$ replaced by $m$.
\end{step}

\begin{proof}
This formula follows from the fact that, for all finite $\L\Subset\Z^d$,
\begin{equation*}
H_\L(P_{\L}|\P_\L)=|\L|\sum_{k\in\N}\Big(q_k-m_k+m_k\log\frac{m_k}{q_k}\Big)+H_\L(P_{\L}|\P^{m}_\L),
\end{equation*}
where we recall that $P_\L$ is the projection of $P$ on the set of configurations in the set $\L$, and we write now $H_\L$ for the entropy on the set $\Mcal_1(\Omega_\L)$. Now use \eqref{Igriddef}. 
\end{proof}

Now we turn to the upper bound for $\chi(\rho)=\chi(\rho,0)$. Recall that $\chi^{\ssup{v=0}}(\rho,0)=\inf\{H(m|q)\colon m\in[0,\infty)^\N,\sum_k k m_k=\rho\}$ is the free energy of the non-interacting model.

\begin{step}\label{step-uppboundchi}We can upper bound $\chi$ as $\chi(\rho)\leq \chi^{\ssup{v=0}}(\rho,0)+\v \rho^2+\v \rho$ for any $\rho\in[0,\infty)$. 
\end{step}

\begin{proof}Fix $\rho\in[0,\infty)$ and take some $m\in [0,\infty)^\N$, such that $\sum_{k \in \N} k m_k=\rho$. Consider the reference process $\mathbb{P}^m\in \mathcal{M}_1^{\ssup{\rm s}}(\Omega)$ with $q$ replaced by $m$. It has density $\mathbb{P}^m(N_0^{\ssup\ell})=\rho$, relative entropy $I(\mathbb{P}^m)=H(m|q)$ and energy
\begin{align*}
    \mathbb{P}^m(\Phi_{0, \mathbb{Z}^d})
& = \sum_{y\in \mathbb{Z}^d}\sum_{k,l\in\mathbb{N}}\mathbb{E}^m\big[ \xi^{\ssup k}(0) \xi^{\ssup l}(y) \big]T_{0,y}(G_k,G_l) \\
& = \sum_{y\in \mathbb{Z}^d,\, k,l\in\mathbb{N}\colon (y,l)\neq(0,k)} m_k m_l T_{0,y}(G_k,G_l)+\sum_{k\in\mathbb{N}}(m_k^2+m_k) T_{0,0}(G_k,G_k)\\
& = \sum_{y\in \mathbb{Z}^d}\sum_{k,l\in\mathbb{N}} m_k m_l T_{0,y}(G_k,G_l)+\sum_{k\in\mathbb{N}}m_kt_k,
\end{align*}
where we used the independence of the Poisson point process of $\xi^{\ssup k}(x)$ in $k$ and $x$. Now, for the first summand, we carry out the summation over $y$ and recall that $\v=\sum_{y\in \Z^d}v(y)$ and $|G_k|=k$ and $|G_l|=l$. This gives 
\begin{align*}
    \mathbb{P}^m(\Phi_{0, \mathbb{Z}^d})=  \v \sum_{k,l\in\mathbb{N}} m_k m_l k l +\sum_{k\in\mathbb{N}}m_kt_k = \v \rho^2 + \sum_{k\in \N} m_kt_k,
\end{align*}
and hence, 
$$\varphi(m, \delta_0)\leq H(m|q)+\bar{v} \rho^2+ \sum_{k\in \N} m_kt_k.
$$
Using $t_k\leq k\v$ and minimising over $m$, the claim follows.
\end{proof}

\begin{step}$\chi(0)=H(0|q)=\sum_{k\in \N} q_k$ and $\chi'(0)=-\infty$.
\end{step}

\begin{proof}The first statement is clear since only the void process $P$ fits the constraint $P(N_0^\ssup{\ell})=0$. 
Next, using  Step \ref{step-uppboundchi} we see that
$$
\chi(\eps k)\leq \varphi(\eps \delta_k, \delta_0)\leq H(\eps \delta_k|q)+\v(\eps k)^2 +t_k \eps,\qquad k\in\N,\eps\in(0,1).
$$
Furthermore, $H(\eps \delta_k|q)=\sum_{l\in \N} q_l -\eps +\eps \log\frac{\eps}{q_k}$, so that 
$$
\frac{\chi(\eps k)-\chi(0)}{\eps k}\leq\frac {1}{k} \log\frac{\eps}{q_k}+O(1),\qquad \eps\downarrow0,
$$
which implies the second statement. 
\end{proof}

\begin{step} $\chi(\cdot,\cdot)$ is convex and in particular also $\rho\mapsto\chi(\rho)$ is convex. Further, $\chi(\cdot,\cdot)$ is  non-decreasing in $\rho_{\rm ma}$ and $\chi$ is continuous in $[0,\infty)^2$.
\end{step}
\begin{proof}
The convexity of $\chi$ follows from Step \ref{step-convex}. It implies that $\chi$ is continuous on the interior of its domain, i.e., on $(0,\infty)\times(0,\infty)$, and that $\rho\mapsto \chi(0,\rho)$ and $\rho\mapsto \chi(\rho,0)$ are continuous on $(0,\infty)$. The latter is also continuous at $0$ since $\liminf_{\rho\downarrow0}\chi(\rho)\ge \chi(0)=\sum_kq_k$. Indeed, note that $\chi(\rho)\ge \chi^{\ssup{v=0}}(\rho,0)$, which tends to $\chi^{\ssup{v=0}}(0,0)=\sum_kq_k=\chi(0)$ as $\rho\downarrow0$, as can be shown using standard variational calculus. This gives the desired result. 

From Step \ref{step-stochasnondec} we have that $\chi$ is non-decreasing in $\rho_{\rm ma}$. 
This implies continuity of $\chi$ at any point $(\rho,0)$ with $\rho\in[0,\infty)$, since it is upper semi-continuous there by convexity, and the lower semicontinuity follows from
$$\liminf_{(\rho_{\rm mi}, \rho_{\rm ma})\to (\rho, 0)} \chi(\rho_{\rm mi}, \rho_{\rm ma})\geq \liminf_{\rho_{\rm mi}\to \rho} \chi(\rho_{\rm mi},0)=\chi(\rho,0).
$$

For $\rho\in [0,\infty)$, it is evident that
$$
\chi(0,\rho)=\chi(0)+\inf_{\psi\in\Mcal_1(\N_0)\colon \sum_{a\in\N_0}a\psi(a)=\rho}\v a^2\psi(a).
$$
Now, for any $(\rho_{\rm mi}, \rho_{\rm ma})\in[0,\infty)^2$ we have, by dropping two of the three energy terms,
$$
\chi(\rho_{\rm mi}, \rho_{\rm ma})\geq \chi^{\ssup{v=0}}(\rho_{\rm mi})+\inf_{\psi\in\Mcal_1(\N_0)\colon \sum_{a\in\N_0}a\psi(a)=\rho_{\rm ma}}\v a^2\psi(a)=\chi^{\ssup{v=0}}(\rho_{\rm mi})-\chi(0)+\chi(0,\rho_{\rm ma}).
$$
We deduce that for $\rho\in[0,\infty)$,
$$
\liminf_{(\rho_{\rm mi}, \rho_{\rm ma})\to (0, \rho)} \chi(\rho_{\rm mi}, \rho_{\rm ma})\geq\chi(0, \rho).
$$
Hence the continuity of $\chi$ at $(0,\rho)$ follows, which finishes the proof of the continuity of $\chi$.
\end{proof}

Let us also point out that for fixed $m$, the map $ \rho_{\rm ma} \mapsto \inf_{\psi\in\Mcal_1(\N_0)\colon \sum_{a\in\N_0}a\psi(a)= \rho_{\rm ma}} \varphi(m,\psi)$ is convex and non-decreasing, and for fixed $\psi$, the map  $ \rho_{\rm mi}\mapsto \inf_{m\in [0,\infty)^\N,\sum_{k\in\N}km_k= \rho_{\rm mi}} \varphi(m,\psi)$ is convex. Also the version $\widetilde \chi$  of $\chi$ with the infimum ranging over compactly supported $m$ and $\psi$ (defined in Step~\ref{step:compactsupp} in Section~\ref{sec-lowbound}) is convex and continuous in each coordinate.

\begin{step}
If $2v(0)\geq \v$, i.e., $v(0)\geq \sum_{z\in \Z^d\setminus \{0\}} v(z)$, then $\varphi(m,\delta_0)\geq H(m|q)+\v \rho^2$ for any $\rho\in[0,\infty)$ and any $m\in[0,\infty)^\N$ satisfying $\sum_k km_k=\rho$. In particular, $\chi(\rho)\geq \chi^{\ssup{v=0}}(\rho,0)+\v\rho^2$.
\end{step}

\begin{proof}
Pick any $P\in\Mcal_1^{\ssup{\rm s}}(\Omega)$ satisfying $P(N_0^{\ssup{\delta_k}})=m_k$ for any $k\in\N$. In particular, $P(N_0^{\ssup\ell})=\rho =P(\widetilde N_0)$. Now we use that $P(\Phi_{0,\Z^d})=\sum_{z\in\Z^d}v(z) P(\widetilde N_0  \widetilde N_z)$, where $\widetilde N_x$ denotes the number of particles located at $x$. Indeed, using the shift-invariance of $P$, we see that
$$
\begin{aligned}
P(\Phi_{0,\Z^d})&=\sum_{k\in\N,i\in G_k}\sum_{z\in\Z^d}v(z)  P\big(N_0^{\ssup{\delta_k}}\widetilde N_{i+z}\big)
&=\sum_{k\in\N,i\in G_k}\sum_{z\in\Z^d}v(z)  P\big(N_{-i}^{\ssup{\delta_k}}\widetilde N_{z}\big)
=\sum_{z\in\Z^d}v(z) P(\widetilde N_0  \widetilde N_z).
\end{aligned}
$$

Hence, we have
$$
\begin{aligned}
P(\Phi_{0,\Z^d})&=\sum_{z\in\Z^d} v(z) P(\widetilde N_0\widetilde N_z)\\
&=v(0)P(\widetilde N_0^2)+\sum_{z\in\Z^d\setminus\{0\}}\frac{v(z)}2P\big((\widetilde N_0 + \widetilde N_z)^2-\widetilde N_0^2-\widetilde N_z^2\big)\\
&=\sum_{z\in\Z^d\setminus\{0\}}\frac{v(z)}2P\big((\widetilde N_0 + \widetilde N_z)^2\big)+\Big(v(0)-\sum_{z\in\Z^d\setminus\{0\}}v(z)\Big)P(\widetilde N_0^2).
\end{aligned}
$$
Now use the Cauchy--Schwarz inequality to estimate $P(\widetilde N_0^2)\geq P(\widetilde N_0)^2=\rho^2$ and $P((\widetilde N_0 + \widetilde N_z)^2)\geq P(\widetilde N_0 + \widetilde N_z)^2=4\rho^2$, to deduce that $P(\Phi_{0,\Z^d})\geq \v\rho^2$. Here we used our assumption on $v$. This implies that $\varphi(m,\delta_0)=H(m|q)+I_m(P)+\v\rho^2\geq H(m|q)+\v \rho^2$, using Step~\ref{step_3a} and the non-negativity of $I_m(P)$. Proceeding with the infimum over all $m$ satisfying $\sum_ k km_k=\rho$, we obtain that $\chi(\rho)\geq \chi^{\ssup{v=0}}(\rho)+\v\rho^2$.
\end{proof}

This finishes the proof of Lemma \ref{lem-properties}.

\subsection{Existence of minimising $\rho_{\rm mi}$ and $\rho_{\rm ma}$: proof of Lemma \ref{lem-Varformident}}\label{sec-VarformProof}
In this section we prove that \lq$\geq$\rq\ in \eqref{chiidentvarphi} holds; this implies that Lemma \ref{lem-Varformident} holds.

For this sake, fix a small threshold $\eps>0$ and pick $\rho_{\rm mi}$ and $\rho_{\rm ma}=\rho-\rho_{\rm mi}$ such that $\inf\chi\geq -\eps+ \chi(\rho_{\rm mi},\rho_{\rm ma})$. Then pick $\psi\in\Mcal_1(\N_0)$ satisfying $\sum_{a\in\N_0}a\psi(a)= \rho_{\rm ma}$, and $m=(m_k)_{k\in\N}$ satisfying $\sum_k k m_k= \rho_{\rm mi}$ such that
$$
\chi( \rho_{\rm mi},  \rho_{\rm ma})\geq -2\eps+\varphi(m,\psi).
$$
Note that, by Step \ref{step:compactsupp} in Section \ref{sec-lowbound}, we may and shall assume that $m$ and $\psi$ have compact supports contained in  $[K]=\{1,\dots,K\}$ and $\{0,\dots,A\}$, respectively, for some $K\in\N$ and $A\in\N_0$.

Using part of the proof of the lower bound of Theorem \ref{thm-freeenergygrid} (see Lemma~\ref{lem-varphirepr} or alternatively \eqref{lowboundZdeltavarphidelta}  and Step \ref{step-lowbounddeltato0} in Section \ref{sec-lowbound}), we see that
\begin{equation}
\begin{aligned}\label{varphilimitZNdelta}
\varphi(m,\psi)&\geq -\lim_{\delta\downarrow 0}\liminf_{N\to\infty}\frac1{|\Lambda_N|}\log\E_{\L_N}^{\ssup{\leq K}}\Big[\e^{-\Phi_{\Lambda_N,\Lambda_N}(\cdot+\omega_\psi)}\prod_{k\in[K]}\1\{N_{\Lambda_N}^{\ssup{\delta_k}}\in m_k|\Lambda_N|(1-\delta,1+\delta)\}\Big],
\end{aligned}
\end{equation}
where $\omega_\psi=\sum_{a=1}^A\delta_{(0,G_{N_a})}$, and $(N_a)_{a\in[A]}$ is any deterministic non-increasing sequence in $\N$ (depending on $N$) such that $\frac{N_a}{|\L_N|}\underset{N\to \infty} {\longrightarrow}\psi([a,\infty))$ for all $a\in[A]$. In particular, the marks in $\omega$ have size at most $K$, and $\sum_{a\in[A]}N_a \sim \rho_{\rm ma} |\Lambda_N|$. 

For simplicity, we assume now that $N_a=(2 L_a+1)^d$ are $d$-th powers of odd integers for $a\in[A]$ and we observe that $G_{N_a}=[-L_a,L_a]\cap\Z^d$. We pick now some $\widetilde K=(2L+1)^d>K$ with some integer $L$ and note that $G_{\widetilde K}=[-L,L]^d\cap\Z^d$. We also assume that $2L+1$ is a divisor of each of the numbers $2L_1+1,\dots,2L_A+1$.

Now  we replace the configuration $\omega+\omega_\psi$ by the configuration $\widetilde\omega=\omega+\sum_{a\in[A]}\sum_{x\in (2L+1)\Z^d\cap G_{N_a}}\delta_{(x,G_{\widetilde K})}$. In words, we re-organise all the particles in the $A$ boxes with cardinalities $N_1,\dots, N_A$ into a number of boxes of side length $2L$ without changing any of the locations of the particles. For any $a\in[A]$, these smaller boxes are mutually disjoint and their union is $G_{N_a}$.
In this way, we add to the PPP $\omega$ the marked point process with  $\widetilde N=\sum_{a\in[A]}N_a/\widetilde K$ points and mark $G_{\widetilde K}$ at each of these points. In particular, the energy remains unchanged, i.e., 
$$
\Phi_{\Lambda_N,\Lambda_N}(\omega+\omega_\psi)=\Phi_{\Lambda_N,\Lambda_N}(\widetilde \omega).
$$
Introducing $\widetilde m=(\widetilde m_k)_{k\in \N}$ by putting $\widetilde m_k=m_k$ for $k\in[K]$, and $\widetilde m_{\widetilde K}=\rho_{\rm ma}/\widetilde K$ and $\widetilde m_k=0$ for $k\in\N\setminus([K]\cup\{\widetilde K\})$, we have that $\sum_k k \widetilde m_k=\rho$. Note that $N_{\Lambda_N}^{\ssup{\delta_{\widetilde K}}}(\widetilde \omega)=\widetilde N = \sum_{a\in [A]} N_a/\widetilde K\in\widetilde m_{\widetilde K}|\Lambda_N|(1-\delta,1+\delta)$ for any large $N$. 

We now insert the configuration $\widetilde \omega$ instead of $\omega+\omega_\psi$ and conceive $\widetilde\omega$ as the random variable under $\P_{\L_N}^{\ssup{\leq \widetilde K}}$ instead of $\P_{\L_N}^{\ssup{\leq K}}$. We drop the fixation of the locations of the points with marks $G_{\widetilde K}$ and keep only the event that no mark of cardinalities $K+1,K+2,\dots,\widetilde K-1$ appears, precisely $\widetilde N$ points of cardinality $\widetilde K$ appear, located precisely at the mentioned locations. Denote by $p$ the probability of this. Then we upper estimate the indicator on this event by $\prod_{k=K+1}^{\widetilde K}\1\{N_{\L_N}^{\ssup{\delta_k}}\in \widetilde m_k|\L_N|(1-\delta,1+\delta)\}$. Therefore, we obtain that
\begin{equation}\label{phichicompare1}
\begin{aligned}
\E_{\L_N}^{\ssup{\leq K}}&\Big[\e^{-\Phi_{\Lambda_N,\Lambda_N}(\cdot+\omega_\psi)}\prod_{k\in[K]}\1\{N_{\Lambda_N}^{\ssup{\delta_k}})\in m_k|\Lambda_N|(1-\delta,1+\delta)\}\Big]\\
&\leq \E_{\L_N}^{\ssup {\leq \widetilde K}}\Big[\e^{-\Phi_{\Lambda_N,\Lambda_N}}\prod_{k\in[\widetilde K]}\1\{N_{\Lambda_N}^{\ssup{\delta_k}}\in  \widetilde m_k|\Lambda_N|(1-\delta,1+\delta)\}\Big]\frac 1 p,
\end{aligned}
\end{equation}
where
$$
p=\Big(\prod_{k=K+1}^{\widetilde K-1}\Poi_{q_k |\L_N|}(0)\Big)\Poi_{q_{\widetilde K}|\Lambda_N|}(\widetilde N)\, \frac{\widetilde N!}{|\Lambda_N|^{\widetilde N}}\geq \exp\Big\{-|\L_N|\sum_{k=K+1}^{\widetilde K}q_k\Big\}q_{\widetilde K}^{|\L_N|(\rho_{\rm ma}/\widetilde K + o(1))},
$$
where we recall that $N_{\L_N}^{\ssup{\delta_k}}$ is Poisson-distributed with parameter $q_k|\L_N|$, and the  quotient is the probability to put these $\widetilde N$ points at particular places in $\Lambda_N$. This shows that, for all sufficiently large $N$,
$$
\frac 1p\leq  \exp\Big\{|\Lambda_N|\Big[\sum_{k>K}q_k+o(1)+
\frac {\rho_{\rm ma}}{\widetilde K} \log q_{\widetilde K}\Big\}\leq \e^{|\Lambda_N| \eta_{K,\widetilde K}},
$$
with some $\eta_{K,\widetilde K}>0$ that vanishes as $K,\widetilde K\to\infty$, since $q_k=e^{o(k)}$ as $k\to\infty$. Now we apply Lemma~\ref{lem-varphirepr} (note that we do not have to make $K\to\infty$ nor $A\to\infty$ for $m$ and $\psi$ having compact supports) to see that
$$
\limsup_{\delta\downarrow0}\limsup_{N\to\infty}\frac 1{|\Lambda_N|}\log(\mbox{r.h.s.\ of \eqref{phichicompare1}})\leq -\varphi(\widetilde m,\delta_0)+\eta_{K,\widetilde K}.
$$
Now recall that $\sum_{k\in\N}k \widetilde m_k=\rho$ and hence $(\widetilde m, 0)$ is admissible in the variational formula for $\chi(\rho,0)$, so we have 
$$
\inf\chi\geq -2\eps+\varphi(m,\psi)\geq -2\eps+\varphi(\widetilde m,\delta_0)-\eta_{K,\widetilde K}\geq -2\eps+\chi(\rho,0)-\eta_{K,\widetilde K}.
$$
Taking $K,\widetilde K\to\infty$ and $\eps\downarrow0$, we get $\inf\chi\geq \chi(\rho,0$), which finishes the proof of Lemma \ref{lem-Varformident}.

\subsection{Minimisers of $\varphi$: proof of Lemma \ref{lem-ProofExistMin}}\label{sec-minimivarphi}
Recall from the end of Section~\ref{sec-uppbound} that we proved already the existence of a minimiser $(m,\psi)$ of $\varphi$ with density $\rho=\sum_{a=0}^\infty a\psi(a)+\sum_{k=1}^\infty km_k$.

\setcounter{step}{0}

    \begin{step} For any $(m,\psi)$, there exists at least one minimising family $(P_a)_{a\in\N_0}$ for the variational formula \eqref{varphidef} in the definition of $\varphi(m,\psi)$.
    \end{step}
    
\begin{proof}    
Note that the level sets of $I$ are compact in the local tame topology and the map $P\mapsto P(N_0^{\ssup{\delta_k}})$ is continuous in this topology, and the three other functionals in the first line of \eqref{varphidef} are lower semi-continuous. The difficulty now lies in the fact that the mapping $(P_a)_{a\in \N_0}\mapsto \sum_{a=0}^\infty \psi(a)P_a(N_0^{\ssup{\delta_k}})$ is a priori only lower semi-continuous. However, for any $A>0$,
\begin{equation}\label{eq_minP}
\sum_{a>A}\psi(a)P_a(N_0^{\ssup{\delta_k}})\le \frac{1}{2\v A}\sum_{a=0}^\infty\psi(a)2\v aP_a(N_0^{\ssup{\ell}})\le \frac{1}{2\v A}\Xi\big((P_a)_{a}\big),
\end{equation}
where $\Xi\big((P_a)_{a}\big):=\sum_{a=0}^\infty\psi(a)\big[ I(P_a)+2\v aP_a(N_0^{\ssup{\ell}})+a^2\v\big]$. Now, assume that $(P^{\ssup{n}}_a)_a$ is an approximate minimiser for $\varphi(m,\psi)$, where $\sum_{a=0}^\infty\psi(a)P^{\ssup{n}}_a(N_0^{\ssup{\delta_k}})=m_k$ for all $k,n$. That is, $\lim_{n\to\infty}\Xi\big((P^{\ssup{n}}_a)_{a}\big)=\varphi(m,\psi)$. Then, $(I(P^{\ssup{n}}_a))_n$ is bounded for any $a$. Hence, we have a subsequence, which we call also $P^{\ssup{n}}$ for convenience, such that $\lim_{n\to\infty}P^{\ssup{n}}_a=P_a$ for any $a$. By lower semicontinuity of $\Xi$, we have that $\Xi\big((P_a)_{a}\big)\le \varphi(m,\psi)$. In particular, for any $k,a$, we have that $\lim_{n\to\infty}P^{\ssup{n}}_a(N_0^{\ssup{k}})=P_a(N_0^{\ssup{k}})$. Now, it is easy to see, using \eqref{eq_minP}, that for any $k$ we have that $\sum_{a=0}^\infty\psi(a)P_a(N_0^{\ssup{k}})=m_k$. Hence, $(P_a)_a$ is indeed a minimiser in the formula~\eqref{varphidef}. 
\end{proof}

\begin{step}\label{step-chiinequality}
$$
\chi(\rho+a)\leq \chi(\rho,a)\leq \chi(\rho)+2\v a \rho+\v a^2,\qquad \rho\in[0,\infty),a\in \N_0.
$$ 
\end{step}
\begin{proof} The first inequality comes from Lemma \ref{lem-Varformident}. For the second one note that for any $m$ such that $\sum_{k\in \N} k m_k=\rho$, we have $\chi(\rho,a)\leq \varphi(m,\delta_a)=\varphi(m,\delta_0)+2\v a \rho+\v a^2$, according to the definitions. Now take the infimum over $m$.
\end{proof}

Now we turn to a closer analysis of the crucial variational problem in \eqref{crucialVP}. Recall that Lemma \ref{lem-ProofExistMin} claims the existence of a minimiser; it was proved at the end of Section \ref{sec-proofgrid}. We say that $(m,\psi)$ has density $\rho$ if $\sum_{k\in\N}k m_k+ \sum_{a\in\N_0}a\psi(a)=\rho$.

\begin{step}\label{step-twoatoms}
For any $\rho\in[0,\infty)$, there is a minimiser $(m,\psi)$ of $\varphi$ with density $\rho$ such that $\psi$ has at most two atoms. If $\chi(\cdot)$ is strictly convex at $\rho$, then there is even a minimiser such that $\psi$ has precisely one atom.
\end{step}

\begin{proof}
Fix  $\rho\geq 0$ and a minimiser $(m, \psi)$ of $\varphi$ with density $\rho$. Consider an associated minimising family of associated processes $(P_a)_{a\in \mathbb{N}_0}$ and denote $m^a=(P_a(N_0^\ssup{\delta_k}))_{k\in\N}$ and $\rho_a=\rho_{\rm mi}(m^a)=P_a(N_0^{\ssup \ell})$.
We have that $\rho=\sum_{a\in \mathbb{N}_0} \psi(a) (\rho_a+a)$.

Now use Step \ref{step-chiinequality} and the definition of $\chi$ to see that
\begin{equation}\label{chiinequality}
\begin{aligned}
\chi(\rho_a+a)
& \leq \chi(\rho_a)+2 \v a \rho_a+\v a^2 \\
& \leq  \varphi(m^a, \delta_a)=\varphi(m^a, \delta_0)+2 \v a \rho_a+\v a^2\\
& \leq I(P_a)+P_a(\Phi_{0,\mathbb{Z}^d})+2 \v a \rho_a+\v a^2,\qquad a\in \N_0.
\end{aligned}
\end{equation}
Furthermore, use the minimality of $(P_a)_{a\in\N_0}$ and of $(m,\psi)$ and then the convexity of $\chi(\cdot)$ to see that
\begin{equation}\label{equasandwich}
    \sum_{a \in \mathbb{N}_0} \psi(a) \left[I(P_a)+P_a(\Phi_{0,\mathbb{Z}^d})+ 2  \v a \rho_a +  \v a^2\right]  = \varphi(m,\psi)=\chi(\rho)
 \leq \sum_{a \in \mathbb{N}_0} \psi(a) \chi(\rho_a+a).
\end{equation}
In view of \eqref{chiinequality},  there is in fact equality everywhere in \eqref{equasandwich}. Using this and \eqref{chiinequality}, we get for $a\in \supp(\psi)$:
\begin{equation}\label{descrminimisereq}
\chi(\rho_a+a)=\chi(\rho_a)+2 \bar{v} a \rho_a+\bar{v} a^2,
\end{equation}
and $m^a$ is a minimiser in the definition of $\chi(\rho_a, 0)$.
Moreover, we have the equality $\chi(\rho)= \sum_{a \in \mathbb{N}_0} \psi(a) \chi(\rho_a+a)$. 

Now we treat the cases of strict and non-strict convexity of $\chi(\cdot)$ at $\rho$ separately:
\begin{itemize}
    \item If $\chi(\cdot)$ is strictly convex at $\rho$, then the equality  $\chi(\rho)= \sum_{a \in \mathbb{N}_0} \psi(a) \chi(\rho_a+a)$ in the convexity inequality implies that $\rho_a+a=\rho$ for any $a\in \text{supp}(\psi)$. Therefore, for any such $a$, \eqref{descrminimisereq} shows that $(m^a, \delta_a)$ is a minimiser of $\varphi$ with density $\rho$.
    
    \item If $\chi(\cdot)$ is not strictly convex at $\rho$, then it is affine on a non-trivial interval containing $\rho$ in its interior. We denote by $I_{\text{aff}}(\rho)$ the biggest such interval. Then $\rho_a+a\in I_{\text{aff}}(\rho)$ for all $a\in \text{supp}(\psi)$. Recall that $\rho=\sum_{a\in \mathbb{N}_0} \psi(a) (\rho_a+a)$, so we can pick $a_1$ and $a_2$ in $\supp(\psi)$ such that $\rho_{a_1}+a_1\leq \rho\leq \rho_{a_2}+a_2$. Then we can build a minimiser $(m,\psi)$ with density $\rho$ by taking a suitable convex combination of $(m^{a_1},\delta_{a_1})$ and $(m^{a_2},\delta_{a_2})$, using affinity of $\chi(\cdot)$ on $I_{\text{aff}}(\rho)$. Then $\psi$ has no more than two atoms.
\end{itemize}

\end{proof}

As a complement, let us generalise the above result to give a full description of the minimisers of $\varphi$ at given density.

\begin{step}
Description of all the minimisers $(m,\psi)$ of $\varphi$.
\end{step}

First we consider the case that $\chi(\cdot)$ is strictly convex at $\rho$. Consider $a\in\N$. If $a\leq \rho$ and $\chi(\rho-a)+ 2 \v a (\rho-a) + \v a^2=\chi(\rho)$ denote by $M(\rho,a)$ the set of all the minimisers of $\varphi(\cdot, 0)$ at density $\rho-a$. Otherwise, put $M(\rho,a)=\emptyset$. By convexity of $\varphi$, $M(\rho,a)$ is convex (possibly a singleton or empty). Then it is easy to see from the proof of Step \ref{step-twoatoms} that the convex hull of the set of all the $(m,\delta_a)$ with $a\in\N_0$ and $m\in M(\rho,a)$ is equal to the set of minimisers of $\varphi$ with density $\rho$, i.e., any minimiser with density $\rho$ is a convex combination of over $a$ of such pairs, and conversely any convex combination aver $a$ of such pairs is a minimiser with density $\rho$.

Now we consider the case that $\chi(\cdot)$ is affine on a non-trivial segment $J=[ \rho_1,  \rho_2]$ and this segment is maximal for this property. Consider, for fixed $a\in\N_0$, the set $\widetilde M(J,a)$ of all the $m\in[0, \infty)^\N$ such that $\rho_{\rm mi}(m)+a\in J$, $\varphi(m, \delta_0)=\chi(\rho_{\rm mi}(m))$ and $\chi(\rho_{\rm mi}(m))+ 2 \v a \rho_{\rm mi}(m) + \v a^2=\chi(\rho_{\rm mi}(m)+a)$. Again by convexity of $\varphi$, every $\widetilde M(J, a)$ is convex. The densities of the pairs $(m,\delta_a)$ lie in $J$, on which $\chi(\cdot)$ is affine.
Then, for any $\rho\in [ \rho_1,  \rho_2]$, any minimiser of $\varphi$ with density $\rho$  is a  convex combination over $a$ of the pairs $(m,\delta_a)$ with $m\in \widetilde M(a,J)$, such that the resulting average density is $\rho$. Conversely, any combination over $a$ of such pairs is a minimiser at density the average density of the combination.
  
\section{Differentiability and phase transition}\label{sec-ProofBEC}

We prove the differentiability of $\varphi(\cdot,\psi)$ (Lemma \ref{lem-DerivativeMicro}) in Section \ref{sec-varphidiffble}, the Euler--Lagrange equation (Lemma \ref{lem-ELeq}) and the differentiability of $\chi(\cdot)$ (Corollary \ref{cor-chidiff}) in Section \ref{sec-EL-eqs}, we give an alternative formula for the derivative of $\varphi(\cdot,\psi)$ in Section~\ref{sec-thinning}, and we discuss the nature of the possible phase transition (Lemma \ref{lem-saturation}) in Section \ref{sec-naturephasetrans}.

\subsection{Differentiability of $\varphi$: proof of Lemma \ref{lem-DerivativeMicro}}\label{sec-varphidiffble}

In this section, we give the proof of Lemma \ref{lem-DerivativeMicro}, i.e., of the differentiability of the map $m\mapsto\varphi(m,\psi)$ for fixed $\psi\in\Mcal_1( \N_0)$. We are not going to carry out this proof using the variational formula, but we will be starting from the characterisation of $\varphi(m,\psi)$ in Lemma~\ref{lem-varphirepr} in terms of the exponential rate of a restricted partition function. We fix $m\in[0,\infty)^\N$ an $\psi\in\Mcal_1(\N_0)$.

\setcounter{step}{0}

\begin{step} For any $k\in\N$ and small $\eps>0$, we derive an alternate approximate variational formula for $\varphi(m+\eps\delta_k,\psi)$.
\end{step}

\begin{proof}
We put $ \rho_{\rm ma}=\sum_{a\in\N_0}a \psi(a)$ and $\rho=\sum_a a\psi(a)+\sum_{l}l m_l$. Take $\eps>0$ and put $m(\eps)=m+\eps\delta_k$. We apply \eqref{widetildevarphiformula} for $m(\eps)$ instead of $m$ and with $\rho+\eps k$ instead of $\rho$. On the right-hand side, we replace $N$ by $N(\eps)=N(1+\eps k/\rho)$ and  note that the box  $\L_N$ is the same box with $(N,\rho)$ replaced by $(N(\eps),\rho+\eps k)$. 

We now derive an alternative variational formula as an upper bound for the right-hand side of \eqref{widetildevarphiformula} by explicitly carrying out the integration over the $\eps |\L_N|$ additional Poisson points with mark $G_k$ and describing their influence on the expectation as a functional of the empirical stationary distribution of the Poisson process. (For better readability, we drop the integer-part brackets $\lfloor\cdot\rfloor$  in the following.) For this, we are going to use that the Poisson points of the reference process $\omega_{\rm P}^{\ssup k}$, given their number, are i.i.d.~uniformly over $\L_N$ distributed sites. 

Indeed, assuming that $K>k$ and $\delta<m_k/2$, in \eqref{ZNmpsidef} for $Z_{N(\eps),\L_N}^{\ssup{A,K,\delta}}(m(\eps),\psi)$, we carry out the expectation with respect to $\P^{\ssup k}$ by first taking the Poisson probability $\Poi_{q_k|\L_N|}(l_k)$ for having $l_k\in[m_k(\eps)-\delta,m_k(\eps)+\delta]|\L_N|$ Poisson points with mark $G_k$, then handling $l_k-\eps \L_N$ of them again as the number of Poisson points in the process $\omega_{\rm P}^{\ssup k}$ and treating the remaining $\eps|\L_N|$ of them as i.i.d.~uniformly over $\L_N$ distributed sites $X_1,\dots,X_{\eps|\L_N|}$, each of which carries the mark $G_k$. We write $\omega_{\rm U}^{\ssup k}=\sum_{i=1}^{\eps|\Lambda_N|}\delta_{(X_i,G_k)}$ for the arising marked random point process under the measure $\Ucal_{\L_N}^{\otimes \eps\L_N}$, where we recall that $\Ucal_{\L_N}$ denotes the uniform distribution on $\L_N$. This implies that
\begin{equation}\label{Zperturbedident}
\begin{aligned}
Z_{N(\eps),\L_N}^{\ssup{A,K,\delta}}(m(\eps),\psi)&= \E_{\L_N}^{\ssup{\leq K}}\otimes \Ucal_{\L_N}^{\otimes \eps|\L_N|}\Big[\e^{-\Phi_{\Lambda_N,\Lambda_N}(\cdot+\omega_{\psi^{\ssup A}}+\omega^{\ssup k}_{\rm U})}\prod_{j=1}^K\1\{|N_{\L_N}^{\ssup{\delta_j}}-m_j |\L_N||\leq \delta|\L_N|\}\Big]\\
&\quad \times 
\sum_{l_k= [m_k(\eps)-\delta]|\L_N|}^{[m_k(\eps)+\delta]|\L_N|}\frac{\Poi_{q_k|\L_N|}(l_k)}{\Poi_{q_k|\L_N|}(l_k-\eps\L_N|)}.
\end{aligned}
\end{equation}
The large-$N$ asymptotics of last term is easily identified. Indeed, for any summand $l_k$, 
\begin{equation}\label{Poissonprobeps}
\begin{aligned}
\limsup_{\delta\downarrow0}\limsup_{N \to\infty}-\frac 1 {|\L_N |}\log\frac{\Poi_{q_k|\L_N|}(l_k)}{\Poi_{q_k|\L_N|}(l_k-\eps|\L_N|)}
&\leq \eps\Big[\log\frac{m_k+\eps}{q_k}+\frac {m_k}\eps\log\frac{m_k+\eps}{m_k}-1\Big]\\
&\sim \eps\log\frac{m_k}{q_k},\qquad \eps\downarrow0,
\end{aligned}
\end{equation}
and an analogous estimate is derived for the limit inferior. This explains the first term on the right-hand side of \eqref{DerivativeMicro_Eq}.

Next, we  split the energy according to the contributions from $\omega_{\rm P}$ and $\omega^{\ssup k}_{\rm U}$ and the remainder. For this, we extend our notation for the energy from a self-energy to a mutual energy by putting
$$
\Phi_{\L,\L}^{\ssup{\leftrightarrow}}(\omega,\widetilde\omega)=\sum_{x,y\in \L}\sum_{l,\widetilde l\in\N}\xi^{\ssup l}(x)\widetilde\xi^{\ssup {\widetilde l}}(y)T_{x,y}(G_l,G_{\widetilde l}),\qquad\omega,\widetilde\omega\in\Omega,
$$
where $\xi^{\ssup l}(x)$ is the number of $\omega$-points in $x$ having mark $G_l$, and $\widetilde \xi^{\ssup l}(x)$ is the number of $\widetilde \omega$-points in $x$ having mark $G_l$. Note that $\Phi_{\L,\L}^{\ssup{\leftrightarrow}}$ is linear in each of the two arguments and that $\Phi_{\L,\L}(\omega+\widetilde\omega)=\Phi_{\L,\L}(\omega)+2\Phi_{\L,\L}^{\ssup{\leftrightarrow}}(\omega,\widetilde\omega)+\Phi_{\L,\L}(\widetilde \omega)$ for any $\omega,\widetilde\omega\in\Omega$. Then we see that
\begin{equation}\label{mutualenergy}
\begin{aligned}
\Phi_{\Lambda_N,\Lambda_N}&\big(\omega+\omega_{\psi^{\ssup A}}+ \omega_{\rm U}^{\ssup k}\big)\\
&= \Phi_{\Lambda_N,\Lambda_N}\big(\omega+\omega_{\psi^{\ssup A}}\big)+ \Phi_{\Lambda_N,\Lambda_N}(\omega_{\rm U}^{\ssup k})+ 2\Phi^{\ssup{\leftrightarrow}}_{\Lambda_N,\Lambda_N}\big(\omega+\omega_{\psi^{\ssup A}},\omega_{\rm U}^{\ssup k}\big).
\end{aligned}
\end{equation}
The first term on the right-hand side is equal to the energy of the original, unperturbed configuration.  The one-but-last term is equal to the sum of the internal energies of the $\eps|\L_N|$ marks $G_k$ at the sites $X_1,\dots,X_{\eps|\L_N|}$ (i.e., $t_k \eps |\L_N|$ where we recall that $t_k= T_{0,0}(G_k,G_k)$) plus the mutual interaction between any two of these marked points $\delta_{(X_i,G_k)}$; that is,
\begin{equation}\label{selfinterbounds}
\begin{aligned}
t_k \eps|\L_N|\le \Phi_{\Lambda_N,\Lambda_N}(\omega_{\rm U}^{\ssup k})&\leq t_k \eps|\L_N|+ 2\sum_{1\leq i<j\leq \eps|\L_N|}\Phi^{\ssup{\leftrightarrow}}_{\Lambda_N,\Lambda_N}(\delta_{(X_i,G_k)},\delta_{(X_j,G_k)})\\
&\leq t_k \eps|\L_N|+C_{k,v}\sum_{1\leq i<j\leq \eps|\L_N|} \1\{|X_i-X_j|\leq d_k\},
\end{aligned}
\end{equation}
with some constants $C_{k,v}$ and $d_k$, depending only on $k$, $v$ and  the dimension $d$. 

  By $\L_N^a$ we denote the part of $\L_N$ in which precisely $a$ grids overlap. Let us introduce the volume $\widetilde \L=\{x\in \L\colon x+G_k+\supp(v)\subset\L\}\subset\L$ containing all points in $\L$ such that their $k$-marks do not interact with $\L^{\rm c}$. Then, on the event $\{X_1,\dots,X_{\eps|\L_N|}\in\widetilde\L_N\}$, the last term of the right-hand side of \eqref{mutualenergy} can be expressed as
\begin{equation}\label{mutualenergy2}
\begin{aligned}
\Phi^{\ssup{\leftrightarrow}}_{\Lambda_N,\Lambda_N}\big(\omega+\omega_{\psi^{\ssup A}},\omega_{\rm U}^{\ssup k}\big)&=\sum_{i=1}^{\eps|\L_N|}f_{\L_N}(X_i), \mbox{ with }\\
f_\L(x)&=\Phi^{\ssup{\leftrightarrow}}_{\Lambda,\Lambda}\big(\theta_{x}(\omega),\delta_{(0,G_k)}\big)+\v k\sum_a a \1\{x\in \L^a\},
\end{aligned}
\end{equation}
where we recall the shift operator $\theta_x$ by $x\in\Z^d$. We note that, for each $i\in\{1,\dots,\eps|\L_N|\}$, each of the $k$ particles of $\delta_{(X_i,G_k)}$ has the interaction $\v$ with each of the $a$ grids in the area $\L_N^a$ in which precisely $a$ of the macroscopic grids overlap, for $a\in\N_0$. Let us examine the expectation over $\omega_{\rm U}^{\ssup k}$ in \eqref{Zperturbedident}, conditional on $\omega+\omega_{\psi^{\ssup A}}$. 

We start by deriving a lower bound. We first claim that 
\begin{equation}\label{omegaUintegrate}
\frac{\Ucal_{\L_N}^{\otimes \eps|\L_N|}\big[\e^{ - 2\Phi^{\ssup{\leftrightarrow}}_{\Lambda_N,\Lambda_N}(\omega+\omega_{\psi^{\ssup A}},\omega_{\rm U}^{\ssup k})-\Phi_{\Lambda_N,\Lambda_N}(\omega_{\rm U}^{\ssup k})}\big]}
{\Ucal_{\widetilde\L_N}[\e^{-2 f_{\L_N}}]^{\eps |\L_N|}\e^{- t_k \eps|\L_N|}}\geq\e^{-\eps|\L_N|( \eps C+\log(|\widetilde\L_N|/|\L_N|))},
\end{equation}
where $C$ does not depend on $\omega,\omega_{\psi^{\ssup A}},\L_N,A$ and $\eps$. To see this, we start by inserting the indicator on the event $\{X_1,\dots,X_{\eps|\L_N|}\in\widetilde\L_N\}$ and use~\eqref{mutualenergy2} to bound
$$
\begin{aligned}
\Ucal_{\L_N}^{\otimes \eps|\L_N|}&\big[\e^{ - 2\Phi^{\ssup{\leftrightarrow}}_{\Lambda_N,\Lambda_N}(\omega+\omega_{\psi^{\ssup A}},\omega_{\rm U}^{\ssup k})-\Phi_{\Lambda_N,\Lambda_N}(\omega_{\rm U}^{\ssup k})}\big]\ge \Big(\tfrac{|\widetilde\L_N|}{|\L_N|}\Big)^{\eps|\L_N|}\Ucal_{\widetilde\L_N}^{\otimes \eps|\L_N|}\big[\e^{ -2\sum_{i=1}^{\eps|\L_N|}f_{\L_N}(X_i)-\Phi_{\Lambda_N,\Lambda_N}(\omega_{\rm U}^{\ssup k})}\big]\\
&\ge \Big(\tfrac{|\widetilde\L_N|}{|\L_N|}\Big)^{\eps|\L_N|}\e^{- t_k \eps|\L_N|}\Ucal_{\widetilde\L_N}^{\otimes \eps|\L_N|}\big[\e^{ -2\sum_{i=1}^{\eps|\L_N|}f_{\L_N}(X_i)-C_{k,v}\sum_{1\leq i<j\leq \eps|\L_N|} \1\{|X_i-X_j|\leq d_k\}}\big],
\end{aligned}
$$
where we also used \eqref{selfinterbounds}. Next, we use Jensen's inequality to see that
$$
\begin{aligned}
\frac{\Ucal_{\widetilde\L_N}^{\otimes \eps|\L_N|}\big[\e^{-2\sum_{i=1}^{\eps|\L_N|}f_{\L_N}(X_i)}\e^{-C_{k,v}\sum_{i<j}\1\{|X_i-X_j|\leq d_k\}
}\big]}
{\Ucal_{\widetilde\L_N}^{\otimes \eps|\L_N|}\big[\e^{-2\sum_{i=1}^{\eps|\L_N|}f_{\L_N}(X_i)}\big]}
&\geq \exp\Big\{-C_{k,v}\sum_{i<j}\widehat \Ucal_N(|X_i-X_j|\leq d_k)\Big\}\\
&\geq \exp\big\{-C_{k,v}\eps^2|\L_N|^2 \widehat \Ucal_N(|X_1-X_2|\leq d_k)\big\},
\end{aligned}
$$
where we wrote $\widehat \Ucal_N$ for the measure with density $\e^{-2\sum_{i=1}^{\eps|\L_N|}f_{\L_N}(X_i)}$ with respect to $\Ucal_{\widetilde\L_N}^{\otimes \eps|\L_N|}$, properly normalised. In the last term, $\widehat \Ucal_N(|X_1-X_2|\leq d_k)$, using the product structure of this measure, one can carry out the integration with respect to $X_3,X_4,...,X_{\eps |\L_N|}$ in the numerator and in the denominator, and they cancel each other. Given $X_1$, the integration over $X_2$ is limited to the $d_k$-box around $X_1$, which implies that $\widehat \Ucal_N(|X_1-X_2|\leq d_k)$ is of order $1/|\L_N|$ as $N\to\infty$. Indeed, by dropping $f_{\L_N}(X_2)$ in the exponent in the numerator,
$$
\begin{aligned}
\widehat \Ucal_N(|X_1-X_2|\leq d_k)
 &=\frac{\Ucal_{\widetilde\L_N}^{\otimes 2}\big(\e^{-2f_{\L_N}(X_1)-2f_{\L_N}(X_2)}\1\{|X_1-X_2|\leq d_k\}\big)}{\Ucal_{\widetilde\L_N}(\e^{-2f_{\L_N}(X)})^2} \\
 &\leq \frac{d_k^d}{|\widetilde\L_N|}\Ucal_{\widetilde\L_N}\big(\e^{-2f_{\L_N}(X)}\big)^{-1}\leq \frac{d_k^d}{|\widetilde\L_N|}\exp\big\{2\, \Ucal_{\widetilde\L_N}(f_{\L_N}(X))\big\},
\end{aligned}
$$
where we again used Jensen's inequality in the last step. Finally, using the fact that the particle number in $\omega$ is subject to a constraint in~\eqref{Zperturbedident}, we have that 
$$
\begin{aligned}
\Ucal_{\widetilde\L_N}(f_{\L_N}(X))&=|\widetilde\L_N|^{-1}\sum_{x\in\widetilde\L_N}\Phi^{\ssup{\leftrightarrow}}_{\Lambda_N,\Lambda_N}\big(\theta_{x}(\omega),\delta_{(0,G_k)}\big)+\v k\sum_a a\psi^{\ssup A}(a)\\
&\le \v k\sum_{l=1}^K l (m_l+\delta)+\v k\sum_a a\psi(a)\le \v k (\rho+\delta K^2).
\end{aligned}
$$
Thus we have proved~\eqref{omegaUintegrate}.

Next, we introduce the empirical stationary field of a point process $\omega$ with parameters $a\in\{0,\dots,A\}$:
$$
\Rcal^{\ssup a}_{\L_N}(\omega)=|\L^a_N|^{-1}\sum_{x\in \L^a_N}\delta_{\theta_x(\omega)}\qquad\mbox{and}\qquad \Rcal_{\L_N}=\sum_{a=0}^A\frac{|\L^a_N|}{|\L_N|}\Rcal^{\ssup a}_{\L_N}.
$$ 
Then, we can represent 
\begin{align*}
\Ucal_{\widetilde \L_N}(\e^{-2 f_{\L_N}(X)})= \sum_{b=0}^A\frac{|\widetilde\L^a_N|}{|\widetilde\L_N|}\e^{-2\v ak}\big\langle \Rcal^{\ssup a}_{\widetilde\L_N}(\omega), \e^{- 2\Phi^{\ssup k}_{\Lambda_N}}\big\rangle, 
\end{align*}
where $\Phi^{\ssup k}_{\Lambda_N}(\omega)=\Phi^{\ssup{\leftrightarrow}}_{\Lambda_N,\Lambda_N}(\omega,\delta_{(0,G_k)})$ denotes the interaction of a marked configuration $\omega$ in $\Lambda_N$ with the origin, carrying a mark $G_k$. 
In particular, since there is a maximal size of the marks, the mapping $P\mapsto \big\langle P, \e^{- 2\Phi^{\ssup k}}\big\rangle$ is continuous. From now on we can carry out the same steps as we did in the proof of Theorem~\ref{thm-freeenergygrid} for performing the limit as $N\to\infty$, followed by $\delta\downarrow 0$, which implies
\begin{equation}\label{varphiderivative}
\begin{aligned}
&\liminf_{\delta\downarrow 0}\liminf_{N\to\infty}\frac 1{|\L_N|}\log Z_{N(\epsilon),\L_N}^{\ssup{A,K,\delta}}(m(\epsilon),\psi)\ge-\eps\log\frac{m_k}{q_k}-\eps t_k-\varphi_{K,A,\eps}(m,\psi)-O(\eps^2), 
\end{aligned}
\end{equation}
where $O(\eps^2)$ is independent of $A$, we used  $\lim_{N\to\infty}|\L_N|/|\widetilde\L_N|=1$, and 
\begin{equation}\label{varpiAKdeltadef}
\begin{aligned}
\varphi_{K,A,\eps}&(m,\psi)=\inf\Big\{\sum_{a=0}^A\psi^{\ssup A}(a)\Big[I^{\ssup {\leq K}}(P_a)+P_a(\Phi_{0,\Z^d})+2 {\v} aP_a(N_0^{\ssup \ell})+{\v}a^2\Big]\\
&-\eps \log\Big(\sum_{a=0}^A \psi^{\ssup A}(a)\e^{-2\v ak}P_a(\e^{-2\Phi^{\ssup k}})\Big)\colon P_0,P_1,P_2,\dots\in  \Mcal_1^{\ssup{\rm s}}(\Omega^{\ssup {\leq K}}), \\ 
&\qquad\qquad\sum_{a=0}^A\psi^{\ssup A}(a)P_a(N_0^{\ssup {\delta_l}})=m_l\,\forall l\le K\Big\},\quad A,K\in \N\cup\{\infty\}.
\end{aligned}
\end{equation}
Further, for $K\to\infty$ and $A\to \infty$, $\varphi_{K,A,\eps}(m,\psi)$ converges to $\varphi_{\infty,\infty,\eps}(m,\psi)$.

Before we make the connection to $\varphi(m+\eps\delta_k,\psi)$, we first derive also an upper bound for $Z_{N(\epsilon),\L_N}^{\ssup{A,K,\delta}}(m(\epsilon),\psi)$. For this, instead of \eqref{omegaUintegrate} we claim that 
\begin{equation}\label{omegaUintegrate_1}
\frac{\Ucal_{\L_N}^{\otimes \eps|\L_N|}\big[\e^{ - 2\Phi^{\ssup{\leftrightarrow}}_{\Lambda_N,\Lambda_N}(\omega+\omega_{\psi^{\ssup A}},\omega_{\rm U}^{\ssup k})-\Phi_{\Lambda_N,\Lambda_N}(\omega_{\rm U}^{\ssup k})}\big]}
{\Ucal_{\L_N}[\e^{-2 \tilde f_{\L_N}}]^{\eps |\L_N|}\e^{- t_k \eps|\L_N|}}\leq1,
\end{equation}
where $\tilde f_\L(x)=\one\{x\in \widetilde \L\}\Phi^{\ssup{\leftrightarrow}}_{\Lambda,\Lambda}\big(\theta_{x}(\omega),\delta_{(0,G_k)}\big)+\v k\sum_a a \1\{x\in \L^a\}$. Indeed, this follows from the lower bound in \eqref{selfinterbounds} and the fact that $\Phi^{\ssup{\leftrightarrow}}_{\Lambda_N,\Lambda_N}\big(\omega+\omega_{\psi^{\ssup A}},\omega_{\rm U}^{\ssup k}\big)\ge\sum_{i=1}^{\eps|\L_N|}\tilde f_{\L_N}(X_i)$. 
Then, using the same steps as in the lower bound, with slide changes in the approximations, we arrive at
\begin{equation}\label{varphiderivative_1}
\begin{aligned}
&\limsup_{K,L\to\infty}\limsup_{\delta\downarrow 0}\limsup_{N\to\infty}\frac 1{|\L_N|}\log Z_{N(\epsilon),\L_N}^{\ssup{A,K,\delta}}(m(\epsilon),\psi)\le -\eps\log\frac{m_k}{q_k}-\eps t_k-\varphi_{\infty,\infty,\eps}(m,\psi). 
\end{aligned}
\end{equation}

Finally, recall that, according to Lemma~\ref{lem-varphirepr}, the left-hand side is equal to $-\varphi(m+\eps\delta_k,\psi)$.
Hence, combining the lower and upper bounds in \eqref{varphiderivative} and \eqref{varphiderivative_1} we have that
$$
\Big| \varphi(m+\eps\delta_k,\psi)-\Big(\eps\log\frac{m_k}{q_k}+\eps t_k+\varphi_{\infty,\infty,\eps}(m,\psi)\Big)\Big|\le O(\eps^2),\qquad\eps> 0.
$$
This is the announced approximate variational formula for $\varphi(m+\eps\delta_k,\psi)$.
\end{proof}

\begin{step}For any $k\in\N$, the right-derivative of $m_k\mapsto \varphi(m,\psi)$ satisfies \lq$\leq$\rq\ in \eqref{DerivativeMicro_Eq}.
\end{step}

\begin{proof}
As in the preceding step, we fix $k\in\N$ and $m\in[0,\infty)^\N$ satisfying $m_k>0$. Observe that $\varphi(m,\psi)$ is by definition equal to $\varphi_{\infty,\infty,0}(m,\psi)$ defined in \eqref{varpiAKdeltadef}. Hence, we obtain for the directional right-derivative
\begin{equation}\label{derivvarphiuppbound}
\begin{aligned}
\overline\partial_{m_k}^+\varphi(m,\psi)&:=\limsup_{\eps\downarrow 0}\eps^{-1}\big(\varphi(m+\eps \delta_k,\psi)-\varphi(m,\psi)\big)\\
&\leq \log\frac{m_k}{q_k}+ t_k-\log\sum_{a\in \N_0} \psi(a)\e^{-2\v ak}P_a(\e^{-2\Phi^{\ssup k}}), 
\end{aligned}
\end{equation}
by plugging in any minimiser $(P_a)_a$ of the formula \eqref{varphidef} of $\varphi(m,\psi)$ into the infimum on the right-hand side of \eqref{varphiderivative}. Since this upper bound holds for any such minimiser, we arrive at the claimed upper bound in \eqref{DerivativeMicro_Eq} for the right-derivative.
\end{proof}

\begin{step}For any $k\in\N$, the right-derivative  of $m_k\mapsto \varphi(m,\psi)$ satisfies also the opposite inequality in \eqref{DerivativeMicro_Eq}.
\end{step}

\begin{proof}We pick
a minimiser $(P^{\ssup{\eps}}_b)_{b\in\N_0}$ for the infimum on the right-hand side of  \eqref{varphiderivative}. Then, as $\eps\downarrow 0$, each family $(P_b^{\ssup{\eps}})_{\eps>0}$ with $b\in\N_0$ possesses at least one accumulation point $\widetilde P_b$, by non-negativity of each of the minimised terms and compactness of the level sets of $I$. Using the bounded-convergence theorem $\sum_a\psi(a)P_a^{\ssup{\eps}}$ also converges to $\sum_a\psi(a)\widetilde P_a$ in the local tame topology. Since the map $P\mapsto P(N_0^{\ssup{\delta_l}})$ is continuous, $(\widetilde P_b)_{b\in\N_0}$ is admissable in the variational formula for $\varphi(m,\psi)$. Again by lower semicontinuity, $(\widetilde P_b)_{b\in\N_0}$ is a minimiser for that formula, and we obtain the following lower bound for the right-derivative:
\begin{equation}\label{derivvarphilowbound}
\begin{aligned}
\underline\partial_{m_k}^+\varphi(m,\psi)
&:=\liminf_{\eps\downarrow 0}\eps^{-1}\big(\varphi(m+\eps\delta_k,\psi)-\varphi(m,\psi)\big)\\
&\geq \log\frac{m_k}{q_k}+t_k+
\liminf_{\eps\downarrow0}\big[-\log \sum_{b\in \N_0} \psi(b)\e^{-2\v b k}P^{\ssup{\eps}}_b(\e^{-2\Phi^{\ssup k}})\big]\\
&\geq \log\frac{m_k}{q_k}+t_k -\log \sum_{b\in \N_0} \psi(b)\e^{-2\v bk}\widetilde P_b(\e^{-2\Phi^{\ssup k}})\\
&\geq \log\frac{m_k}{q_k}+t_k -\sup_{(P_a)_a}\log \sum_{b\in \N_0} \psi(b)\e^{-2\v bk}P_b(\e^{-2\Phi^{\ssup k}}),
\end{aligned}
\end{equation}
where the supremum is on all minimisers $(P_a)_a$ in the formula for $\varphi(m,\psi)$. In the third line, we used first Fatou's lemma and then the local tame convergence with an additional spatial-truncation for $\e^{-2\Phi^{\ssup k}}$, see our argument around~\eqref{continuitylower} below for details. Now we see that both right-hand sides of \eqref{derivvarphilowbound} and \eqref{derivvarphiuppbound} coincide, and we have proved the lower bound  in  \eqref{derivvarphilowbound} for the right-derivative instead of the derivative.
\end{proof}

So far, we have proved that the right-derivative  of $m_k\mapsto \varphi(m,\psi)$ exists and is given by the right-hand side of \eqref{derivvarphilowbound}.

\begin{step} $\varphi(\cdot,\psi)$ is differentiable, and \eqref{derivvarphilowbound} holds.
\end{step}

\begin{proof}
By convexity, see Lemma~\ref{lem-properties}(1), it suffices to show that $\partial_{m_k}^+\varphi(m,\psi)$ is left-continuous in $m_k$. Take a sequence $(\eps_n)_{n\in\N}$ in $(0,1)$ that converges to zero as $n\to\infty$, and take a sequence $(P^{\ssup{\eps_n}}_b)_{b\in\N_0}$, $n\in\N$, of minimisers for the  formula \eqref{varphidef} of $\varphi(m-\eps_n\delta_k,\psi)$ that asymptotically optimises the term $\sum_{b\in \N_0} \psi(b)\e^{-2\v bk}P_b(\e^{-2\Phi^{\ssup k}})$. Again by the compactness of the level sets of $I$ and by lower-semicontinuity of $P\mapsto I(P)+P(\Phi)$ and continuity of $P\mapsto P(N_0^{\ssup{\delta_l}})$ for any $l\in\N$, we see that, along some subsequence, $(P^{\ssup{\eps_n}}_b)_{b\in\N_0}$ converges as $n\to\infty$ towards some $(\widetilde P_b)_{b\in\N_0}$, and the latter is minimal in the formula for $\varphi(m,\psi)$. Now we see that
\begin{equation}\label{continuitylower}
\liminf_{n\to\infty}\big[-\log P^{\ssup{\eps_n}}_b(\e^{-2\Phi^{\ssup k}})\big]\geq -\log \widetilde P_b(\e^{-2\Phi^{\ssup k}}),
\end{equation}
Indeed, for any centred box $Q$, we introduce $\Phi^{\ssup k}_{Q}(\omega)$ as the interaction between $\delta_{(0,G_k)}$ with all the particles that belong to points in $Q$. It is a local bounded functional, and we have $\Phi^{\ssup k}_{Q}(\omega)\leq \Phi^{\ssup k}(\omega)$, hence:
$$\widetilde P_b (\e^{-2\Phi^{\ssup k}_{Q}})= \lim_{n\to \infty} P^{\ssup{\eps_n}}_b(\e^{-2\Phi^{\ssup k}_{Q}})\geq \limsup_{n\to \infty} P^{\ssup{\eps_n}}_b(\e^{-2\Phi^{\ssup k}}).$$
Letting $Q\to \Z^d$, we get
$$\widetilde P_b (\e^{-2\Phi^{\ssup k}})\geq \limsup_{n\to \infty} P^{\ssup{\eps_n}}_b(\e^{-2\Phi^{\ssup k}}),
$$
implying \eqref{continuitylower}.
Further, using first \eqref{derivvarphiuppbound} and then \eqref{continuitylower} together with the asymptotic optimality, we have that
$$ 
\begin{aligned}
\partial_{m_k}^+\varphi(m,\psi) &\leq \log\frac{m_k}{q_k}+t_k -\log \sum_{b\in \N_0} \psi(b)\e^{-2\v bk}\widetilde P_b(\e^{-2\Phi^{\ssup k}}) \leq \liminf_{n\to \infty} \partial_{m_k}^+\varphi(m-\eps_n \delta_k,\psi)\\
&\leq \limsup_{n\to \infty} \partial_{m_k}^+\varphi(m-\eps_n \delta_k,\psi)\leq \partial_{m_k}^+\varphi(m,\psi),
\end{aligned}
$$
where the last inequality comes from the monotonicity of $m_k\mapsto \partial_{m_k}^+\widetilde\varphi(m,\psi)$. This concludes the proof of the left-continuity of $m_k\mapsto \partial_{m_k}^+\widetilde\varphi(m,\psi)$ as desired.
\end{proof}

This ends the proof of Lemma \ref{lem-DerivativeMicro}. Let us draw a corollary from the preceding proof:

\begin{lemma}
For any sequence $(\eps_K)_{K\in\N}$ in $(0,\infty)$ tending to zero, the map $P\mapsto P(\e^{-2\Phi^{\ssup k}})$ is continuous in the local tame topology on the set of all $P\in\Mcal_1^{\ssup{\rm s}}(\Omega)$ such that $\sum_{l\geq K } lP(N_0^{\ssup\ell})\leq \eps_K$ for any $K\in\N$.
\end{lemma}    
   
\begin{proof}
The upper semicontinuity was shown below \eqref{continuitylower}. In order to prove the lower semicontinuity, we obtain a lower bound for $P(\e^{-2\Phi^{\ssup k}})$ by inserting the indicator on the event $\cap_{l\in\N}\{M_{Q^c,V_k}^{\ssup {\delta_l}}=0\}$ that no particle attached to a point outside the box $Q$ lies in $V_k=G_k+\supp(v)$. On this event, we can replace $\Phi^{\ssup k}$ by $\Phi^{\ssup k}_{Q}(\omega)$, the interaction between $\delta_{(0,G_k)}$ with all the particles that belong to points in $Q$. This gives
$$
P(\e^{-2\Phi^{\ssup k}})\geq P\big(\e^{-2\Phi_Q^{\ssup k}}\prod_{l\in\N} \1\{M_{Q^c,V_k}^{\ssup {\delta_l}}=0\}\big)
\geq P(\e^{-2\Phi^{\ssup k}_{Q}})-P\big(\bigcup_{l\in\N}\{M_{Q^c,V_k}^{\ssup {\delta_l}}\not=0\}\big).
$$
Now we estimate
$$
\begin{aligned}
P\big(\bigcup_{l\in\N}\{M_{Q^c,V_k}^{\ssup {\delta_l}}\not=0\}\big)
&\leq\sum_{z\in V_k}\sum_{l\in\N}\sum_{x\in Q^{\rm c}\colon z\in x+G_l} P(N_x^{\ssup{\delta_l}})\not=0)\\
&\leq \sum_{z\in V_k}\sum_{l\in\N}\sum_{x\in Q^{\rm c}\cap(z-G_l)}P(N_x^{\ssup{\delta_l}})\\
&\leq \sum_{z\in V_k}\sum_{l\in\N}|Q^{\rm c}\cap(z-G_l)| P(N_0^{\ssup{\delta_l}})\\
&\leq \sum_{z\in V_k}\sum_{l>K_Q} lP(N_0^{\ssup{\delta_l}}) \leq |V_k|\, \eps_{K_Q},
\end{aligned}
$$
where we picked a large $K_Q\in\N$ such that $Q^{\rm c}\cap(z-G_l)$ is empty for any $z\in V_k$ and any $l\leq K_Q$. Since we can choose $K_Q$ such that $K_Q\to\infty$ as $Q\uparrow \Z^d$, and since $\eps_K\to0 $ as $K\to\infty$, we see that $P\big(\bigcup_{l\in\N}\{M_{Q^c,V_k}^{\ssup {\delta_l}}\not=0\}\big)$ vanishes uniformly in these $P$ as $Q\uparrow \Z^d$. This makes it easy to finish the proof.
\end{proof}

\subsection{Differentiability of $\chi$: proofs of Lemma \ref{lem-ELeq} and Corollary \ref{cor-chidiff}}\label{sec-EL-eqs}

In this section, we prove Lemma \ref{lem-ELeq} and Corollary \ref{cor-chidiff}, i.e., the Euler--Lagrange analysis of the minimiser $m$ of $\varphi(m,\psi)$ defined in \eqref{varphidef}, and the resulting differentiability of $\chi(\cdot)$.

\begin{lemma}[Positivity of minimising $m$]\label{lem-mkpos}
Fix $\rho\in(0,\infty)$ and $\psi\in\Mcal_1(\N_0)$ with $\rho_{\rm ma}=\sum_{a\in\N_0} a\psi(a)\in[0,\rho)$ and assume that $m=(m_k)_{k\in\N}$ is a minimiser of $\varphi(\cdot,\psi)$ under the assumption that $\sum_{k\in\N} k m_k=\rho- \rho_{\rm ma}= \rho_{\rm mi}$. Then $m_k>0$ for any $k\in\N$.
\end{lemma}

\begin{proof}
This is a well-known argument that is based on the fact that the slope of $x\mapsto x\log x$ at zero is equal to $-\infty$. Recall from lemma \ref{lem-properties} that $I(P)=H(m|q)+I_m(P)$ for any $P\in\Mcal_1^{\ssup{\rm s}}(\Omega)$ satisfying $P(N_0^{\ssup{\delta_k}})=m_k$ for any $k\in\N$, where we wrote $I_m$ for the entropy density function $I$ defined in \eqref{Igriddef} with $q$ replaced by $m$, where we recall that
\begin{equation}
H(m|q)=\sum_{k\in\N}\Big(q_k-m_k+m_k\log\frac{m_k}{q_k}\Big)
\end{equation}
is the relative entropy of $m$ with respect to $q$.

As usual, the convention $0\log 0=0$ is in force and makes $m\mapsto H(m|q)$ continuous coordinate-wise in $[0,\infty)^\N$. If now $m$ is a minimiser of $\varphi(\cdot,\psi)$ under the constraint $\sum_k k m_k= \rho_{\rm mi}$ and has a zero $m_k=0$, then one can construct $\widetilde m$ from $m$ by putting some small positive mass $\delta$ at $k$ (i.e., $\delta=\widetilde m_k$) and subtracting at some $\widetilde k$ with $m_{\widetilde k}>0$ some mass in such a way that $\sum_l l \widetilde m_l= \rho_{\rm mi}$. The resulting entropy difference is $H(m|q)-H(\widetilde m|q)=\delta(C-\log  \delta)$ for some $C$, depending on $m_{\widetilde k}$ and $q_{\widetilde k}$ and $q_k$. For $\delta$ sufficiently small, this is positive. Since $\partial_{m_k}[\varphi(m,\psi)-H(m|q)]$ is finite, as we have seen in Lemma~\ref{lem-DerivativeMicro}, we see that $\varphi(\widetilde m,\psi)<\varphi(m,\psi)$ for sufficiently small $\delta>0$, in contradiction to the minimality of $m$.  
\end{proof}

We employ the Euler--Lagrange formalism only for perturbations in the direction $m$ and keep $\psi$ fixed. We thus solve the following system of equations
\begin{equation}\label{EulerLangrangeEqn}
\partial_{m_k}\Big[\varphi(m,\psi)-\alpha\sum_{l\in\N}l m_l \Big]=0,\qquad k\in\N,
\end{equation}
where $\alpha\in\R$ is the Lagrange multiplier, to be adjusted such that the constraint $\sum_{k\ge 1}km_k= \rho_{\rm mi}$ is satisfied. Clearly, the conclusion is then that 
\begin{equation}\label{ELequationsvarphi}
\alpha k =\partial_{m_k}\varphi(m,\psi),\qquad k\in\N.
\end{equation}
By Lemma \ref{lem-DerivativeMicro}, this finishes the proof of Lemma \ref{lem-ELeq}.

Now we can also give the proof of Corollary \ref{cor-chidiff}:

Fix $\rho\in(0,\infty)$. We show that $\chi(\cdot)$ is differentiable with $\chi'(\rho)=\alpha$, the Euler--Lagrange multiplier of Lemma \ref{lem-ELeq}.

According to Lemma \ref{lem-ProofExistMin}, we can take a minimiser $(m,\psi)$ of $\varphi$ with density $\rho$. Then, using Lemma~\ref{lem-mkpos}, $m_k>0$ for any $k\in\N$ and using Lemma \ref{lem-DerivativeMicro},  $m_k\mapsto \varphi(m,\psi)$ is differentiable. According to the Euler--Lagrange equations in \eqref{ELequationsvarphi}, there is some $\alpha$ such that $\partial_{m_k}\varphi(m,\psi)=\alpha k$ for any $k\in \N$.   As $\chi(\cdot)$ is convex, it has left- and right-derivatives at $\rho$. Note that, for $\eps\to0$,
$$\chi(\rho+\eps k)\leq \varphi(m+\eps \delta_k, \psi)=\varphi(m, \psi) + \alpha k \eps +o(\eps)=\chi(\rho)+ \alpha k \eps +o(\eps).
$$
Using this first for $\eps\downarrow0$ and then for $\eps\uparrow 0$, we get
$$
\chi'^+(\rho)\leq \alpha\qquad\mbox{and}\qquad \chi'^-(\rho)\geq \alpha.
$$
By convexity,  $\chi'^-(\rho)\leq \chi'^+(\rho)$, hence we get that $\chi(\cdot)$ is differentiable at $\rho$ with derivative equal to $\alpha=\frac 1k \partial_{m_k} \varphi(m,\psi)$.

\subsection{An alternate formula for the derivative}\label{sec-thinning}
Let us give a brief heuristic derivation of the formula in \eqref{DerivativeMicro_EqAlter} for the derivative of $m\mapsto \varphi(m,\psi)$. We will do this only for $\psi=\delta_0$. We are sure that a full proof can be given using Lemma \ref{lem-varphirepr} and arguments similar to those that we carried out in the proof of Lemma \ref{lem-DerivativeMicro} in Section \ref{sec-EL-eqs}.

We start from the formula for $\varphi(m, \delta_0)$ as the negative exponential rate of 
\begin{equation}\label{Zthinning}
\begin{aligned}
Z_{N,\L_N}(m)
& =\E\Big[\e^{-\Phi_{\Lambda_N,\Lambda_N}(\omega_{\rm P})}\prod_{k\in\N}\1\{N_{\L_N}^{\ssup{\delta_k}}(\omega_{\rm P})=m_k |\L_N|\}\Big]\\
& = \e^{-H(m|q)|\L_N|} \bigotimes_{k\in \N} \Ucal_{\L_N}^{\otimes m_k\L_N}\Big[\e^{-\Phi_{\Lambda_N,\Lambda_N}(\omega_{\rm U})}\Big].
\end{aligned}
\end{equation}
(For simplicity, we do not introduce auxiliary parameters $A,K$ and $\delta$ for this heuristics.)

Now take $\widetilde{m}, m \in [0,\infty)^\N$ such that $\widetilde{m}_k\leq m_k$ for all $k\in \N$. Let $K$ be the set of indices $k$ such that $\widetilde m_k < m_k$. Let us first state a general formula that writes $Z_{N,\L_N}(\widetilde m)$ as a thinning of $Z_{N,\L_N}(m)$. We will see the process with $k$-densities $\widetilde m_k$ as a thinning of the process with $k$-densities $m_k$. Let us write the expectation with respect to the reference PPP by first sampling the Poisson number of points with marks $G_k$ for each $k$ and then sampling the locations with the uniform distribution $\Ucal_{\L_N}$ on $\L_N$.  We  use  the symbol $\widetilde\omega_{\rm U}=\sum_{k\in\N}\sum_{i=1}^{\widetilde m_k|\Lambda_N|}\delta_{(X_{k,i},G_k)}$ for the arising marked random point process under the measure $\bigotimes_{k\in \N} \Ucal_{\L_N}^{\otimes \widetilde m_k|\L_N|}$, and we use the analogous notation if $\widetilde m$ is replaced by $m$. We define for any $k\in K$ a random uniform subset $B_k$ of  the index set $[m_k |\L_N|]$ with size $(m_k-\widetilde m_k)|\L_N|$. The law of $B_k$ is denoted by $\mathcal{B}_k$, and the point process that is selected by $(B_k)_{k\in K}$ is denoted by $\omega_{\rm B}$. Then the distribution of $\widetilde \omega_{\rm U}$ under $\bigotimes_{k\in \N} \Ucal_{\L_N}^{\otimes \widetilde m_k|\L_N|}$ is equal to the distribution of $\omega_{\rm U}-\omega_{B}$ under $\bigotimes_{k\in \N}\Ucal_{\L_N}^{\otimes m_k|\L_N|}\otimes \bigotimes_{k\in K} \Bcal_k$.
Hence, asymptotically as $N\to\infty$, we obtain
\begin{equation}\label{eqthinning}
\begin{aligned}
Z_{N,\L_N}(\widetilde{m})
&=  \e^{-H(\widetilde{m}|q)|\L_N|} \bigotimes_{k\in \N} \Ucal_{\L_N}^{\otimes \widetilde m_k|\L_N|}\Big[\e^{-\Phi_{\Lambda_N,\Lambda_N}(\widetilde\omega_{\rm U})}\Big] \\
&= \e^{-H(\widetilde{m}|q)|\L_N|} \bigotimes_{k\in \N}\Ucal_{\L_N}^{\otimes m_k|\L_N|}\otimes \bigotimes_{k\in K} \Bcal_k\Big[\e^{-\Phi_{\Lambda_N,\Lambda_N}(\omega_{\rm U}-\omega_{\rm B})}\Big] \\
&=   \e^{-(H(\widetilde m|q)-H(m|q))|\L_N|}\E \Big[\bigotimes_{k\in K} \mathcal{B}_k\big[\e^{-\Phi_{\Lambda_N,\Lambda_N}(\omega_{\rm P}-\omega_{\rm B})}\big]\prod_{k\in\N}\1\{N_{\L_N}^{\ssup{\delta_k}}(\omega_{\rm P})=m_k|\L_N|\} \Big].
\end{aligned}
\end{equation}

Now we compute the left-derivative of $\varphi(m,\delta_0)$ using formula \eqref{eqthinning} for a special choice of $\widetilde m$.
We fix $k\in \N$. We assume $m_k>0$, take a small $\varepsilon$ and put $\widetilde m=m(\eps)=m-\eps \delta_k$. Formula \eqref{eqthinning} gives:
\begin{equation*}
    \begin{aligned}
    Z_{N,\L_N}(m(\eps))
    &=\e^{-(H(m(\eps)|q)-H(m|q))|\L_N|}\\
    & \times \E \Big[\mathcal{B}_k^{\ssup\eps}\big[\e^{-\Phi_{\Lambda_N,\Lambda_N}(\omega_{\rm P}-\omega_{\rm B})}\big]\prod_{l\in\N}\1\{N_{\L_N}^{\ssup{\delta_l}}(\omega_{\rm P})=m_l|\L_N|\} \Big],
    \end{aligned}
\end{equation*}
where $\mathcal B_k^{\ssup{\eps}} $ is the uniform law over all subsets of $[m_k |\L_N|]$ of size $\eps  |\L_N|$. Explicitly, $\omega_{\rm B}=\sum_{i=1}^{\eps  |\L_N|} \delta_{(X_{k,U_i},G_k)}$, where $U_1,\dots,U_{\eps  |\L_N|}$ are picked according to $\mathcal{U}_{[m_k|\L_N|]}^{\otimes \eps |\L_N|}$, conditioned on the event $\{\forall i\neq j, U_i\neq U_j\}$. Note that the probability of the latter has an exponential rate that is $o(\epsilon)$ as $\epsilon$ goes to 0. Hence we will assume that we can remove this conditioning, with a cost $e^{|\L_N|o(\epsilon)}$: for the upper bound, this is immediate; but the lower bound requires some fine work, we do not elaborate on that.

Observe that 
$$
\Phi_{\Lambda_N,\Lambda_N}(\omega_{\rm P}-\omega_{\rm B})=\Phi_{\Lambda_N,\Lambda_N}(\omega_{\rm P})-2\Phi^{\ssup{\leftrightarrow}}_{\Lambda_N,\Lambda_N}(\omega_{\rm P},\omega_{\rm B})+\Phi_{\Lambda_N,\Lambda_N}(\omega_{\rm B}).
$$ 
Furthermore, we also assume in this heuristics that all the boxes $U_i+G_k$ for $i\in[\eps |\L_N|]$ do not overlap each other (we demonstrated in Section \ref{sec-EL-eqs} how to control the complement of this event). Then  $\Phi_{\Lambda_N,\Lambda_N}(\omega_{\rm B})= \eps  t_k |\L_N|$. Therefore, we can proceed with
\begin{equation}
    \begin{aligned}
    Z_{N,\L_N}(m(\eps ))
        &=  \e^{-(H(m(\eps )|q)-H(m|q))|\L_N|}\e^{-\eps  t_k |\L_N|}\e^{|\L_N| o(\eps)}\\
        & \times \E \Big[\e^{-\Phi_{\Lambda_N,\Lambda_N}(\omega_{\rm P})}  U_{[m_k|\L_N|]}^{\otimes \eps  |\L_N|} \big[\e^{2\Phi^{\ssup{\leftrightarrow}}_{\Lambda_N,\Lambda_N}(\omega_{\rm P},\omega_{\rm B})}\big]\prod_{l\in\N}\1\{N_{\L_N}^{\ssup{\delta_l}}(\omega_{\rm P})=m_l|\L_N|\} \Big].
    \end{aligned}
\end{equation}

Now we introduce the empirical individual field 
$$
{\mathcal{R}}^\circ_{\L_N}(\omega)=\frac{1}{m_k|\L_N|}\sum_{i=1}^{m_k |\L_N|} \delta_{\theta_{X_{k,i}}(\omega)}=\frac{1}{m_k|\L_N|}\sum_{x\in \L_N} \xi^\ssup{k}(x)\delta_{\theta_x(\omega)}=\frac{1}{|\L_N|}\sum_{x\in \L_N} \frac{N_0^\ssup{\delta_k}(\theta_x(\omega))}{m_k} \delta_{\theta_x(\omega)},
$$
which is the Palm version of the empirical stationary field $\Rcal_{\Lambda_N}$ that we introduced in the proof of 
Theorem \ref{thm-freeenergygrid}. Then we have that
\begin{align*}
\mathcal{U}_{[m_k|\L_N|]}^{\otimes \eps  |\L_N|}\big[\e^{2 \Phi^{\ssup{\leftrightarrow}}_{\Lambda_N,\Lambda_N}(\omega_{\rm P},\omega_{\rm B})}\big]
& = \Big( \mathcal{U}_{[m_k|\L_N|]}\big[\e^{2 \Phi^{\ssup{\leftrightarrow}}_{\Lambda_N,\Lambda_N}(\omega_{\rm P},\delta_{(X_{k,U},G_k)})}\big]\Big)^{\eps  |\L_N|}\\
& = \Big( \big\langle R_{\L_N}^{\circ}(\omega_{\rm P}), \e^{ 2\Phi^{\ssup k}_{\Lambda_N}}\big\rangle \Big)^{\eps  |\L_N|}= \Big( \Big\langle R_{\L_N}(\omega_{\rm P}), \frac{N_0^\ssup{\delta_k}}{m_k}\e^{ 2\Phi^{\ssup k}_{\Lambda_N}}\Big\rangle \Big)^{\eps  |\L_N|}.
\end{align*}

Proceeding as in Section \ref{sec-varphidiffble}, we get:
$$
\begin{aligned}
    \varphi(&m(\eps ),\delta_0)  =
    -\eps  \log \frac{m_k}{q_k}+\eps  t_k+o(\eps )\\
    & + \inf \Big\{I(P)+P(\Phi_{0,\Z^d})-\eps  \log \frac{P(N_0^\ssup{\delta_k}\e^{2\Phi^\ssup{k}})}{m_k}\colon P\in \mathcal{M}_1^\ssup{{\rm s}}(\Omega), P(N_0^\ssup{\delta_k})=m_k\ \forall k \in\N\Big\}.
\end{aligned}
$$
Likewise, we obtain:
$$\partial_{m_k}^- \varphi(m, \delta_0)=\log\frac{m_k}{q_k}-t_k+ \sup_P \log \frac{P(N_0^\ssup{\delta_k}\e^{2\Phi^\ssup{k}})}{m_k},$$
where the maximum is taken over all minimisers $P$ in the definition of $\varphi(m, \delta_0)$. This ends our heuristic derivation of \eqref{DerivativeMicro_EqAlter} for $\psi=\delta_0$.

\subsection{Qualitative description in case of a phase transition: proof of Lemma \ref{lem-saturation}}\label{sec-naturephasetrans}
    
In this section, we prove Lemma \ref{lem-saturation}. Recall that we assume that $\rho_{\rm c}$ is finite. Also recall that \eqref{A(psi)def} is an alternative representation of $\chi$ as defined in~\eqref{chi(m,rho2)def}, and that the minimisers $\psi$ coincide. 

\setcounter{step}{0}

\begin{step}\label{step-varformident}
We have the following alternative representation of  the free energy
\begin{equation}\label{infimumforchi}
\begin{aligned}
\chi(\rho)=\inf\Big\{\sum_{a\in\N_0}& \psi(a)[\chi(\rho_a)+2\v a \rho_a+\v a^2] \colon \\
&\psi\in \Mcal_1(\N_0), (\rho_a)_{a\in\N_0}\in [0,\rho_{\rm c}]^{\N_0}, \sum_{a\in \N_0} \psi(a) [\rho_a+a]=\rho\Big\}.
\end{aligned}
\end{equation}
In other words, the epigraph of $\chi(\cdot)$ is the convex hull of the epigraphs of the functions
$$
f_a\colon [a, \rho_{\rm c}+a]\to [0,\infty),\qquad \rho+a\mapsto \chi(\rho)+2{\v} a \rho + {\v} a^2,\qquad a\in \N_0, \rho\in[0,\infty).
$$
Furthermore, we can restrict in \eqref{infimumforchi} to $(\rho_a)_{a\in\N_0}\in [0,\rho_{\rm c}]\times[\rho_{\rm c}-1,\rho_{\rm c}]^{\N}$.
\end{step}

\begin{proof}
Based on the considerations made in the proof of Step~\ref{step-twoatoms} in Section~\ref{sec-minimivarphi}, we see that formula~\eqref{infimumforchi} coincides with formula~\eqref{A(psi)def}, when we let $(\rho_a)_{a\in\N_0}\in [0,\infty)^{\N_0}$. Further, we see from the argument presented around \eqref{descrminimisereq} that $\chi(\rho_a)$ admits a minimiser $(m^a,\delta_0)$, and hence, $\rho_a\le \rho_{\rm c}$ for any $a\in \N$ by definition of $\rho_{\rm c}$. Hence the formula \eqref{infimumforchi} follows. We rewrite it as
\begin{equation}
\begin{aligned}\label{infimumforchi2}
\chi(\rho)=\inf\Big\{\sum_{a\in\N_0}& \psi(a)f_a(\widetilde \rho_a) \colon \psi\in \Mcal_1(\N_0), (\widetilde \rho_a)_{a\in\N_0}\in [0,\infty)^{\N_0},\\
& \forall a\in\N, \widetilde\rho_a \in[a, \rho_c+a], \sum_{a\in \N_0} \psi(a) \widetilde\rho_a=\rho\Big\}.
\end{aligned}
\end{equation}
Now, observe that Step \ref{step-chiinequality} in Section \ref{sec-minimivarphi} implies that $a\leq b$ implies $f_a(\rho)\leq f_b(\rho)$ whenever this is defined. Therefore, using convexity of the $f_a$'s, we can restrict in \eqref{infimumforchi2} to $\widetilde \rho_0\in[0, \rho_c]$ and $\widetilde\rho_a\in[\rho_c+a-1, \rho_c+a]$ when $a\in\N$. Correspondingly, in \eqref{infimumforchi} we can restrict to $(\rho_a)_{a\in\N_0}\in [0,\rho_{\rm c}]\times[\rho_{\rm c}-1,\rho_{\rm c}]^{\N}$.
\end{proof}

\begin{step}The saturation hypothesis~\eqref{saturation} is false.
\end{step}

\begin{proof}
Observe that, if the saturation hypothesis~\eqref{saturation} were true, then for any minimiser in the \eqref{infimumforchi} at density $\rho\geq \rho_{\rm c}$, we would have $\rho_a=\rho_{\rm c}$ for any $a$, since $\sum_{a\in\N_0}\psi(a)\rho_a=\rho-\sum_{a\in\N_0}\psi(a)a=\rho-\rho_{\rm ma}=\rho-[\rho-\rho_{\rm c}]_+=\rho_{\rm c}$ and $\rho_a\le \rho_{\rm c}$. So from \eqref{infimumforchi} we have
$$\chi(\rho)=\chi(\rho_{\rm c})+2\v (\rho-\rho_{\rm c})\rho_{\rm c}+\v \inf\Big\{\sum_{a\in\N_0} \psi(a) a^2 \colon \psi\in \Mcal_1(\N_0), \sum_{a\in \N_0} \psi(a) a=\rho-\rho_{\rm c}\Big\}.$$
Evaluating the infimum explicitly, one sees that $\chi$ is a non-trivial polygon line on $[\rho_{\rm c},\infty)$ (it is equal to $\chi^{\rm sat}$ defined in next step) and thus not differentiable. This contradicts Corollary~\ref{cor-chidiff}.
\end{proof}

\begin{step}\label{non_sat_1}
$\chi'(\rho_{\rm c})>(2\rho_{\rm c}-1)\v$.
\end{step}

\begin{proof}
Let us define $\chi^{\rm sat}\colon \R^+\to \R^+$ coinciding with $\chi(\cdot)$ on $[0,\rho_{\rm c}]$, and equal on $[\rho_{\rm c},\infty)$ to the linear polygon line interpolating the values $\chi^{\rm sat}(\rho_{\rm c}+a)=f_a(\rho_{\rm c}+a)=\chi(\rho_{\rm c})+2\v a \rho_{\rm c}+\v a^2$, for $a\in\N_0$. The graph of $\chi^{\rm sat}$ is depicted in Figure~\ref{Pix_chidiff}.
\begin{figure}[!htpb]
	\centering
\includegraphics[scale=0.34]{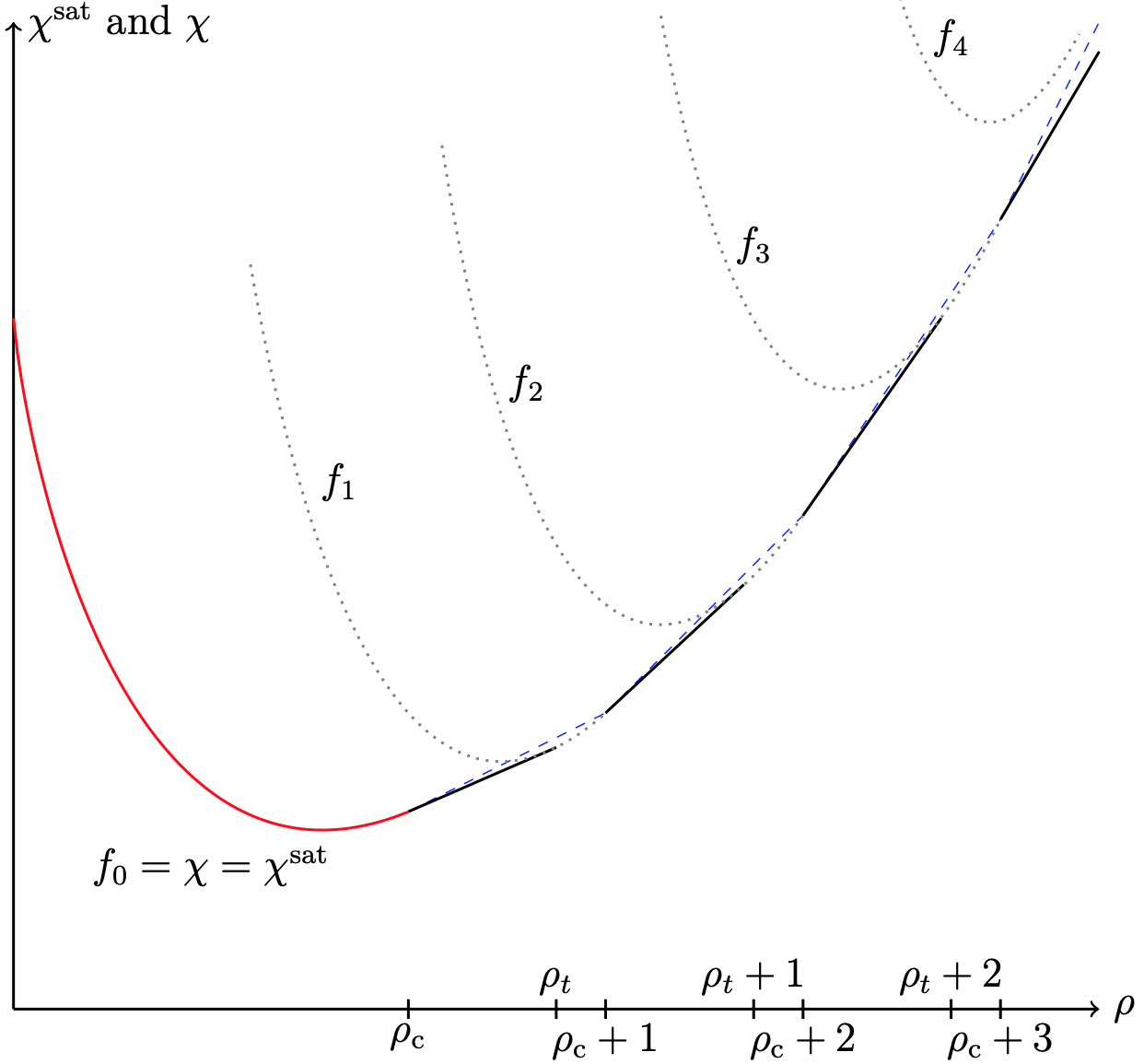}
	\caption{Illustrations of $\chi^{\rm {sat}}$ and $\chi$. The graph of $\chi^{\rm {sat}}$ is given by the red continuous line in $[0,\rho_{\rm c}]$ and continued by the blue dashed line segments to the right of $\rho_{\rm c}$. The grey dotted lines are the graphs of the functions $f_a\colon [a, \rho_{\rm c}+a]\to [0,\infty), f_a(\rho+a)= \chi(\rho)+2{\v} a \rho + {\v} a^2$ for $a\in \{1,\dots,4\}$. The four blue dashed lines are the segments joining the points $(\rho_{\rm c}+a-1, f_{a-1}(\rho_{\rm c}+a-1))$, for $a \in \{1,\dots, 5\}$. Observe that $\chi^{\rm {sat}}$ is not differentiable in $\rho_{\rm c}+\N$.
	The graph of $\chi$ coincides with the graph of $\chi^{\rm {sat}}$ on $[0,\rho_{\rm c}]$. Then it is continued by the alternating solid black line segments and grey dotted segments. The left-most black line is the tangent both, to $\chi$ in $(\rho_{\rm c}, \chi(\rho_{\rm c}))$ and to $f_1$ in $\rho_{\rm t}\in[\rho_{\rm c}, \rho_{\rm c}+1)$. The graph of $\chi$ coincides with this tangent line on $[\rho_{\rm c}, \rho_{\rm t}]$, with $f_1$ on $[\rho_{\rm t}, \rho_{\rm c}+1]$. Analogous assertions hold on $[\rho_{\rm c}+a, \rho_{\rm t}+a]$ and $[\rho_{\rm t}+a, \rho_{\rm c}+a+1]$ for $a\in\{2,3,4,5\}$.}
	\label{Pix_chidiff}
\end{figure}
Observe that $\chi^{\rm sat}$ is obtained by choosing $\psi=\delta_a$ and $\rho_a=\rho_{\rm c}$ in~\eqref{infimumforchi}, so by convexity of $\chi(\cdot)$, we have $\chi^{\rm sat}(\rho)\geq\chi(\rho)$. Also note that $\chi^{\rm sat}$ is convex, since $\N_0\ni a\mapsto \chi^{\rm sat}(\rho_{\rm c}+a)$ is convex, $\chi^{\rm sat}(\rho)=\chi(\rho)$ on $[0,\rho_{\rm c}]$ and $\chi^{\rm sat}(\rho)\geq\chi(\rho)$.

Now, assume that $\chi'(\rho_{\rm c})\leq (2\rho_{\rm c}-1)\v$. Our goal is then to show that  $\chi^{\rm sat}=\chi$. For this, note that for any $\rho\in [\rho_{\rm c}-1, \rho_{\rm c}]$, we have
\begin{align*}
    \chi^{\rm sat}(\rho+1)
    & =\chi(\rho_{\rm c})+(2\rho_{\rm c}+1)(\rho+1-\rho_{\rm c})\bar{v}\\
    & =\chi(\rho_{\rm c})+(2\rho_{\rm c}-1)(\rho-\rho_{\rm c})\bar{v}+(2\rho+1)\bar{v}\\
    & \leq \chi(\rho_{\rm c})+\chi'(\rho_{\rm c})(\rho-\rho_{\rm c})+(2\rho+1)\bar{v}\\
    & \leq \chi(\rho_{\rm c})+\chi'(\rho)(\rho-\rho_{\rm c})+(2\rho+1)\bar{v}\\
    & \leq \chi(\rho)+(2\rho+1)\bar{v},
\end{align*}
using the assumption in line three and convexity of $\chi(\cdot)$ in lines four and five. Using this, for any $a\in \N$, we derive
\begin{align*}
\chi^{\rm sat}(\rho+a)
& =\chi^{\rm sat}(\rho+1)+(2(a-1)\rho+(a^2-1))\bar{v}\\
& \leq  \chi(\rho)+(2\rho+1)\bar{v}+(2(a-1)\rho+(a^2-1))\bar{v}\\
& = \chi(\rho)+(2\rho+a)\bar{v}a.
\end{align*}
Now, using this and convexity of $\chi^{\rm sat}$, we have  for any $\rho\in[0,\infty]$ that
$$\chi^{\rm sat}(\rho)
\leq \sum_{a\in\N_0} \psi(a)\chi^{\rm sat}(\rho_a+a) 
\leq \sum_{a\in\N_0} \psi(a)[\chi(\rho_a)+2\v a \rho_a+\v a^2],$$
for any $\psi\in \Mcal_1(\N_0), (\rho_a)_{a\in\N_0}\in [0,\rho_{\rm c}]\times[\rho_{\rm c}-1,\rho_{\rm c}]^{\N}$ such that $\sum_{a\in \N_0} \psi(a) [\rho_a+a]=\rho$. But, taking the infimum over such $\psi$ and $(\rho_a)_{a\in\N_0}$, we see that $\chi^{\rm sat}(\rho)\leq \chi(\rho)$ for any $\rho$, by~\eqref{infimumforchi}.

As a consequence $\chi$ and $\chi^{\rm sat}$ coincide and thus $\chi$ is non-differentiable at any $\rho_{\rm c}+a$, for $a\in \N$, which contradicts Corollary~\ref{cor-chidiff}.
\end{proof}

So far, we have derived  that 
$$
(2\rho_{\rm c}-1)\v<\chi'(\rho_{\rm c})\leq (2\rho_{\rm c}+1)\v,
$$ 
where the second inequality comes from the convexity of $\chi$ and from $\chi(\rho_{\rm c}+1)\leq \chi(\rho_{\rm c})+(2\rho_{\rm c}+1)\v$, according to Step \ref{step-chiinequality} in Section \ref{sec-minimivarphi}.

\begin{step} There exists $\rho_{\rm t}\in[\rho_{\rm c}\lor 1, \rho_{\rm c}+1)$ such that
$$
\chi(\rho_{\rm t}-1)+(2\rho_{\rm t}-1)\v+(\rho_{\rm c}-\rho_{\rm t})(\chi'(\rho_{\rm t}-1)+2\v)=\chi(\rho_{\rm c})
$$
i.e. 
$$
f_1(\rho_{\rm t})+(\rho_{\rm c}-\rho_{\rm t})f_1'(\rho_{\rm t})=\chi(\rho_{\rm c}).
$$
\end{step}

\begin{proof}
We consider the function 
\begin{align*}
g(\rho)=\chi(\rho-1)+(2\rho-1)\v+(\rho_{\rm c}-\rho)(\chi'(\rho-1)+2\v),
\end{align*}
defined on $[1, \infty)$. If $\rho_{\rm c}\geq 1$, note that $g(\rho_{\rm c})=\chi(\rho_{\rm c}-1)+(2\rho_{\rm c}-1)\v\ge \chi(\rho_{\rm c})$ and $g(\rho_{\rm c}+1)=\chi(\rho_{\rm c})+(2\rho_{\rm c}-1)\v-\chi'(\rho_{\rm c})<\chi(\rho_{\rm c})$, where we used Step~\ref{step-chiinequality} in Section~\ref{sec-minimivarphi} and Step~\ref{non_sat_1} above.

Applying the intermediate value theorem on $[\rho_{\rm c}, \rho_{\rm c}+1]$ to $g$ (which is continuous since $\chi'$, the derivative of a differentiable convex function, is continuous), we obtain the existence of $\rho_{\rm t} \in [\rho_{\rm c}, \rho_{\rm c}+1)$ such that $g(\rho_{\rm t})=\chi(\rho_{\rm c})$. If $\rho_{\rm c}<1$, we may also apply the intermediate value theorem to $g$ on the interval $[1,\rho_{\rm c}+1]$, using that $\chi'(0)=-\infty$.
\end{proof}

Now, we define $\bar\chi$ as on the right-hand side of \eqref{chiidentnew}, i.e., for $\rho\in [0, \rho_{\rm c}]$, $\bar\chi(\rho)=f_0(\rho)=\chi(\rho)$, and, for any $a\in\N_0$ and $\rho\in[\rho_{\rm c}, \rho_{\rm c}+1]$,
\begin{equation*}
\bar\chi(\rho+a)=\begin{cases}
\frac{\rho_{\rm t}-\rho}{\rho_{\rm t}-\rho_{\rm c}}f_a(\rho_{\rm c}+a)+\frac{\rho-\rho_{\rm c}}{\rho_{\rm t}-\rho_{\rm c}}f_{a+1}(\rho_t+a+1)&\mbox{if }\rho\leq\rho_{\rm t},\\
f_{a+1}(\rho+a)&\mbox{if }\rho\geq  \rho_{\rm t}.
\end{cases}
\end{equation*}

Our aim is to show that $\chi=\bar{\chi}$. First observe that $\bar{\chi}(\rho)$ is obtained by some particular choice in~\eqref{infimumforchi}, so $\bar{\chi}(\rho)\geq\chi(\rho)$. Also note that for any $\rho\geq \rho_{\rm t}-1$, for any $a\in \N_0$, \begin{equation}\label{eqsimbarchi}
\bar\chi(\rho+a)=\bar\chi(\rho)+2\v a \rho +\v a^2.
\end{equation}

Figure \ref{Pix_chidiff} depicts the graph of $\bar\chi$, which will turn out in the end of the proof to be identical with $\chi$.

\begin{step}\label{stepbarchiconv}
The function $\bar\chi$ is convex.
\end{step}

\begin{proof}
$\bar\chi$ is obviously convex on the intervals $[\rho_{\rm c}+a,\rho_{\rm t}+a]$ and $[\rho_{\rm t}+a,\rho_{\rm c}+a+1]$ for $a\in\N_0$. The only difficulty is to show that $\partial^- \bar\chi\leq\partial^+ \bar\chi$ at $\rho_{\rm c}+a$ and $\rho_{\rm t}+a$, for all $a\in \N_0$. Using \eqref{eqsimbarchi}, it is enough to show this for $a=0$. 

First, we have, using convexity of $\chi(\cdot,0)$ and Step \ref{step-chiinequality} in Section \ref{sec-minimivarphi}:
$$
\begin{aligned}
    \partial^- \bar\chi(\rho_{\rm c},0)&=\partial^-_{ \rho_{\rm mi}} \chi(\rho_{\rm c},0)\leq \partial^+_{ \rho_{\rm mi}} \chi(\rho_{\rm c},0)\leq \frac{\chi(\rho_{\rm t})-\chi(\rho_{\rm c})}{\rho_{\rm t}-\rho_{\rm c}}\\
    &\leq \frac{\chi(\rho_{\rm t}-1)+(2\rho_{\rm t}-1)\v-\chi(\rho_{\rm c})}{\rho_{\rm t}-\rho_{\rm c}}=\partial^+ \bar\chi(\rho_{\rm c}).
\end{aligned}
$$

Secondly, by convexity of $\rho\mapsto \chi(\rho-1)+(2\rho-1)\v$, we have for all $\rho\in[\rho_{\rm t}, \rho_{\rm c}+1]$, 
$$
\begin{aligned}
\bar\chi(\rho_{\rm t})+\partial^-\bar\chi(\rho_{\rm t})(\rho-\rho_{\rm t})&=\chi(\rho_{\rm t}-1)+(2\rho_{\rm t}-1)\v+(\rho-\rho_{\rm t})(\chi'(\rho_{\rm t}-1)+2\v)\\
&\leq \chi(\rho-1)+(2\rho-1)\v=\bar\chi(\rho),
\end{aligned}
$$
and therefore $\partial^- \bar\chi(\rho_{\rm t})\leq\partial^+ \bar\chi(\rho_{\rm t})$.
\end{proof}

\begin{step}\label{stepbarchileqchi}
$\chi(\rho)+2\v a \rho+\v a^2\geq \bar\chi(\rho+a)$ for any $\rho\in[\rho_{\rm c}-1, \rho_{\rm c}]$ and any $a\in \N$.
\end{step}

\begin{proof}
It is enough to prove this for $a=1$ in regard of \eqref{eqsimbarchi}. Then, if $\rho\in[\rho_{\rm t}-1, \rho_{\rm c}]$, we actually have the equality by definition of $\bar\chi$. The case when $\rho\in[\rho_{\rm c}-1, \rho_{\rm t}-1]$ follows from the definition of $\rho_{\rm t}$ and the convexity of $\chi(\cdot)$.
\end{proof}

Now the proof of point (3) in Lemma \ref{lem-saturation} easily follows. Indeed,  from Steps \ref{stepbarchiconv} and \ref{stepbarchileqchi} and \eqref{infimumforchi}, we deduce that $\bar\chi\leq \chi$; the details of this are the same as the ones of Step \ref{non_sat_1}. Hence, the functions $\chi$ and $\bar\chi$ coincide. The observation that $\chi'(\rho_{\rm c})=\chi'(\rho_{\rm t}-1)+2\v$ is immediate using the continuity of $\chi'$.

The point (4) in Lemma \ref{lem-saturation} is straight forward.

\section{Appendix: The interacting quantum Bose gas}\label{sec-Bosegas}

\noindent For comparison to the model and the results of this paper, let us recall here the interacting Bose gas and the description of its free energy from the viewpoint of random point processes using large-deviation analysis. We are citing from \cite{ACK10}, to which we also refer for more references.

We consider an interacting bosonic many-body system in a large box in $\R^d$ at positive temperature $1/\beta\in(0,\infty)$ with fixed particle density $\rho\in(0,\infty)$ in the thermodynamic limit. Denote by
$$
\Hcal_N=-\sum_{i=1}^N \Delta_i+\sum_{1\leq i<j\leq N}v(|x_i-x_j|),\qquad x_1,\dots,x_N\in\R^d,
$$
the $N$-particle Hamilton operator with kinetic energy and pair-interaction given by an  interaction functional $v\colon [0,\infty)\to[0,\infty]$ satisfying some properties that we state later.  Since we do not want to exclude the possibility that $v$ has a singularity at $0$ satisfying $\lim_{r\downarrow 0} v(r)=\infty$, we cannot include the self-interactions; furthermore we register each pair of particles only once. We are interested in {\em bosons} and introduce a symmetrisation, i.e., we project the operator $\Hcal_N$ on the set of symmetric, i.e., permutation invariant, wave functions. Furthermore, we consider the particle system at positive temperature $1/\beta\in(0,\infty)$ in a centred box $\Lambda\subset \R^d$ with some boundary conditions, to be detailed also later. In other words, we consider the trace of the operator $\e^{-\beta\Hcal_N}$ in $\Lambda$ with symmetrisation:
$$
Z_N^{\ssup{\rm bc}}(\beta,\Lambda)={\rm Tr}_{\Lambda,+}^{\ssup{\rm bc}}(\e^{-\beta\Hcal_N}),
$$
where the index $+$ denotes the symmetrisation. This is the so-called {\em partition function} of the system, the main object of the study in this  model. We introduce the {\em particle density} $\rho\in(0,\infty)$, the number of particles per unit volume. Fix a centred box $\Lambda_N$ of volume $N/\rho$, and consider the {\em free energy},
$$
f(\beta,\rho)=-\frac 1\beta\lim_{N\to\infty}\frac 1{|\Lambda_N|}\log Z_N^{\ssup{\rm bc}}(\beta,\Lambda_N).
$$
The existence of this limit and the fact that it is independent of the boundary condition are well-known for many decades, but an explicit or even interpretable formula is still lacking, with the exception of the main result of \cite{ACK10}, which holds only for all small $\rho$, see below. 

In brevity, let us state here the main conjecture about the occurrence of {\it Bose--Einstein condensation (BEC)}: One expects that, in dimensions $d\geq 3$ but not in dimensions $d\in\{1,2\}$, the map $\rho\mapsto f(\beta,\rho)$ has a non-analyticity at some unique $\rho_{\rm c}(\beta)\in(0,\infty)$. However, much more interesting than this fact is the underlying interpretation and explanation in terms of the underlying particle process; see below.

\begin{figure}[!htpb]
	\centering
\includegraphics[scale=0.4]{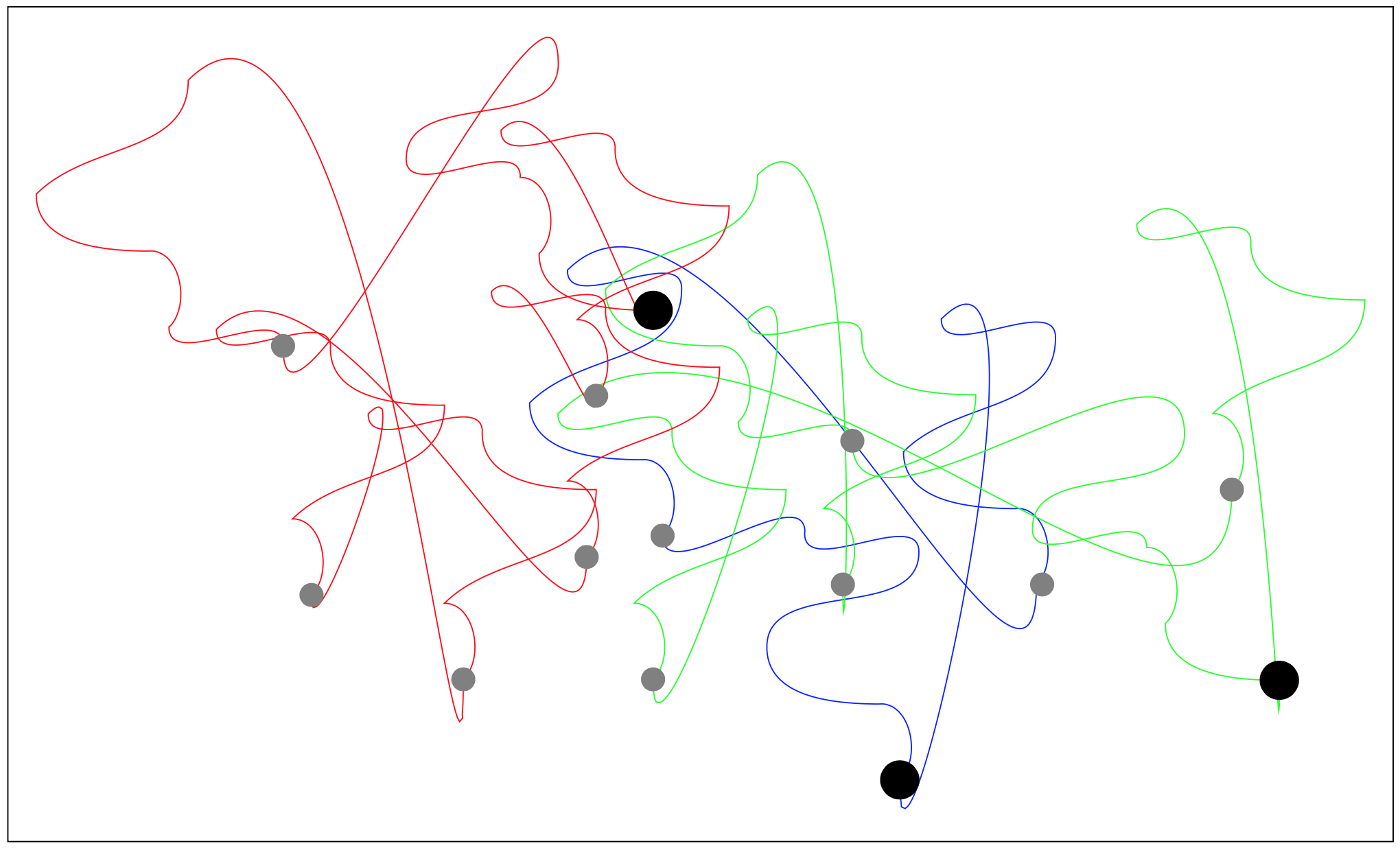}
	\caption{Illustration of a realisation of the Bose gas with 14 particles (grey and black bullets), organised in three Brownian bridges, attached to three Poisson points (black bullets). The red cycle has six particles, the blue and green ones each four.}
	\label{Pix_BoseGas}
\end{figure}

In \cite{ACK10}, a description of the model in terms of a marked PPP was developed. The marks are here random cycles of Brownian motions.  See Figure \ref{Pix_BoseGas} for an illustration of the marked point process. 

For $k\in\N$, let
$$
\omega^{\ssup k}_{ \rm P}=\sum_{x \in \xi^{\ssup k}_{ \rm P}} \delta_{(x,B_x)} \qquad\mbox{and}\qquad \omega_{\rm P}=\sum_{k\in\N}\omega^{\ssup k}_{ \rm P},
$$
then $\omega_{ \rm P}$ is the independent superposition of  PPPs $\omega^{\ssup k}_{ \rm P}$ over $k\in\N$ on $\R^d\times \Ccal_k$, where the mark space  $\Ccal_k$ is the set of continuous functions $[0,k\beta]\to\R^d$. The intensity measure $ \nu_k $ of $\omega^{\ssup k}_{ \rm P}$ is given by 
\begin{equation}\label{nudef}
\nu_k(\d x,\d f)=\frac{1}{k}\Lambdaeb(\d x)\otimes\mu_{x,x}^{\ssup{k\beta}}(\d f),
\end{equation}
where $\mu_{x,y}^{\ssup\beta}$ is the unnormalised canonical measure for a Brownian bridge from $x$ to $y$ on the time interval $[0,\beta]$; for $x=y$ it has total mass equal to $(4\pi\beta)^{-d/2}$. Alternatively, we can think of $\omega^{\ssup k}_{\rm P}$ as of an independently marked PPP  on $\R^d$, based on some standard homogeneous PPP $\xi^{\ssup k}_{\rm P}$ on $\R^d$, and a family $(B_x)_{x\in \xi^{\ssup k}_{\rm P}}$ of  i.i.d.~marks, given $\xi^{\ssup k}_{\rm P}$. The intensity of $\xi^{\ssup k}_{\rm P}$ is 
\begin{equation}\label{q*def}
q_k=\frac 1k \mu_{x,x}^{\ssup{k\beta}}(\Ccal_k)=\frac 1{(4\pi\beta)^{d/2}k^{1+d/2}}.
\end{equation}
Elements $f$ of the mark space $\Ccal_k$ have the length $\ell(f)=k$, which should also be seen as the number of particles in the mark. Indeed, a cycle $f\in\Ccal_k$ contains the $k$ particles $f(\beta), f(2\beta), f(3\beta),\dots, f(k\beta)$. 
Conditionally on  $\ell(B_x)=k$,  $B_x$ is in distribution equal to a Brownian bridge with time horizon  $[0,k\beta]$, starting and ending at $x$. Put 
\begin{equation}\label{qbardef}
\overline q=\sum_{k\in\N}q_k=(4\pi\beta)^{-d/2}\zeta(1+d/2),
\end{equation}
where $\zeta$ is the Riemann zeta function. We denote by $\Omega$ the state space of $\omega_{\rm P}$, i.e., the set of all marked point processes $\omega=\sum_{x\in\xi}\delta_{(x,f_x)}$ with point set $\xi\subset\R^d$ and marks $f_x\in\Ccal=\bigcup_{k\in\N}\Ccal_k$, starting and ending at $x$. We call the sites $f_x(k\beta)$ with $k\in\N_0$ the {\it particles} of $\omega$; each point $x\in\xi$ has precisely $\ell(f_x)$ particles.

We introduce a functional on $\Omega$ that expresses the pair interaction between  any two particles belonging to a mark in $ \Lambda\subset\R^d$.  Define the interaction between  $\Lambda$ and $\Lambda'\subset\R^d$ by
\begin{equation}\label{Hamiltonian}
\Phi_{\Lambda,\Lambda'}(\omega)=\sum_{x\in\xi\cap\Lambda, y\in\xi\cap \Lambda'}T_{x,y}(f_x,f_y),\qquad\omega\in\Omega,
\end{equation}
where we abbreviate
\begin{equation}\label{Tdef}
T_{x,y}(f_x,f_y)=\frac{1}{2}\sum_{i=1}^{\ell(f_x)}\sum_{j=1}^{\ell(f_y)} \1_{\{(x,i)\not=(y,j)\}}V(f_{x,i},f_{y,j}) \quad x,y\in\xi,\, f_x,f_y\in \Ccal,
\end{equation}
and $f_{x,i}(\cdot)=f_x((i-1)\beta+\cdot)|_{[0,\beta]}$
is the $i$-th {\em leg} of a function $f_x\in\Ccal$, and $V(f,g)=\int_0^\beta v(|f(s)-g(s)|)\,\d s$. Denote by
\begin{equation}\label{Nlength}
N^{\ssup\ell}_{\Lambda}(\omega)=\sum_{x\in\xi\cap\Lambda}\ell(f_x)
\end{equation}
the number of particles in the cloud in marks whose suspension point lies in $\Lambda$. We are going to consider three different boundary conditions in the box $\Lambda$: periodic, zero Dirichlet and open boundary condition, written \lq per\rq, \lq Dir\rq\ and \lq $\emptyset$\rq\ (the latter means that the Poisson points belong to $\Lambda$, but the particles do not have to). The first two boundary conditions are reflected in the definition of the Brownian bridges; they actually need to be adapted, which also necessitates adaptations in the intensity $q_k$ and in the mark measure $\mu_{x,x}^{\ssup{\beta k}}$. We write the superscript \lq bc\rq\ to express the boundary condition and $\P^{\ssup{\rm bc}}_\L$ and $\E^{\ssup{\rm bc}}_\L$ for the corresponding distribution and expectation of the marked PPP. 
The following is Proposition 1.1 in \cite{ACK10}.

\begin{prop}[Rewrite in terms of the marked PPP]\label{lem-rewrite} Fix $\beta\in(0,\infty)$. Let $ v\colon[0,\infty)\to (-\infty,\infty] $ be measurable and bounded from below and let $ \Lambda\subset\R^d$ be measurable with finite volume (assumed to be a torus for periodic boundary condition). Then, for any $N\in\N$, and $ {\rm bc}\in\{\emptyset,{\rm per},{\rm Dir}\} $,
\begin{equation}\label{rewrite}
\begin{aligned}
Z_N^{\ssup{\rm bc}}(\beta,\Lambda)&=\e^{|\Lambda|\overline{q}^{\ssup{\rm bc}}}
\E_\Lambda^{\ssup{\rm bc}}\big[{\rm e}^{-\Phi_{\Lambda,\Lambda}(\omega_{\rm P})}\1\{N^{\ssup\ell}_\Lambda(\omega_{\rm P})=N\}\big].
\end{aligned}
\end{equation}
\end{prop}

That is, up the non-random term $|\Lambda|\overline{q}^{\ssup{\rm bc}}$, the partition function  is equal to the expectation over the Boltzmann factor $ {\rm e}^{-\Phi_{\Lambda,\Lambda}} $ of a marked PPP restricted to a fixed total length of marks of the particles. Here we see the motivation for the box-version of the model that we introduced in Section \ref{sec-simplified}: the marks that are here random Brownian cycles are boxes in the box-version.

Now that we have revealed a characterisation of the free energy in terms of a point process with Brownian cycles as marks, we can give another, more descriptive, interpretation of BEC: for sufficiently large $\rho$, a main part of the contribution to the expectation on the right-hand side of \eqref{rewrite} should come from realisations of the point process in which a number $\asymp N$ of particles (i.e., a macroscopic part of the $N$ particles) are in long cycles, i.e., in cycles whose lengths depend on $N$ and diverge as $N\to\infty$.

Now we explain how to use large-deviation theory to derive asymptotic assertions in the thermodynamic limit, i.e., in the limit $N\to\infty$ with the box $\L=\L_N$ having volume equal to $N/\rho$. This has much to do with ergodic theory. Let $ \theta_x\colon \R^d \to \R^d$ denote the  shift operator by $x\in\R^d$; we extend it to an operator $\theta_x\colon \Omega\to\Omega$, where the shifts are performed with respect to the suspension points and their associated marks. By $\Mcal_1^{\ssup{\rm s}}(\Omega)$ we denote the set of all shift-invariant probability measures on $\Omega$; note that the distribution $\P$ of the reference process $\omega_{\rm P}$ belongs to $\Mcal_1^{\ssup{\rm s}}(\Omega)$.
We write $U=[-1/2,1/2]^d$ for the centred unit box. 

Next, we introduce an entropy term. For probability measures $\mu, \nu $ on some measurable space, we write
\begin{equation}\label{Hdef}
H(\mu|\nu)=\begin{cases} \int f\log f \, \d \nu &\mbox{if } f=\d \mu/\d \nu\mbox{ exists,} 
\\\infty & \mbox{otherwise,}\end{cases}
\end{equation}
for the relative entropy of $\mu$ with respect $\nu$. It will be clear from the context which measurable space is used. It is easy to see and well-known that $H(\mu| \nu)$ is nonnegative and that it vanishes if $\mu=\nu$.  Now we introduce  the {\em entropy density function}
\begin{equation}\label{Idef}
I(P) = \lim_{N \to \infty} \frac{1}{|\L_{N}|} H\big( P_{\L_N} \big|  \P_{\L_N}\big),\qquad P\in\Mcal_1^{\ssup{\rm s}}(\Omega),
\end{equation} 
where we write $P_{\L}$ for the projection of $P$ to $\L$, i.e., the image measure of $P$ under the map $\omega\mapsto \omega|_{\L}=\sum_{x\in\xi\cap\L}\delta_{(x, f_x)}$. According to \cite[Prop.~2.6]{GZ93}, the  limit  in \eqref{Idef}  exists, and $I$ is an affine and lower-semi-continuous function with compact level sets in the topology of local tame convergence. It  turns out there that $I$ is the rate function of a crucial large deviations principle for the family of the stationary empirical fields, which we do not write out here explicitly.

We introduce an important variational formula:
\begin{equation}\label{chi=def}
\chi(\rho)=\inf\Big\{I(P)+P(\Phi_{U,\R^d})\colon P\in\Mcal_1^{\ssup{\rm s}}(\Omega), P(N_U^{\ssup\ell})=\rho\Big\},
\end{equation}
where we write $P(f)=\int f\,\d P=\langle P,f\rangle$ for the integral of a function $f$ with respect to a measure $P$. This formula ranges over shift-invariant marked point processes $P$ and has three crucial  components: the entropic distance $I(P)$ between $P$ and the reference measure $\P$, the interaction term $P(\Phi_{U,\R^d})$ and the effective particle density per unit volume $P(N_U^{\ssup\ell})$. The main result of \cite{ACK10} that we are interested in here is the following.

\begin{theorem}[Theorem 1.2 in \cite{ACK10}]\label{thm-mainresACK} Let $v\colon[0,\infty)\to[0,\infty]$ be measurable such that $v(r)\leq A r^{-h}$ for some $A\in(0,\infty)$ and some $h\in(d,\infty)$ and all sufficiently large $r$, and assume that $\liminf_{t\downarrow0}v(r)>0$ and that $\alpha(v)=\int_{\R^d}v(|x|)\,\d x$ is finite.  Then, for any $\beta,\rho\in(0,\infty)$ such that $(4\pi\beta)^{-d/2}<\rho\e^{\beta \rho \alpha(v)}$, and for any boundary condition $ {\rm bc}\in\{\emptyset,{\rm Dir}, {\rm per}\} $, 
\begin{equation} \label{freeenlower}
\liminf_{N \to \infty} \frac{1}{|\L_N|} \log Z^{\ssup{\rm bc}}_N(\beta,\L_{N})\geq  \overline q-\chi(\rho),
\end{equation}
\end{theorem}

\begin{figure}[!htpb]
\centering
\includegraphics[scale=0.4]{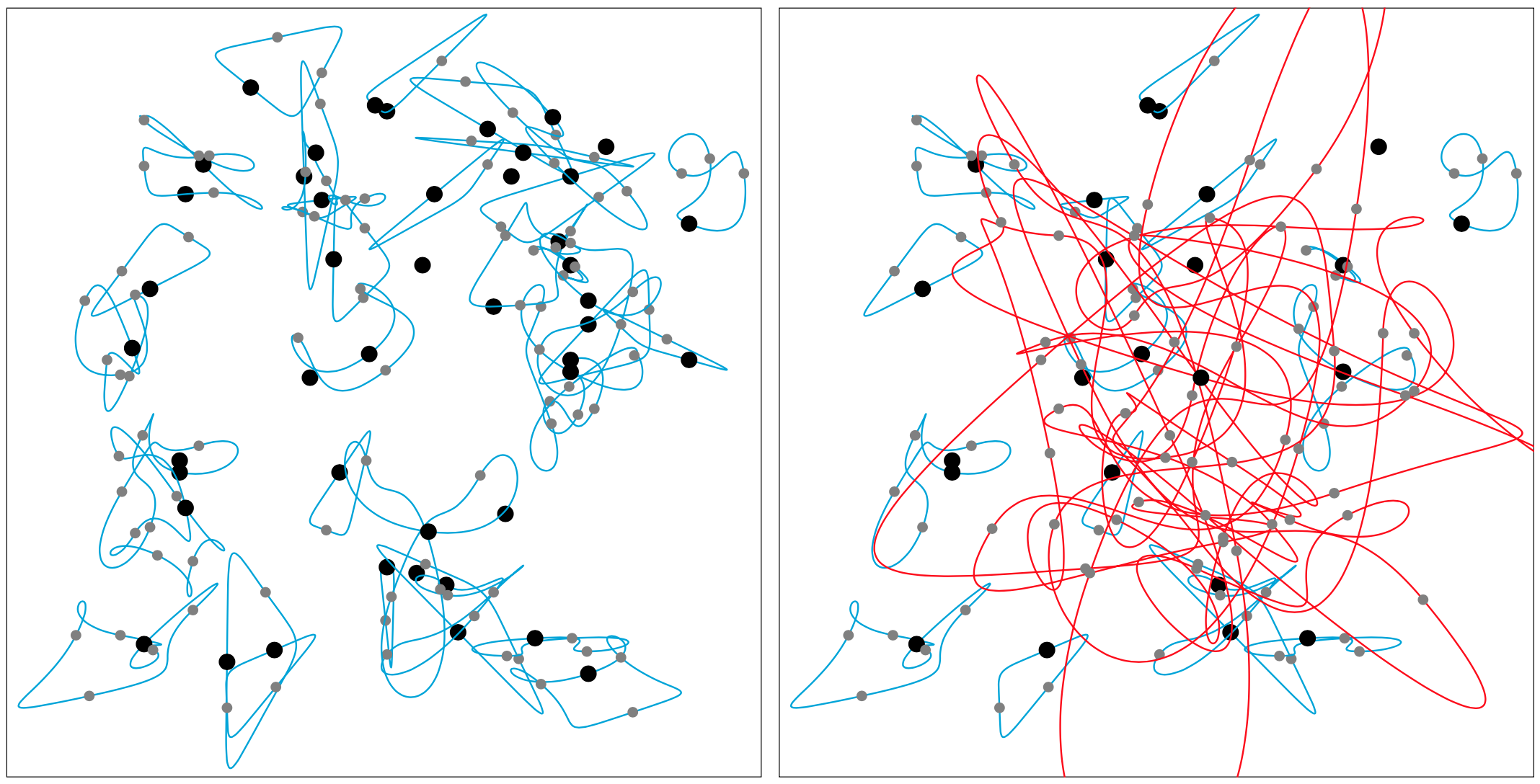}
\caption{Illustration of a subcritical (low $\rho$) Bose gas without condensate (left) and supercritical (large $\rho$) Bose gas with additional condensate (red) (right)}
\label{Pix_BoseGassub}
\end{figure}

This is only a small part of what we are able to prove for the box-version in Section \ref{sec-simplified}. We conjecture that a great deal of that results are true also here. In particular, we conjecture that  \eqref{freeenlower} and the complementary inequality $\leq$ hold true for any $\rho\in(0,\infty)$.  Furthermore, we believe that BEC can be characterised  in terms of existence of minimisers, analogously to the box-version. See Figure~\ref{Pix_BoseGassub} for illustrations of the two phases in terms of Brownian cycle ensembles. 

The main difficulty in the proof of $\leq$ in \eqref{freeenlower} is the discontinuity of the functional $P\mapsto P(N_U^{\ssup\ell})$; it is only semi-continuous from below. This could be partially overcome by another main result of \cite{ACK10}: for sufficiently small $\rho$, the corresponding upper bound holds with $\chi(\rho)$ replaced by the same formula with the condition $P(N_U^{\ssup\ell})=\rho$ replaced by $P(N_U^{\ssup\ell})\leq \rho$; however, this formula is not expected to be accurate at all for large $\rho$. Furthermore, it was proved in \cite{ACK10} that \eqref{freeenlower} holds with this formula instead of $\chi(\rho)$. It is not deep to show that both formulas coincide for small $\rho$, so the main message from \cite{ACK10} is that the free energy is expressed in terms of the natural variational formula $\chi(\rho)$, if $\rho$ is sufficiently small.

Let us briefly discuss BEC in the {\em free Bose gas}, where no interaction is present, i.e., $v=0$. In this case, the formula in \eqref{chi=def} drastically reduces to the formula
\begin{equation}
\chi^{\ssup{v=0}}(\rho)=\inf\Big\{H(m|q)\colon m\in [0,\infty)^{\N}, \sum_{k\in\N}k m_k=\rho\Big\},
\end{equation}
where $H(m|q)=\sum_k(q_k-m_k+m_k\log \frac {m_k}{q_k})$ is the relative entropy of the sequence $m$ with respect to $q$. Compare to the  remark at the end of Section \ref{sec-results}  on  $\chi^{\ssup{v=0}}(\rho)$ and the phase transition that it undergoes as a function of $\rho$.
Indeed, a possible minimiser $m$ is characterised by the Euler--Lagrange equation $m_k=q_k\e^{\alpha k}$ for $k\in\N$, where $\alpha\in\R$ is the Lagrange multiplier. In order to meet the constraint $\sum_k km_k=\rho$, certainly $\alpha$ needs to be non-positive, and the largest value $\rho$ that can be achieved by this is $\rho_{\rm c}(\beta)=\sum_k k q_k=(4\pi\beta)^{-d/2}\zeta(d/2)$ (compare to \eqref{qbardef}), which is finite precisely in $d\geq 3$.  Hence, this is the critical value for existence of a minimiser, i.e., the critical threshold for the occurrence of BEC. The understanding is that, for $\rho>\rho_{\rm c}(\beta)$, it is not possible to arrange all the microscopic particles in finite-size cycles, and all the remainder is the condensate. However, it gives no mathematical expression for the condensate.

\end{document}

%% file: Pix_rho.tex
\begin{tikzpicture}[scale = 1.25]
		\draw[thick, ->] (6,0) -- (10,0);
		\draw[thick, ->] (0,5) -- (0,10);
		\draw[thick, dashed] (0,4) -- (0,5);
		\draw[thick, red, dashed] (4,0) -- (5,0);
		\draw[thick, red] (0,0) -- (4,0);
		\draw[thick, red] (5,0) -- (6,0);
		\draw[thick] (0,0) -- (0,4);
		\draw[thick, red] (6.7,1) -- (7,1);
		\draw[thick, red] (7.7,2) -- (8,2);
		\draw[thick, red] (8.7,3) -- (9,3);
		\draw[thick, red] (9.7,4) -- (10,4);
		\draw[thick, red] (6,0) -- (6.7,1);
		\draw[thick, red] (7,1) -- (7.7,2);
		\draw[thick, red] (8,2) --( 8.7,3);
		\draw[thick, red] (9,3) --( 9.7,4);
		\draw[thick, green] (6,6) -- (10,10);
		\draw[thick, blue] (0,0) -- (4,4);
		\draw[thick, blue, dashed] (4,4) -- (5,5);
		\draw[thick, blue] (5,5) -- (6,6);
		
		\coordinate[label = 0 : ${\left(\rho_{\rm mi}, \rho_{\rm ma}\right)}$](a) at  (0,10);
		\coordinate[label = 0 : $\rho$](a) at  (10,0);
		\coordinate[label = 270 : \text{$\rho_{\rm c}$} ](a) at  (6,0);
		\coordinate[label = 90 : \text{$\rho_{t}$} ](a) at  (6.7,0);
		\coordinate[label = 270 : \text{$\rho_{\rm c}\hspace{-0.1cm}+\hspace{-0.05cm}1$} ](a) at  (7,0);
		\coordinate[label = 270 : \text{$\rho_{\rm c}\hspace{-0.1cm}+\hspace{-0.05cm}2$} ](a) at  (8,0);
		\coordinate[label = 270 : \text{$\rho_{\rm c}\hspace{-0.1cm}+\hspace{-0.05cm}3$} ](a) at  (9,0);
		\coordinate[label = 0 : \text{$\rho_{\rm c}$} ](a) at  (-0.7,6);
		\draw[thick] (-0.1,6) -- (0.1,6);
		\coordinate[label = 0 : \text{$1$} ](a) at  (-0.5,1);
		\coordinate[label = 0 : \text{$2$} ](a) at  (-0.5,2);
		\coordinate[label = 0 : \text{$3$} ](a) at  (-0.5,3);
		\coordinate[label = 0 : \text{$4$} ](a) at  (-0.5,4);
		\coordinate[label = 0 : $\rho_{\rm mi}+\rho_{\rm ma}$](a) at  (10,10);
		\coordinate[label = 0 : $\rho_{\rm ma}$](a) at  (10,4);
		\coordinate[label = 0 : $\rho_{\rm mi}$](a) at  (10,6);
		\foreach \i in {1,...,4}
{
		\draw[thick] (\i+5,0.1) -- (\i+5,-0.1);
		\draw[thick] (-0.1, \i) -- (0.1,\i);
}
		\draw[thick] (6.7,0.1) -- (6.7,-0.1);
		
		\foreach \i in {1,...,4}
{
		\draw[thick, blue] (\i+5,6) -- (\i+5.7,5.7);
		\draw[thick, blue] (\i+5.7,5.7) -- (\i+6,6);
}

%
\end{tikzpicture}

%% file: Pix_LowerBound.tex
\begin{tikzpicture}[decoration={brace}, scale=2.32]
\draw [fill=orange] (1,1) rectangle (4,4);
\draw [fill=red] (2,2) rectangle (3,3);
\draw[step=1mm,help lines] (0,0) grid (50mm,50mm);
\draw[step=10mm] (0,0) grid (50mm,50mm);
\draw [decorate,decoration={brace,amplitude=4pt}]
(5,1)--(5,0) node[midway, right, font=\footnotesize, xshift=2pt] {$2R$};
\draw [decorate,decoration={brace,amplitude=4pt}]
(5.5,5)--(5.5,0) node[midway, right, font=\footnotesize, xshift=2pt] {$|\Lambda_N|^{1/d}$};
\filldraw[black] (2.5,2.5) circle (1pt);
\end{tikzpicture}

%% file: BoxBoseGas_ArXiv.bbl
\begin{thebibliography}{WWW98}

\bibitem[ACK11]{ACK10}
{\sc S.~Adams, A.~Collevecchio} and  {\sc W.~K\"onig},
\newblock A variational formula for the free energy of an interacting many-particle system,
\newblock {\it Ann.~Probab.} {\bf 39:2}, 683--728 (2011).

\medskip

\bibitem[AFY19]{AFY19}
\newblock {\sc I.~Armend\'ariz, P.A.~Ferrari,} and {\sc S.~Yuhjtman,}
\newblock Gaussian random permutation and the boson point process,
\newblock {\it Comm.~Math.~Phys.} {\bf 387:3}, 1515--1547 (2021).


\medskip

\bibitem[BKM21]{BKM21}
\newblock {\sc E.~Bolthausen,  W.~König}, and {\sc Ch.~Mukherjee},
\newblock Bose--Einstein condensate and the self-avoiding walk,
\newblock {\em in preparation} (2022).

\medskip

\bibitem[Fe53]{F53} 
\newblock {\sc R.P.~Feynman},
\newblock Atomic theory of the $ \lambda $ transition in Helium,
\newblock {\it Phys. Rev.} {\bf 91}, 1291--1301 (1953).

\medskip

\bibitem[F91]{F91}
\newblock {\sc K.-H.~Fichtner},
\newblock On the position distribution of the ideal Bose gas,
\newblock {\it Math. Nachr.} {\bf 151}, 59--67 (1991). 

\medskip

\bibitem[FKSS20]{FKSS20}
\newblock {\sc J.~Fröhlich, A.~Knowles, B.~Schlein,} and {\sc V.~Sohinger,}
\newblock  A path-integral analysis of interacting Bose gases and loop gases. 
\newblock {\em J. Stat. Phys.} {\bf 180:1-6}, 810–831 (2020).
 
\medskip
 
\bibitem[G88]{G88}
\newblock {\sc H.-O.~Georgii},
\newblock {\it Gibbs Measures and Phase Transitions},
\newblock Berlin: de Gruyter (2011).

\medskip

\bibitem[G93]{G93}
\newblock {\sc H.-O.~Georgii,}
\newblock Large deviations and maximum entropy principle for interacting random fields on $\Z^d$,
\newblock {\it Ann.~Probab.} {\bf 21:4},  1845-1875 (1993).

\medskip

\bibitem[GZ93]{GZ93}
\newblock {\sc H.-O.~Georgii} and {\sc H.~Zessin},
\newblock Large deviations and the maximum entropy principle for marked point random fields,
\newblock {\it Prob.~Theory Relat. Fields} {\bf 96}, 177--204 (1993).

\medskip

\bibitem[G94]{G94}
\newblock {\sc H.-O.~Georgii,}
\newblock Large deviations and the equivalence of ensembles for Gibbsian particle systems with superstable interaction,
\newblock {\it Prob.~Theory Relat. Fields} {\bf 99}, 171--195 (1994).

\medskip

\bibitem[G70]{G70}
\newblock {\sc J.~Ginibre,}
\newblock {\it Some Applications of Functional Integration in Statistical Mechanics, and Field Theory\/}, 
\newblock C. de Witt and R. Storaeds, Gordon and Breach, New York   (1970).

\medskip

\bibitem[J20]{J20}
\newblock {\sc S.\ Jansen,}
\newblock Thermodynamics of a hierarchical mixture of cubes. 
\newblock {\it J. Stat. Phys.} {\bf 179:2}, 309-340  (2020).

\medskip

\bibitem[NPZ13]{NPZ13}
\newblock {\sc B.~Nehring, S.~Poghosyan,} and {\sc H.~Zessin,}
\newblock On the construction of point processes in statistical mechanics. 
\newblock {\em J. Math. Phys.} {\bf 54:6}, 063302 (2013).

\medskip

\bibitem[RZ20]{RZ20}
\newblock {\sc S.~R\oe lly} and {\sc A.~Zass,}
\newblock Marked Gibbs point processes with unbounded interaction: an existence result,
\newblock  {\em J. Stat. Phys.} {\bf 179:4}, 972–996 (2020). 

\medskip

\bibitem[TI06]{TI06}
\newblock {\sc H.~Tamura} and {\sc K.R.A.~Ito,}
\newblock A canonical ensemble approach to the fermion/boson random point processes and its applications. 
\newblock{\em Comm. Math. Phys.} {\bf 263:2}, 353-380 (2006). 

\medskip

\bibitem[U06a]{U06a}
\newblock{\sc D.~Ueltschi,}
\newblock Feynman cycles in the Bose gas. 
\newblock {\em J. Math. Phys.} {\bf 47:12}, 123303, 15 pp., (2006).

\medskip

\bibitem[V21]{V20}
\newblock {\sc Q.~Vogel},
\newblock Emergence of interlacements from the finite volume Bose soup,
\newblock preprint, arXiv:2011.02760 (2021).


\end{thebibliography}
